\documentclass[12pt]{amsart}
\usepackage{amssymb,amsmath,amsfonts,amscd,euscript}
\usepackage{tikz}
\usepackage{tikz-cd}
\usepackage{fullpage}
\usepackage{stmaryrd}
\usepackage[all]{xy}
\usepackage{color}

\usepackage[backref=page]{hyperref}
\usepackage{todonotes}

\usepackage{amsthm,enumitem,bm,pict2e,xcolor}
\usepackage{indentfirst}

\newcommand{\nc}{\newcommand}

\numberwithin{equation}{section}
\newtheorem{thm}{Theorem}[section]
\newtheorem*{thm*}{Theorem}

\newtheorem{prop}[thm]{Proposition}
\newtheorem{lem}[thm]{Lemma}
\newtheorem{cor}[thm]{Corollary}

\newtheorem{dfn}[thm]{Definition}

\theoremstyle{definition}
\newtheorem{rem}[thm]{Remark}
\newtheorem{example}[thm]{Example}

\nc{\gl}{\mathfrak{gl}}
\nc{\GL}{\mathsf{GL}}
\nc{\g}{\mathfrak{g}}
\nc{\gh}{\widehat\g}
\nc{\h}{\mathfrak{h}}
\nc{\la}{\lambda}
\nc{\al}{\alpha }
\nc{\be}{\beta }
\nc{\ve}{\varepsilon }
\nc{\om}{\omega }
\nc{\br}{{\bf r}}
\nc{\tr}{{\rm tc}}
\nc{\sh}{{\rm sh}}
\nc{\ta}{\theta}
\nc{\ch}{{\mathop {\rm ch}}}
\nc{\Tr}{{\mathop {\rm Tr}\,}}
\nc{\Id}{{\mathop {\rm Id}}}
\nc{\ad}{{\mathop {\rm ad}}}
\nc{\bra}{\langle}
\nc{\ket}{\rangle}
\nc{\pa}{\partial}
\nc{\ld}{\ldots}
\nc{\cd}{\cdots}
\nc{\hk}{\hookrightarrow}
\nc{\T}{\otimes}
\nc{\gr}{\mathrm{gr}}
\nc{\ov}{\overline}

\nc{\cO}{\mathcal O}
\nc{\msl}{\mathfrak{sl}}
\nc{\mgl}{\mathfrak{gl}}
\nc{\U}{\mathrm U}
\nc{\V}{\EuScript V}
\nc{\cL}{\mathcal{L}}
\nc{\La}{\Lambda}
\newcommand{\im}{\mathrm{im }}
\nc{\fQ}{\mathfrak{Q}}
\newcommand{\Ga}{\Gamma}
\nc{\cA}{\mathcal{A}}
\nc{\bc}{{\bf c}}

\newcommand{\bC}{{\Bbbk}}

\newcommand{\bZ}{{\mathbb Z}}
\newcommand{\bN}{{\mathbb N}}
\newcommand{\bP}{{\mathbb P}}

\newcommand{\fh}{{\mathfrak h}}

\newcommand{\fg}{{\mathfrak g}}

\newcommand{\Hom}{\mathrm{Hom}}

\newcommand{\mL}{\mathcal L}

\newcommand{\mJ}{\mathcal J}

\newcommand{\mC}{\mathcal C}
\nc{\mfd}{\mathfrak d}
\nc{\red}{\mathrm{red}}

\nc{\cat}{\mathcal{C}}
\nc{\RepGL}[2]{{\mathsf{Rep}}(\GL_{#1})^{(#2)}}
\nc{\RepGLt}[2]{{\mathsf{Rep}}(\mgl_{#1}[t])^{(#2)}}
\nc{\RepSL}[2]{{\mathsf{Rep}}(\SL_{#1})^{(#2)}}
\nc{\RepSLt}[2]{{\mathsf{Rep}}(\msl_{#1}[t])^{(#2)}}
\nc{\Rep}{\mathsf{Rep}}
\nc{\de}{{\text -}}
\nc{\id}{\mathrm{id}}

\DeclareMathOperator{\grad}{grad}
\DeclareMathOperator{\Spec}{Spec}
\DeclareMathOperator{\Proj}{Proj}
\DeclareMathOperator{\Der}{Der}

\DeclareMathOperator{\spanned}{span}
\DeclareMathOperator{\In}{in}
\DeclareMathOperator{\pr}{pr}

\mathcode`l="8000
\begingroup
\makeatletter
\lccode`\~=`\l
\DeclareMathSymbol{\lsb@l}{\mathalpha}{letters}{`l}
\lowercase{\gdef~{\ifnum\the\mathgroup=\m@ne \ell \else \lsb@l \fi}}%
\endgroup

\DeclareFontFamily{U}{mathx}{\hyphenchar\font45}
\DeclareFontShape{U}{mathx}{m}{n}{
      <5> <6> <7> <8> <9> <10>
      <10.95> <12> <14.4> <17.28> <20.74> <24.88>
      mathx10
      }{}
\DeclareSymbolFont{mathx}{U}{mathx}{m}{n}
\DeclareMathSymbol{\bigtimes}{1}{mathx}{"91}

\begin{document}

\title[On reduced arc spaces of toric varieties]
{On reduced arc spaces of toric varieties}

\author{Ilya Dumanski}
\address{I. Dumanski: \newline
HSE University\\ Faculty of Mathematics\\ Ulitsa Usacheva 6\\Moscow 119048\\Russia}
\email{ilyadumnsk@gmail.com}
\author{Evgeny Feigin}
\address{E. Feigin: \newline
School of Mathematical Sciences, Tel Aviv University, Tel Aviv, 69978, Israel}
\email{evgfeig@gmail.com}

\author{Ievgen Makedonskyi}
\address{Ie. Makedonskyi:\newline
 Friedrich Schiller Jena University, Faculty of Mathematics and Computer Science; 
Ernst-Abbe-Platz 2, 07743 Jena, Germany and Beijing Institute of Mathematical Sciences and Applications}
\email{makedonskii\_e@mail.ru}

\author{Igor Makhlin}
\address{I. Makhlin:\newline
Technische Universit\"at Berlin, Berlin, Germany}
\email{iymakhlin@gmail.com}

\begin{abstract}
 An arc space of an affine cone over a projective toric variety
 is known to be non-reduced in general. 
 It was demonstrated recently that the reduced scheme structure 
 of arc spaces 
 is very meaningful from algebro-geometric, representation-theoretic and combinatorial points of view.   
 In this paper we develop a general machinery for the description of the reduced arc spaces of affine cones over toric varieties. We apply our techniques to a number of classical cases and explore some connections with representation theory of current algebras. 
\end{abstract}

\maketitle
\tableofcontents

\section*{Introduction}
Let $V$ be the affine cone over a  normal projective  toric variety.  
In this paper we study the arc space over $V$. Roughly speaking, the arc space consists of $\Bbbk[[t]]$-points of $V$  \cite{Na,EM,Fr,Mus2}, where $\Bbbk$ is the base field.
If $V$
is embedded into an affine space with coordinates $X_1,\dots,X_n$, then  each variable $X_i$ produces a sequences of variables 
$X_i^{(j)}$ with generating series $X_i(s)=\sum_{j\ge 0} X_i^{(j)}s^j$. The relations cutting out $V$ inside the affine space are replaced with an infinite
number of relations appearing as coefficients of powers of $s$. 
For example, the relation $X_1X_3-X_2^2$ defining the degree two Veronese curve 
gives $X_1(s)X_3(s)-X_2^2(s)$ and one is interested in the ideal, generated by the $s$-coefficients of this relation. 
We denote the arc space by $J^\infty(V)$. Arc spaces  
attracted a lot of attention in the last decades (see e.g. \cite{AS, AM, BS2, dFD, G, Ish1})  due to the beauty of the theory and
various connections with 
singularity theory, theory of motivic integration, differential algebra, representation theory, vertex algebras and other mathematical subjects. The particular case of arc spaces of toric varieties was also studied by several authors. For instance, in \cite{KI} the Nash problem is positively solved for toric varieties; in \cite{Ish2} the orbits of the $J^\infty((\Bbbk^\times)^n)$-action are studied; in \cite{BS1} the finite formal model for toric singularities is constructed; for other results see also \cite{R, Mou1}.

One interesting feature of arc spaces is that $J^\infty(V)$ may be non-reduced even
if $V$ is. The simplest example relevant to our studies is the cone over the degree three Veronese curve. It is defined inside ${\mathbb A}^4$ by three quadratic equations on variables 
$X_0,X_1,X_2,X_3$, so the arc space
$J^\infty(V)$ is defined by three series of relations. However, in order to get the reduced
scheme structure of $J^\infty(V)$ one has to add another series of quadratic relations.    

Studying the reduced scheme $J_\red^\infty(V)$ is an important problem, related to singularity theory \cite{Mus1}, differential algebra \cite{BLOP}, representation theory \cite{FM2, DF}, invariant theory \cite{LS}, vertex algebras \cite{FM1}, and other areas.
Our main motivation for the study of the reduced schemes comes from geometric representation theory due to the connection with the theory of semi-infinite flag varieties (see   \cite{DF,DFF,FM2,Mak}).
In particular, in \cite{FM2, Mak, DF} the complete list of relations is found for the arc space of certain varieties (proving, in particular, that the radical is differentially finitely generated).
Let us also mention several other results.
In \cite{BLOP} an algorithmic approach to the problem of finding the radical of a differential ideal is developed; in \cite{Se, KS} the question of reducedness of the arc space of a plane curve is investigated; in \cite{FM1} a vertex algebra structure on a reduced arc-ring for flag varieties is introduced; see also \cite{BH, CS1, CS2} for other results.

Let $V={\rm Spec}(R)$.
We denote by $J^\infty(R)$ the ring of functions on $J^\infty(V)$ and by $J_\red^\infty(R)$ the ring of functions on $J_\red^\infty(V)$. 
The goal of this paper is two-fold. First, we develop a technique for the study of reduced arc rings for toric varieties in terms of certain spaces of symmetric polynomials.
Second, we apply our construction to the description of $J_\red^\infty(R)$
for several classical examples. Let us describe our results in more details. 

We note that since $V$ is an affine cone over a projective variety, there are two
natural gradings on $J_\red^\infty(R)$: the $z$-grading defined by attaching degree $1$ to 
each variable $X_i^{(j)}$ and the $q$-grading attaching degree $j$ to $X_i^{(j)}$.  These gradings  define a graded character of the arc rings, which are series in $z$ and $q$. Note that this graded character of the rings $J^\infty(R)$ and $J_\red^\infty(R)$ is of great importance in singularity theory and combinatorics, see \cite{MP, BGK, Mou2, Mou3, BMS}.
The two main problems we address in this paper are as follows:
\begin{itemize}
    \item describe $J_\red^\infty(R)$ in terms of generators and relations,
    \item find the graded character of $J_\red^\infty(R)$.
\end{itemize} 

These problems are resolved for affine cones over certain (not necessarily toric) projective varieties in \cite{FM2,Mak,SS}. For instance, for Veronese embeddings of the projective line the $z$-homogeneous components of the 
reduced arc rings are isomorphic to the dual global Demazure modules (see \cite{DF}).  More generally, in favourable situations the following nice properties hold true:
\begin{itemize}
    \item the ideal of relations is generated in degree two by explicitly given relations, 
    \item for each $L\ge 0$ the $q$-character (the Hilbert series with respect to the grading by $q$-degree) of the $z$-degree $L$ homogeneous component of the arc ring
    admits the factorization $\chi_L(q)/(q)_L$ for certain polynomial $\chi_L(q)$
\end{itemize}
(here $(q)_L=(1-q)\dots (1-q^L)$). We note that the factor $1/(q)_L$ is the $q$-character of 
space of symmetric functions in $L$ variables. In a favourable situation the character has the above form because
\begin{itemize}
    \item the $L$-th $z$-homogeneous component of the arc ring admits a free action of the ring of symmetric polynomials in $L$ variables.  
\end{itemize}

For a normal lattice polytope $P\subset {\mathbb R}^n$ we denote by $X(P)$ the corresponding  projective toric variety equipped with a projective embedding (see Subsection \ref{subsection toric varieties} for definitions and constructions). The variety
$X(P)$ is the projective spectrum of the toric ring $R(P)$.
The ring $R(P)$ is known to be generated inside a polynomial ring in $z_1,\dots,z_n$ and $w$ 
by the elements $Y_{\bar\al^1}, \dots,Y_{\bar\al^m}$, where $Y_{\bar\al^i}$ are the exponents in $z$ variables of the integer points $\bar\al^i$ of $P$ multiplied by an extra variable $w$. One also has a presentation of $R(P)$ as the quotient of the
polynomial ring in variables $X_{\bar\al^1},\dots,X_{\bar\al^m}$ corresponding to the elements $Y_{\bar\al^i}$.

Similarly to the finite picture, 
the reduced arc ring $J_\red^\infty(R(P))$ we are interested in is generated by the coefficients of the currents $Y_{\bar\al^i}(s)$. Thus, the reduced arc ring
has two realizations: as a subring of the polynomial ring in variables $z_i^{(j)}, w^{(j)}$ and as a quotient ring of the polynomial ring in variables $X^{(j)}_{\bar\al^i}$. 

We introduce a  family of subspaces $A_{\bar r} \subset  J_\red^\infty(R(P))$ labeled by tuples ${\bar r}=(r_1,\dots,r_m)$ of nonnegative integers. The space $A_{\bar r}$ is realized inside 
the polynomial ring in variables $z_i^{(j)}$; it is spanned as a vector space by elements of the form
\[
Y_{\bar \alpha^1}^{(j_1^1)}\dots Y_{\bar \alpha^1}^{(j_{r_1}^1)} \dots
Y_{\bar \alpha^m}^{(j_1^m)}\dots Y_{\bar \alpha^m}^{(j_{r_m}^m)}
\]
for all $j_u^v \in \bN$. 
Since the ring $J_\red^\infty(R(P))$ is spanned by such elements, it can be presented as the sum of subspaces $A_{\bar r}$. This sum is
clearly not direct. In order to control the ring  $J_\red^\infty(R(P))$ we introduce an order $\preceq$ on the set of $m$-tuples ${\bar r}$ and describe the  subquotients
\begin{equation}\label{subq}
\sum_{\bar r' \preceq \bar r}A_{\bar r'}\left/\sum_{\bar r' \prec \bar r}A_{\bar r'}\right..
\end{equation}
in terms of certain symmetric polynomials. More precisely, let   $\Lambda_{\bar r}(\mathbf{t})=\Bbbk[t_i^{(j)}]_{i=1,\dots,m;j=1,\dots,r_i}^{\bigtimes\mathfrak{S}_{r_i}}$ be the ring of partially symmetric polynomials, see \eqref{LambdaDefinition} and discussions after it for the precise definition. We construct  embeddings
\[
\left(\sum_{\bar r' \preceq \bar r}A_{\bar r'}\left/\sum_{\bar r' \prec \bar r}A_{\bar r'}\right.\right)^* \subset \Lambda_{\bar r}(\mathbf{t})
\]
(here we mean the restricted dual with respect to the natural grading, see Section \ref{SymPol}).
In Section \ref{SymPol} (see Proposition \ref{DualSpaceDerivedQuagraticMon}) we prove
\begin{thm}
There is a homogeneous element $K \in \Lambda_{\bar r}(\mathbf{t})$ such that 
\[\left(\sum_{\bar r' \preceq \bar r}A_{\bar r'}\left/\sum_{\bar r' \prec \bar r}A_{\bar r'}\right.\right)^* \subset K\Lambda_{\bar r}(\mathbf{t}).\]
\end{thm}
In Section \ref{ToricArcs} (see Corollary \ref{DualToSubauotient} and Lemma \ref{Split}) we prove
\begin{thm}
There is a homogeneous element $\Gamma \in \Lambda_{\bar r}(\mathbf{t})$  such that 
\[\left(\sum_{\bar r' \preceq \bar r}A_{\bar r'}\left/\sum_{\bar r' \prec \bar r}A_{\bar r'}\right.\right)^* \supset \Gamma \Lambda_{\bar r}(\mathbf{t}).\]
\end{thm}

We also prove the following corollary:

\begin{cor}
If $K=\Gamma$, then for all $L\ge 0$ the degree $L$ $z$-homogeneous component of the reduced arc ring admits a free action of the algebra of symmetric polynomials in $L$ variables. 
\end{cor}

The proofs are constructive, i.e.\ we compute  elements $K$ and $\Gamma$. In the rest of the paper we provide several examples when $K=\Gamma$, and we call such polytopes and corresponding varieties favourable. In particular, starting from a one-dimensional case (with $P$ being a segment and $X(P)$ isomorphic to a projective line) we develop an inductive procedure providing a family of favourable polytopes.
Namely, for a convex function $\zeta: P\rightarrow \mathbb R_{\ge0}$ with $\zeta(\bar \alpha^i)=\zeta_i\in\bZ$, set $\zeta_{\max}=\max\{\zeta_i\}$.
Denote by $P^{\zeta}$ the polytope in $\mathbb{R}^{n+1}$ obtained as the convex hull of $P\times \zeta_{\max}$ and all $\bar \alpha^i\times(\zeta_{\max}-\zeta_i)$.

We prove the following theorem (see Theorem \ref{main} for a precise statement). 
\begin{thm}
Let $P$ be a $d$-dimensional lattice polytope such that $X(P)$ is favourable. Then 
$X(P^{\zeta})$ is also favourable provided certain conditions are satisfied.
\end{thm}

As a corollary we derive the following.
\begin{cor}
The toric varieties $X(P)$ are favourable for $P$ being a simplex, a parallelepiped, a Hirzebruch trapezoid. In all these cases the radical of the ideal of the corresponding arc scheme is finitely generated as a differential ideal.
\end{cor}

In Appendix~\ref{appendix} we also develop a completely different technique to study the case when $P$ is the product of simplices. This technique uses the representation theory of current algebras and the geometry of the semi-infinite flag variety.

The paper is organized as follows. In Section \ref{Prelim} we collect the main definitions and basic properties of the main objects of study: toric varieties, arc schemes, arc rings, supersymmetric polynomials. In Section \ref{SymPol} we develop a general
machinery for the description of arc rings in terms of certain subspaces
of symmetric polynomials. The approach is applied to the case of toric varieties in Section \ref{ToricArcs}. In Section \ref{Inductive} we present an inductive procedure for constructing a large family of favourable toric arc spaces. 
The construction is applied in Section \ref{Examples} where we work out explicitly several
classical examples. Finally, in Appendix~\ref{appendix} we provide a representation theoretic description
of the arc rings of the Veronese--Segre embeddings.   

{\it Acknowledgements}. The authors are indebted to Alexander Popkovich for useful discussions and his ability to keep calm in any situation. Ie.\ Makedonskyi is grateful to the Max Planck Institute for Mathematics for help in leaving Russia. The work was partially supported by the Basic Research Program at the HSE University.

\section{Preliminaries and generalities}\label{Prelim}
We work over an algebraically closed field $\Bbbk$ of characteristic zero. For variables $z_1, \dots, z_m$ and vector
$\bar a =(a_1, \dots, a_m)\in \mathbb{Z}^m$ we use the multi-index notation:
\[z^{\bar a}=z_1^{a_1}\dots z_m^{a_m}.\]
\subsection{Toric varieties} \label{subsection toric varieties}
The basic properties of toric varieties recalled here can be found in~\cite{CLS}.

\begin{dfn}
A toric variety $V$ is an algebraic variety $V$ containing the algebraic torus $(\Bbbk^\times)^k$ as an open dense subset such that the natural action of the torus on itself can be extended to an action on the variety.
\end{dfn}

\begin{dfn}
An ideal in a polynomial ring is called {\it toric} if it is prime and generated by binomials. A ring is toric if it is the quotient of a polynomial ring by a toric ideal. 
\end{dfn}
\begin{prop}\label{ToricIdeal}
Affine toric varieties are precisely those varieties which have the form $\Spec R$ for a toric ring $R$.
\end{prop}

For a rational cone $\sigma \subset  \mathbb{R}^m$ its set of integer points $\sigma \cap \mathbb{Z}^m$ forms a semigroup. 
\begin{prop}
Normal affine toric varieties are precisely those varieties which have the form
\[\Spec(\Bbbk[\sigma \cap \mathbb{Z}^m])\] for some $m$ and rational cone $\sigma \subset  \mathbb{R}^m$ of dimension $m$.
\end{prop}

A convex polytope $P \subset \mathbb{R}^m$ is said to be normal if for any integer $k>0$ every integer point in its dilation $kP$ can be expressed as the sum of $k$ (not necessarily distinct) integer points in $P$. In other words, this means that the set $kP\cap\mathbb{Z}^m$ is the $k$-fold Minkowski sum of $P\cap\mathbb{Z}^m$ with itself.

Consider a normal convex polytope $P \subset \mathbb{R}^m$.
Define a polyhedral cone $C(P) \subset \mathbb{R}^{m+1}$ to be generated by vectors
$(a,1) \in \mathbb{R}^{m+1}, a \in P$. The semigroup of this cone $C(P) \cap \mathbb{Z}^{m+1}$ consists of pairs $(b,k)$ where $b$ is an integer point in $kP$.
By normality this semigroup is generated by the points
$(a,1)$, $a \in P \cap \mathbb{Z}^{m}$.
The ring $R(P):=\Bbbk[C(P)\cap \mathbb{Z}^{m+1}]$ is naturally graded, the $k$th graded component is equal to the linear span of $z^{\bar{a}} w^k$, $\bar a \in kP \cap  \mathbb{Z}^{m}$.

\begin{dfn}
The projective toric variety of $P$ is $\Proj(R(P))$.
\end{dfn}

We denote the generators $z^{\bar a}w$ of the ring $R(P)$ by
$Y_{\bar a}$, $\bar a \in P \cap \mathbb{Z}^m$.

\begin{example}
Let $P\subset \mathbb{R}^2$ be the convex hull of points $(0,0), (1,0), (1,1)$, i.e. a unimodular simplex. Then $R(P)\subset \Bbbk[z_1,z_2,w]$ is generated by monomials $w; z_1w; z_1z_2w$. Clearly this algebra is freely generated by these monomials, so it is isomorphic to the polynomial ring in $3$ variables. All of them have degree 1, so the corresponding variety is isomorphic to a projective plane.
\end{example}

\subsection{Arc spaces}
In this subsection we recall the basic properties of arc spaces, cf.\ \cite{CLNS,EM,Na}.  
\begin{dfn}\label{thedef}
Let $V$ be a variety over $\Bbbk$. The arc (or jet-$\infty$) space $J^{\infty}(V)$ is a scheme  for which there exists a natural isomorphism of functors from $\mathbf{CAlg}_\Bbbk$ (commutative algebras) to $\mathbf{Set}$
\[\Hom(\Spec(-),J^{\infty}(V))\simeq\Hom(\Spec(-[[s]]),V).\]
\end{dfn}

Note that if $V=\Spec(R)$, this definition is equivalent to the existence of a natural isomorphism
\begin{equation}\label{JetDefiningHomIsom}
    \Hom(\Bbbk[J^{\infty}(V)],-)\simeq\Hom(R,-[[s]]).
\end{equation}
We immediately see that $J^{\infty}(V)$ is unique if it exists, since $\Bbbk[J^{\infty}(V)]$ represents the functor $\Hom(R,-[[s]])$. Accordingly, the natural isomorphisms in Definition~\ref{thedef} and in~\eqref{JetDefiningHomIsom} must also be unique. Throughout the paper we mostly deal with affine varieties (actually, with affine cones of projective varieties). The existence of $J^{\infty}(V)$ is deduced from the following explicit construction.

\begin{dfn}\label{ringdef}
For a finitely generated $\Bbbk$-algebra $R$ choose generators and relations:
\[R\simeq\Bbbk[X_1, \dots, X_m]/\langle f_1(X_1, \dots, X_m), \dots,
f_l(X_1, \dots, X_m)\rangle.\]
Consider the polynomial ring $\Bbbk[\{X_i^{(j)}\}]$ where $i=1, \dots, m$ and $j=0,1,\dots$. The \emph{arc ring} $J^{\infty}(R)$ is the quotient of $\Bbbk[\{X_i^{(j)}\}]$ by the ideal generated by all $s$-coefficients in series
\begin{equation}\label{JetRingRelations}
f_k(X_1(s), \dots, X_m(s))
\end{equation}
where $k=1, \dots,l$ and $X_i(s):=\sum_{j=0}^{\infty} X_i^{(j)}s^j$.
\end{dfn}

\begin{prop}\label{ExistenceArc}
For an affine variety $V=\Spec(R)$ the scheme
\[J^{\infty}(V)=\Spec J^{\infty}(R)\]
satisfies Definition~\ref{thedef}. In particular,
$J^{\infty}(R)$ does not depend on the chosen presentation of $R$.
\end{prop}

In the notations of Definition~\ref{ringdef} let $Y_i\in R$ be the image of $X_i$ and denote by $Y_i^{(j)}\in J^\infty(R)$ the image of $X_i^{(j)}$. Consider a $\Bbbk$-algebra $A$ and a homomorphism $\iota: R \to A[[s]]$ with
\begin{equation}\label{homToSeries}
\iota(Y_i)=\sum_{j=0}^{\infty} a_i^{(j)}s^j.
\end{equation}
Then we have a homomorphism $\iota':J^{\infty}(R)\to A$ with $\iota'(Y_i^{(j)})=a_i^{(j)}$
We obtain a morphism from $\Hom(J^{\infty}(R),A)$ to $\Hom(R,A[[s]])$ given by $\iota\mapsto\iota'$, together these morphisms provide the functorial isomorphism from Definition~\ref{thedef}.
We see that $J^\infty$ is a functor from the category of finitely generated commutative algebras to the category $\mathbf{CAlg}_\Bbbk$.
In geometric terms we can view $J^\infty$ as a functor from the category of varieties to the category of schemes
which satisfies $\Spec\circ J^\infty\simeq J^\infty\circ\Spec$.

\begin{example}
Consider $R=\Bbbk[X_{00},X_{01},X_{10},X_{11}]/\langle X_{00}X_{11}-X_{10}X_{01} \rangle$. Then the ring $J^\infty(R)$ is defined by relations which are $s$-coefficients of 
\[X_{00}(s)X_{11}(s)-X_{10}(s)X_{01}(s).\]
In other words, the ideal of relations is generated by the following expressions:
\[\sum_{i=0}^k X_{00}^{(i)}X_{11}^{(k-i)}-\sum_{i=0}^k X_{01}^{(i)}X_{10}^{(k-i)}.\]
\end{example}

Proposition~\ref{ExistenceArc} provides a canonical bijection from $\Hom (J^\infty(R), J^\infty(R))$ to $\Hom(R,J^\infty(R)[[s]])$. Let $\iota$ be the image of the identity under this bijection. For every $r \in R$ we obtain elements $r^{(j)}\in J^\infty(R)$ by setting 
\[\iota(r)=\sum_{j=0}^\infty r^{(j)}s^j=r(s).\]
Explicitly, 
if $R$ is generated by $Y_1,\dots,Y_m$ and
polynomial $p$ is such that $r=p(Y_1,\dots, Y_m)$, then $r(s)=p(Y_1(s),\dots, Y_m(s))$. In particular, $r\mapsto r(s)$ is injective.

Let $\varphi:R\to S$ be a homomorphism of finitely generated $\Bbbk$-algebras. The homomorphism $J^\infty(\varphi): J^\infty(R) \rightarrow J^\infty(S)$ is given by 
\begin{equation}\label{functorial}
J^\infty(\varphi)(r^{(j)})=\varphi(r)^{(j)}.
\end{equation}
It is clear that the functor $J^\infty$ preserves surjections. However, this functor need not preserve injections, such examples will be discussed in this paper. More precisely, we will consider situations where $R$ is embedded into a polynomial ring $S$ but $J^\infty(R)$ has a non-trivial nilradical which, in particular, implies that $J^\infty(R)$ cannot be a subring of $J^\infty(S)$ (which is itself a polynomial ring). See, for example, Lemma~\ref{NilpotentSeriesinCube}.

The fact that the arc space of a reduced ring can be non-reduced is, in a sense, the only obstacle to preserving injections (at least if the rings are nice enough). We have a functor $J^\infty_\red$ from $\mathbf{CAlg}_\Bbbk$ to itself which takes $R$ to its \textit{reduced arc ring}, i.e.\ $J^\infty(R)$ modulo its nilradical. We will now show that  $J^\infty_\red$ preserves injections.

\begin{thm} \label{dominance theorem}
	Suppose $\phi: X \rightarrow Y$ is a dominant morphism of varieties. 
 Then the corresponding morphism of arc spaces $J^{\infty}(\phi): J^{\infty}(X) \rightarrow J^{\infty}(Y)$ is dominant.
\end{thm}

\begin{proof}
We note that the restriction of $\phi$ to the smooth locus  $X^{\mathrm{reg}}$ gives a dominant map $X^{\mathrm{reg}} \rightarrow Y$. 
To prove the theorem it suffices to show that $J^\infty(X^{\mathrm{reg}}) \rightarrow J^\infty(Y)$ is dominant. Hence, we may (and will) assume that $X$ is smooth. 

	By the Chevalley constructibility theorem, the dominance is equivalent to the fact that $\phi(X)$ contains a dense open subset of $Y$, denote it by $U$. Then the restriction of $\phi$ to $\phi^{-1}(U)$ can be decomposed as $\phi^{-1}(U) \xrightarrow{\phi_1} U \xrightarrow{\phi_2} Y$ with surjective $\phi_1$ and open embedding $\phi_2$. Using that $\phi^{-1}(U) \subset X$ is smooth, we can apply the ``generic smoothness on the target'' theorem \cite[Theorem 25.3.3]{Vak} to $\phi_1$, and obtain that there is an open dense $V \subset U$ such that $\phi_1$ is smooth on $\phi_1^{-1}(V)$. Restricting to $\phi_1^{-1}(V)$, we have:
	\[
	\phi_1^{-1}(V) \xrightarrow{\phi_1} V \xrightarrow{\psi} U \xrightarrow{\phi_2} Y,
	\]
	where $\phi_1$ is surjective smooth and $\psi, \phi_2$ are open embeddings. 
	
	Passing to arcs, we get that $J^{\infty}(\phi_1)$ is surjective \cite[Remark 2.10]{EM} and $J^{\infty}(\psi), J^{\infty}(\phi_2)$ are open embeddings \cite[Lemma 2.3]{EM}. Note also that by the Kolchin's theorem \cite[Theorem 3.3]{EM} $J^{\infty}(U), J^{\infty}(Y)$ are irreducible and hence the open embeddings $J^{\infty}(\psi), J^{\infty}(\phi_2)$ are dominant. Therefore, the composition is dominant and the theorem is proven.
\end{proof}


In this paper we study the ring $J_{\red}^\infty(R)$ for a toric ring $R$ embedded into  a polynomial algebra. The following corollary implies that the reduced arc ring $J_{\red}^\infty(R)$ can be described using arc map of the embedding.

\begin{cor}\label{jredinj}
	Suppose $f: R \hookrightarrow L$ is an injective map of integral finitely generated $\Bbbk$-algebras. 
    Then the corresponding map $J_{\red}^{\infty}(R) \rightarrow J_{\red}^{\infty}(L)$ is also injective.
\end{cor}

\begin{proof}
	The corresponding map $f^*: \Spec L \rightarrow \Spec R$ is dominant. Applying Theorem \ref{dominance theorem} we obtain that $J^{\infty}(f)^*: \Spec J^{\infty}(L) \rightarrow \Spec J^{\infty}(R)$ is dominant. But that is equivalent to the fact that $\ker (J^{\infty}(f))$ is contained in the nilradical of $J^{\infty}(R)$ and we are done.
\end{proof}

The next proposition shows how the arc space of a closure is related to the closure of the arc space.

\begin{prop} \label{closure and arcs}
Let $\phi: X \rightarrow Y$ be a morphism of varieties.
Then
\[
J^\infty( \overline{\im \phi}) = \overline{ \im( J^\infty(\phi) )},
\]
where both sides are endowed with the reduced scheme structure.

\end{prop}
\begin{proof}
Note that for any morphism of schemes $A \rightarrow B$ there exists unique decomposition $A \rightarrow C \rightarrow B$ such that the first arrow is dominant and the second is a closed embedding. We prove that both decompositions 
\begin{equation} \label{first decomposition}
J^\infty(X) \xrightarrow{J^\infty(\phi)} \ov{ \im( J^\infty(\phi))} \hookrightarrow J^\infty(Y)
\end{equation}
and
\begin{equation} \label{second decomposition}
J^\infty(X) \xrightarrow{J^\infty(\phi)} J^\infty( \ov{\im \phi}) \hookrightarrow J^\infty(Y)
\end{equation}
satisfy this property. This will imply the proposition.

Decomposition \eqref{first decomposition} satisfies the desired properties just by definition (the first map is dominant and the second is a closed embedding).

In order to obtain the same for \eqref{second decomposition} we first consider the decomposition $X \rightarrow \ov{\im \phi} \rightarrow Y$ of the map $\phi$, this decomposition has the desired properties. Applying the functor $J^\infty$ we obtain that in decomposition \eqref{second decomposition} the first arrow is dominant due to Theorem \ref{dominance theorem} and the second is a closed embedding.
\end{proof}

\begin{cor}\label{ArcHomogenousSpace}
Let $X$ be a variety with an action of an algebraic group $G$. This induces the action of $J^\infty(G) = G[[t]]$ on $J^\infty(X)$. Then for any point $x \in X \hookrightarrow J^\infty(X)$ and orbits $G.x \subset X$, $G[[t]].x \subset J^\infty(X)$ one has: $J^\infty( \ov{G.x} ) = \ov{G[[t]].x}$, where both sides are endowed with the reduced scheme structure.
\end{cor}
\begin{proof}
Use Proposition \ref{closure and arcs} for the map $G \rightarrow X$ defined by $g \mapsto g.x$.
\end{proof}



The  corollary above implies that the arc ring $J^\infty_\red(R)$ of an affine toric variety, which we mainly study in this paper, is isomorphic to the coordinate ring of the cone over the closure of the orbit of the toric arc group.

To end this section we provide an example of how one may find relations in a ring $J^\infty_\red(R)$. We first prove two general (and well-known) lemmas concerning an arbitrary commutative $\Bbbk$-algebra $A$. The lemmas can be seen to be equivalent but it is convenient for us to use them both.
\begin{lem}\label{nilpotentCoefficient}
Consider elements $a_{i} \in A, i=0, 1, \dots$.
Assume that the series $a(s)=\sum_{i=0}^\infty a_{i} s^i$ is nilpotent. Then all the coefficients $a_{i}$
are nilpotent.
\end{lem}
\begin{proof}
Clearly $a_{0}$ is nilpotent. Assume that $a_{i}$, $i<j$ are nilpotent. Then
\[a(s)-\sum_{i=0}^{j-1} a_{i} s^i\]
is nilpotent and thus its leading term $a_{j}s^j$ is niplotent. This completes the proof by induction.
\end{proof}

\begin{lem}\label{NilpotentDerivations}
Let $d$ be a $\Bbbk$-linear derivation of the algebra $A$. If $a\in A$ is nilpotent, then
$d(a)$ is nilpotent.
\end{lem}
\begin{proof}
Suppose $a^n=0$, then $d^n(a^n)=0$. Expanding $d^n(a^n)$ via the Leibniz rule we obtain $n!d(a)^n+ab$ for some $b\in A$. Hence $d(a)^n=-ab/n!$ is nilpotent and so is $d(a)$.
\end{proof}

\begin{example}[\cite{BS2}]
Here and further we write $a \equiv b$ for elements $a,b$ of a ring to denote that $a-b$ is nilpotent.
Consider the ring $R$ generated by $Y_1$ and $Y_2$ satisfying the single relation $Y_1Y_2=0$. Then the ring 
$J^\infty(R)$ is generated by $Y_1^{(j)}, Y_2^{(j)}$ satisfying the relations
\[Y_1(s)Y_2(s)=0.\]
We show by induction that for $n \geq 0$ the series
\begin{equation}\label{ProductWithDerived}
Y_1(s) \frac{\partial^n Y_2(s)}{\partial s^n}
\end{equation}
is nilpotent in $J^\infty(R)[[s]]$.
Indeed, assume that $Y_1(s) \frac{\partial^{n-1} Y_2(s)}{\partial s^{n-1}}$ is nilpotent. Using Lemma \ref{NilpotentDerivations} we have
\[Y_1(s) \frac{\partial^n Y_2(s)}{\partial s^n}\equiv -
\frac{\partial Y_1(s)}{\partial s} \frac{\partial^{n-1} Y_2(s)}{\partial s^{n-1}}.\]
Multiplying both sides by the left-hand side part we obtain
\[\left(Y_1(s) \frac{\partial^n Y_2(s)}{\partial s^n}\right)^2\equiv -Y_1(s) \frac{\partial^n Y_2(s)}{\partial s^n}
\frac{\partial Y_1(s)}{\partial s} \frac{\partial^{n-1} Y_2(s)}{\partial s^{n-1}}.\]
By the induction hypothesis the right-hand side is nilpotent, hence so is $Y_1(s) \frac{\partial^n Y_2(s)}{\partial s^n}$. 
In fact, one may check that the coefficients of the series \eqref{ProductWithDerived} span the same linear space as the elements $Y_1^{(i)}Y_2^{(j)}$ for all $i,j\ge0$. Thus the quotient of $J^\infty(R)$ modulo its nilradical is spanned by the images of all
$Y_1^{(i_1)}\dots Y_1^{(i_M)}$ and $Y_2^{(j_1)}\dots Y_2^{(j_N)}$. Using this and the results of Subsection \ref{subsection arcs of monomial ideals} one can show that the nilradical is generated by coefficients of series \eqref{ProductWithDerived}.
\end{example}

\subsection{Derivations of arc rings}

In this subsection we define two Lie algebra actions on an arbitrary arc ring, these actions will be instrumental to our construction.

Consider any commutative $\Bbbk$-algebra $\Theta$ and an odd formal variable $\tau$ (i.e. $\tau^2=0$). Let  $\varkappa:\Theta[\tau] \twoheadrightarrow \Theta$ be the natural projection, i.e.\ $\varkappa(\Theta\tau)=0$. In the following lemma and throughout the text we only consider $\Bbbk$-derivations.
\begin{lem}\label{DerHom}
The derivations of $\Theta$ are in one-to-one correspondence with the homomorphisms $\delta:\Theta \rightarrow \Theta[\tau]$ for which $\varkappa \circ \delta = 1$. A derivation $d$ corresponds to the homomorphism $1+\tau d$. 
\end{lem}
\begin{proof}
For any derivation $d \in \Der \Theta$ one evidently has \[\varkappa \circ (1+ \tau d) =1.\]
Suppose that we have a homomorphism $\delta:\Theta \rightarrow \Theta[\tau]$ such that $\varkappa \circ \delta = 1$. Then we have
\[\delta=1+\tau d\]
for a linear map $d$. For any $x,y \in \Theta$:
\[xy +\tau d(xy)=(1+\tau d)(xy)=(1+\tau d)(x)(1+\tau d)(y)=xy +\tau(xd(y)+d(x)y).\]
Hence $d$ is a derivation.
\end{proof}

Fix a finitely generated $\Bbbk$-algebra $R$ for the rest of this section. Consider a continuous derivation $d \in \Der\Bbbk[[s]]$, we can extend $d$ to a derivation of $J^{\infty}(R)[[s]]=J^{\infty}(R)\widehat\otimes\Bbbk[[s]]$. Here continuity is with respect to the standard topology on $\Bbbk[[s]]$ and  $\widehat\otimes$ denotes the completed tensor product with respect to this topology, see, for instance,~\cite{SP}. 
Lemma~\ref{DerHom} applied to $\Theta=J^{\infty}(R)[[s]]$ and the obtained derivation provide maps
\[
J^{\infty}(R)[[s]] \rightarrow J^{\infty}(R)[[s]][\tau] \rightarrow J^{\infty}(R)[[s]],
\]
and hence the maps
\[
\Hom(R,J^{\infty}(R)[[s]]) \rightarrow \Hom(R,J^{\infty}(R)[[s]][\tau]) \rightarrow  \Hom(R,J^{\infty}(R)[[s]]),
\]
with composition equal to identity. Applying isomorphism \eqref{JetDefiningHomIsom} we have maps $\varphi$ and $\pi$
\begin{equation}\label{composition}
\Hom(J^{\infty}(R), J^{\infty}(R)) \xrightarrow{\varphi} \Hom(J^{\infty}(R),J^{\infty}(R)[\tau]) \xrightarrow{\pi} \Hom(J^{\infty}(R), J^{\infty}(R))
\end{equation}
with composition equal to identity.
\begin{prop}
The image of the identity under $\varphi$ has the form $1+\tau\rho(d)$ for some $\rho(d)\in\Der J^{\infty}(R)$. The map $-\rho$ defines a Lie algebra action of $\Der^c \Bbbk[[s]]$ (continuous derivations) on $J^{\infty}(R)$.
\end{prop}
\begin{proof}
By construction, the map $\pi$ is obtained by applying the functor $\Hom(J^\infty(R),-)$ to the projection $\varkappa:J^{\infty}(R)[\tau]\twoheadrightarrow J^{\infty}(R)$, i.e.\ it is the composition with $\varkappa$ on the left. In particular, if $\delta=\varphi(\id_{J^\infty(R)})$, then $\id_{J^\infty(R)}=\varkappa\delta$. By Lemma~\ref{DerHom}, we indeed have $\delta=1+\tau\rho(d)$ for some $\rho(d)\in\Der J^{\infty}(R)$.

Let us show that \[-\rho:\Der^c\Bbbk[[s]]\to\Der J^\infty(R)\] is a Lie algebra homomorphism. The derivation $d$ is given by a $(\bZ_{\ge0}\times\bZ_{\ge0})$-matrix $B$ such that \[d(s^j)=\sum_{i\ge 0} B_{i,j}s^i\] for any $j\ge 0$, i.e. $d$ maps a series with coefficients $(c_0,c_1,\dots)$ to the series with coefficients $B (c_0,c_1,\dots)^T$. From the definitions we see that for any $r\in R$ and $i\ge 0$ we then have \[\varphi(\id_{J^\infty(R)})(r^{(i)})=r^{(i)}+\tau\sum_{j\ge 0}B_{i,j}r^{(j)}.\] In other words, $\rho(d)$ maps a linear combination of the $r^{(j)}$ with coefficients $(c_0,c_1,\dots)$ to their linear combination with coefficients $(c_0,c_1,\dots) B$. This shows that $\rho$ is an antihomomorphism, i.e\ $-\rho$ is indeed a homomorphism.
\end{proof}
We will use the shorthand notation $dX=-\rho(d)(X)$ for $d\in\Der^c\Bbbk[[s]]$ and $X\in J^\infty(R)$. Recall that the Lie algebra $\Der^c\Bbbk[[s]]$ is spanned by the derivations $d_k=s^{k+1}\frac{\partial}{\partial s}$ with $k\ge-1$. These derivations act on $J^\infty(R)$ as follows (here and further $r^{(i)}=0$ for any $r\in R$ and $i<0$).
\begin{prop}\label{dkaction}
For any $r\in R$, $i\ge 0$ and $k\ge -1$ we have
\begin{equation}\label{LieDerivationsAction}
d_k(r^{(i)})=(i-k)r^{(i-k)}.
\end{equation}
\end{prop}
\begin{proof}
Setting $d=d_k$ within the proof of the previous proposition we see that the matrix $B$ is given by $B_{i,j}=j$ when $i-j=k$ and $B_{i,j}=0$ otherwise. In particular, the proof then implies that \[d_k(r^{(i)})=\sum_{j>0} B_{i,j}r^{(j)}=(i-k)r^{(i-k)}.\qedhere\]
\end{proof}

From \eqref{functorial} and Proposition \ref{dkaction} we see that $J^{\infty}$ respects the action of $\Der^c\Bbbk[[s]]$. In other words, $J^\infty$ as a functor from $\mathbf{CAlg}_\Bbbk^\mathrm{fin}$ to the category of commutative $\Der^c\Bbbk[[s]]$-differential algebras. 

We next describe another Lie algebra action on $J^\infty(R)$. Fix a presentation \[R\simeq\Bbbk[X_1, \dots, X_m]/\langle f_1(X_1, \dots, X_m), \dots, f_l(X_1, \dots, X_m)\rangle\] and let $Y_i\in R$ be the image of $X_i$.
\begin{dfn}\label{currentDerivationDefinition}
For a derivation $\mathfrak d \in \Der R$ and $k\ge 0$ we consider $\mathfrak d^{(k)}\in\Der J^{\infty}(R)$ defined on generators by
\[\mfd^{(k)}(Y_i^{(j)}):=\mfd(Y_i)^{(j-k)}.\]
\end{dfn}

\begin{lem}
The derivations $\mathfrak d^{(k)}$ are well-defined and independent of the presentation.
\end{lem}
\begin{proof}
Choose polynomials $p_1,\dots,p_m$ in $\Bbbk[X_1, \dots, X_m]$ such that $p_i(Y_1,\dots,Y_m)=\mfd Y_i$ for all $Y_i$. We have a unique derivation $\tilde\mfd$ of $\Bbbk[X_1,\dots,X_m]$ given by $\tilde\mfd X_i=p_i$. We also consider the derivation $\tilde\mfd^{(k)}$ of $\Bbbk[\{X_i^{(j)}\}]$ such that $\tilde\mfd^{(k)}(X_i^{(j)})$ is the coefficient of $s^{j-k}$ in $p_i(X_1(s),\dots,X_m(s))$ or 0 if $j<k$. Furthermore, we extend $\tilde\mfd^{(k)}$ to a derivation of $\Bbbk[\{X_i^{(j)}\}][[s]]$ by applying it coefficientwise. Note that 
\begin{equation}\label{seriesderiv}
\tilde\mfd^{(k)}(X_i(s))=s^kp_i(X_1(s),\dots,X_m(s)).
\end{equation}

Now, for any $p\in\Bbbk[X_1, \dots, X_m]$ its image $\tilde\mfd(p)$ projects to $\mfd(p(Y_1,\dots,Y_m))\in R$. Hence $\tilde\mfd$ preserves the ideal of relations in $R$, i.e.\ for all $f_o$ we have polynomials $a_q$ such that 
\[\tilde\mfd(f_o)=\sum_{q=1}^la_q f_q.\]
From~\eqref{seriesderiv} we then obtain
\[
\tilde\mfd^{(k)}\Big(f_o\big(X_1(s), \dots, X_m(s)\big)\Big)=s^k\sum_{q=1}^la_q\big(X_1(s), \dots, X_m(s)\big) f_q\big(X_1(s), \dots, X_m(s)\big).
\]
This means that $\tilde\mfd^{(k)}\in\Der\Bbbk[\{X_i^{(j)}\}]$ preserves the defining ideal of $J^\infty(R)$ generated by all coefficients of the series $f_q(X_1(s), \dots, X_m(s))$. By projecting onto  $J^\infty(R)$ we obtain a derivation which maps $Y_i^{(j)}$ to the coefficient of $s^{j-k}$ in $p_i(Y_1(s), \dots, Y_m(s))$. The latter coefficient is precisely $\mfd(Y_i)^{(j-k)}$ and the obtained derivation satisfies the definition of $\mfd^{(k)}$.

To show that $\mfd^{(k)}$ is independent of the presentation extend  $\mfd^{(k)}$ to a derivation of $J^\infty(R)[[s]]$ coefficientwise. We have the embedding $\iota$ of $R$ into $J^\infty(R)[[s]]$ with $\iota(r)=r(s)$. Now we note that \[s^k\iota\mfd(Y_i)=\mfd^{(k)}\iota(Y_i),\] i.e. the maps $s^k\iota\mfd\iota^{-1}$ and $\mfd^{(k)}$ agree on the $\iota(Y_i)$, hence on all of $\iota(R)$. In other words, \begin{equation}\label{indepderdef}
\mfd^{(k)}(r^{(j)})=\mfd(r)^{(j-k)}
\end{equation}
for any $r\in R$ and $j\ge 0$ which is independent of the presentation.
\end{proof}

\begin{lem}\label{currentaction}
The current Lie algebra $(\Der R) [s]=\Der R \otimes_{\Bbbk}\Bbbk[s]$ acts by derivations on $J^{\infty}(R)$
by $\mfd s^k$ acting as $\mfd^{(k)}$.
\end{lem}
\begin{proof}
For any two derivations $\mfd_1, \mfd_2 \in \Der(D)$ and $k_1, k_2\ge 0$ we need to prove
\[[\mfd_1^{(k_1)},\mfd_2^{(k_2)}]=[\mfd_1,\mfd_2]^{(k_1+k_2)}.\]
However, from \eqref{indepderdef} we indeed have
for any $r \in R$:
\[[\mfd_1^{(k_1)},\mfd_2^{(k_2)}](r^{(j)})=([\mfd_1,\mfd_2](r))^{(j-k_1-k_2)}=[\mfd_1,\mfd_2]^{(k_1+k_2)}(r^{(j)}).\]
Since the $r^{(j)}$ generate $J^\infty(R)$, the two derivations agree on all of $J^\infty(R)$.
\end{proof}


Let $\tilde \tau:\mathfrak{a}\rightarrow \mathfrak{b}$ be a morphism of Lie algebras. Then there exists the morphism
\[\tilde \tau[s]:\mathfrak{a}[s]\rightarrow \mathfrak{b}[s],~as^k\mapsto \tilde \tau(a)s^k.\]

\begin{lem}\label{equivariantCurrentAction}
Let $R,L$ be associative algebras, $\fg_R,\fg_L$
be Lie algebras acting on them by derivations.
Let $\tau$, $\tilde \tau$ be two morphisms $\tau:R \rightarrow L$,
$\tilde \tau:\fg_R\rightarrow \fg_L$ such that $\tau$ is equivariant with respect to $\tilde \tau$, i.e.\ for any $d \in \fg_R$, $r \in R$
\[\tau(d(r))=\tilde \tau(d)\big(\tau(r)\big).\]
Then $J^{\infty}(\tau)$ is equivariant with respect to $\tilde\tau[s]$.
\end{lem}
\begin{proof}
This follows from 
\eqref{indepderdef}
and the definition of equivariance.
\end{proof}




\subsection{Supersymmetric polynomials}
In this subsection we discuss the properties of super-symmetric polynomials.
Take two groups of variables $p^{(1)},\dots,p^{(b)}$, $q^{(1)},\dots,q^{(c)}$.
\begin{dfn}\label{superSymmetricFunctions}
A polynomial $f \in \Bbbk[p^{(1)},\dots,p^{(b)}, q^{(1)},\dots,q^{(c)}]^{\mathfrak {S}_b \times \mathfrak {S}_c}$ symmetric in the $p^{(i)}$ and the $q^{(j)}$ is supersymmetric if $f|_{p^{(b)}\mapsto Z,q^{(c)}\mapsto Z}$ is
independent of $Z$. We denote the ring of such polynomials by $\Omega_{b,c}$.
\end{dfn}
\begin{prop}\label{resultant}
There exists the map $\psi_{b,c}:\Omega_{b,c}\rightarrow \Omega_{b-1,c-1}$:
\[\psi_{b,c}:f \mapsto f|_{p^{(b)}\mapsto Z,q^{(c)}\mapsto Z}.\]
The kernel of $\psi_{b,c}$ is the ideal 
\[\prod_{1 \leq i \leq b, 1 \leq j \leq c}\left(p^{(i)}-q^{(j)}\right)\Bbbk[p^{(1)},\dots,p^{(b)}, q^{(1)},\dots,q^{(c)}]^{\mathfrak {S}_b \times \mathfrak {S}_c}.\]
\end{prop}
\begin{proof}
For a polynomial $f \in \Omega_{b, c}$ the image $\psi_{b,c}(f)$ is symmetric in $p^{(1)}, \dots, p^{(b-1)}$ and $q^{(1)}, \dots, q^{(c-1)}$. 
We have:
\[\psi_{b,c}(f)|_{p^{(b-1)}\mapsto Z',q^{(c-1)}\mapsto Z'}=f|_{p^{(b)}\mapsto Z,q^{(c)}\mapsto Z,p^{(b-1)}\mapsto Z',q^{(c-1)}\mapsto Z'}.\]
Therefore it is independent of $Z'$. Hence $\psi_{b,c}(f)\in \Omega_{b-1,c-1}$.

Assume now $f \in \ker(\psi_{b,c})$. Since $f$ is divisible by $\left(p^{(b)}-q^{(c)}\right)$ and symmetric in $p^{(1)}, \dots p^{(b)}$ and $q^{(1)}, \dots q^{(c)}$, 
we see that $f$ is divisible by $\prod_{1 \leq i \leq b, 1 \leq j \leq c}\left(p^{(i)}-q^{(j)}\right)$. Conversely 
\[\prod_{1 \leq i \leq b, 1 \leq j \leq c}\left(p^{(i)}-q^{(j)}\right)_{p^{(b)}\mapsto Z,q^{(c)}\mapsto Z}=0.\]
Thus the principal ideal $\prod_{1 \leq i \leq b, 1 \leq j \leq c}\left(p^{(i)}-q^{(j)}\right)\Omega_{b,c}$ is annihilated by $\psi_{b,c}$.
\end{proof}

The proof of the following proposition can be found in \cite[Section 1.3.23]{M}.
\begin{prop}\label{SuperSymGenerators}
The ring $\Omega_{b, c}$ has two following sets of generators:
\[p_k\left(p^{(1)},\dots, p^{(b)};q^{(1)}, \dots, q^{(c)}\right)=\sum_{i=1}^b\left(p^{(i)}\right)^{k}-\sum_{i=1}^c\left(q^{(i)}\right)^{k}\]
and
$h_k(p^{(1)},\dots, p^{(b)};q^{(1)}, \dots, q^{(c)})$, where
\[\frac{\prod_{i=1}^b\left(1-Tp^{(i)}\right)}{\prod_{i=1}^c\left(1-Tq^{(i)}\right)}=
\sum_{k=0}^{\infty} h_k\left(p^{(1)},\dots, p^{(b)};q^{(1)}, \dots, q^{(c)}\right)T^k\]
for a formal variable $T$. In particular $\psi_{b,c}$ is surjective.
\end{prop}

\subsection{Associated graded algebras and initial ideals}\label{initial}
Consider the polynomial ring \linebreak $\Bbbk[X_1, \dots, X_n]$. Let $\prec$
be a monomial order on monomials in $X_1, \dots, X_n$, i.e.\
a total order on monomials such that if $a \prec b$ then for any monomial $c$: $ac 
\prec bc$. In other words, $\prec$ is a total semigroup order on $\mathbb{N}^n$ (see \cite{MS}).

The ring $\Bbbk[X_1, \dots, X_n]$ is naturally $\mathbb{N}^n$-graded. All graded components are one-dimensional and are spanned by monomials:
\[\Bbbk[X_1, \dots, X_n][\bar r]=\Bbbk\prod X_i^{r_i}.\]
The ring $\mathbb J_n=J^{\infty}(\Bbbk[X_1, \dots, X_n])$ is the polynomial ring in variables $X_i^{(j)}$ with $1\le i\le n$ and $j\ge 0$. This ring has an induced $\mathbb{N}^n$-grading where graded component $\mathbb J_n[\bar r]$ has the following basis:
\[\left\{\prod_{1 \leq i \leq n, j=0,1,\dots} \left(X_{i}^{(j)}\right)^{b_i^j}|~\sum_{j=0}^\infty b_i^j=r_i, i=1, \dots, n\right\}.\]

$\mathbb J_n$ has two filtrations by the ordered semigroup $(\bN^n,\prec)$. The first has components
\begin{equation*}
\mathbb J_n[\prec\bar r]=\bigoplus_{\bar r'\prec \bar r} J^{\infty}(\Bbbk[X_1, \dots, X_n])[\bar r'],
\end{equation*}
while the second has components
\begin{equation*}
\mathbb J_n[\preceq\bar r]=\bigoplus_{\bar r'\preceq \bar r} J^{\infty}(\Bbbk[X_1, \dots, X_n])[\bar r'].
\end{equation*}
Both filtrations are multiplicative: for any $\bar r_1,\bar r_2\in\bN^n$ we have 
\[\mathbb J_n[\prec\bar r_1]\mathbb J_n[\prec\bar r_2]\subset\mathbb J_n[\prec\bar r_1+\bar r_2]~ \text{and} ~\mathbb J_n[\preceq\bar r_1]\mathbb J_n[\preceq\bar r_2]\subset\mathbb J_n[\preceq\bar r_1+\bar r_2].\]


Let $\mathcal I \subset \mathbb J_n$ be an ideal.
Denote $R=\mathbb J_n/\mathcal I$ and let $\pi$ be the natural projection
$\pi:\mathbb J_n\twoheadrightarrow R$.
The filtrations on $\mathbb J_n$ induce filtrations on
$R$ by $R[\prec\bar r]=\pi(\mathbb J_n[\prec\bar r])$ and $R[\preceq\bar r]=\pi(\mathbb J_n[\preceq\bar r])$. These filtrations are again multiplicative and we have the associated $\bN^n$-graded algebra \[\gr_\prec R=\bigoplus_{\bar r\in\bN^n} R[\preceq\bar r]/R[\prec\bar r].\]

Next, for any nonzero element $g \in \mathbb J_n$ consider its decomposition
\[g=\sum_{\bar r \in \mathbb{N}^n} g[\bar r],\, g[\bar r] \in \mathbb J_n[\bar r].\]
Let $\bar r'$ be the largest $\bar r$ with respect to $\prec$ such that $g[\bar r]$ is nonzero.
Then $g[\bar r']$ is denoted $\In_{\prec} g$ and is called the \textit{initial part} of $g$. The \textit{initial ideal} of $\mathcal I$ with respect to $\prec$ is the linear span of
$\{\In g| g \in \mathcal I\}$. This space is denoted $\In_\prec\mathcal I$, one easily checks that it is an ideal in $\mathbb J_n$. By construction, the ideal $\In_\prec \mathcal{I}$ is homogeneous with respect to the $\bN^n$-grading. Therefore, the ring $\mathbb J_n/ (\In_{\prec}\mathcal{I})$ is naturally $\mathbb{N}^n$-graded.

Next, let us observe that $\mathbb J_n$ has a further $\bN$-grading $\grad$ given by $\grad X_i^{(j)}=j$. Together with the grading considered above this turns $\mathbb J_n$ into a $\bN^{n+1}$-graded ring. We now assume that $\mathcal I$ is homogeneous with respect to $\grad$. This gives us an induced $\bN$-grading $\grad$ on $R$. Note that every subspace $\mathbb J_n[\prec\bar r]$ and $\mathbb J_n[\preceq\bar r]$ is also $\grad$-homogeneous, hence so is every $R[\prec\bar r]$ and $R[\preceq\bar r]$. As a result we obtain an induced grading $\grad$ on $\gr_\prec R$ which makes the latter $\bN^{n+1}$-graded. Furthermore, since $\mathcal I$ is $\grad$-homogeneous, so is any initial ideal. Consequently, $\In_\prec\mathcal I$ and $\mathbb J_n/\In_\prec\mathcal I$ are $\bN^{n+1}$-graded.

\begin{prop}\label{UpperBoundInitial}
The rings $\gr_\prec R$ and $\mathbb J_n/\In_\prec\mathcal I$ are isomorphic as $\bN^{n+1}$-graded algebras.
\end{prop}

\begin{proof}
We may consider the associated graded algebra $\gr_\prec\mathbb J_n$, note that it is naturally identified with $\mathbb J_n$. Since the filtrations on $R$ are defined as projections of the filtrations on $\mathbb J_n$, we obtain a surjective homomorphism $\pi':\gr_\prec\mathbb J_n\to\gr_\prec R$ preserving the $\bN^{n+1}$-grading. 
In terms of this identification, $p\in\mathbb J_n[\bar r]$ lies in $\ker\pi'$ if and only if $p-q\in\mathbb J_n[\prec\bar r]$ for some $q\in\mathcal I$. In other words: $p\in\ker\pi'$ if and only if $p=\In_\prec q$ for some $q\in\mathcal I$. This proves the proposition.
\end{proof}

\begin{rem}
This proposition is a special case of the well-known and rather general phenomenon of isomorphisms between associated graded rings and quotients by initial ideals. See, for instance,~\cite[Lemma 3.4]{KaM} or~\cite[Propositions 1.1 and 8.1]{Ma} for statements of this form for finitely generated commutative and associative algebras.
\end{rem}

Choose a $\grad$-homogeneous generating set $\{f_1, \dots, f_l\}$ of $\mathcal{I}$. We denote by $\mathcal{J}$ the ideal generated by the polynomials $\In_\prec f_i$. This ideal is $\mathbb{N}^{n+1}$-homogeneous and is contained in $\In_\prec\mathcal I$. This provides a surjection $\mathbb J_n/\mathcal J\twoheadrightarrow \mathbb J_n/\In_\prec\mathcal I$ of $\mathbb{N}^{n+1}$-graded algebras. 
In view of Proposition~\ref{UpperBoundInitial}, we have the following.

\begin{cor}\label{UpperBoundSubquotient}
There exists a surjective homomorphism of $\bN^{n+1}$-graded algebras from $\mathbb J_n/\mathcal J$ to $\gr_\prec R$.
\end{cor}

\subsection{Unit cubes}
\label{hibireduced}

For an integer $l\ge 2$ let $\mathbb R^{[l]}$ be the space of real-valued functions on the set of subsets of $[l]$. Consider the unit cube $\mathcal C_l\subset\mathbb R^{[l]}$, its vertices are the indicator functions $\mathbf 1_I$ with $I\subset [l]$. $R(\mC_l)$ is generated by elements $Y_I$ which satisfy
\begin{equation}
Y_IY_J-Y_{I \cap J}Y_{I \cup J}.
\end{equation}
Consider the arc ring $J^{\infty}(R(\mC_l)$). It is generated by variables
$Y_I^{(j)}, j \ge 0$, satisfying the relations
\begin{equation}\label{archibiboolean}
Y_I(s)Y_J(s)-Y_{I \cap J}(s)Y_{I \cup J}(s)=0.
\end{equation}


\begin{lem}\label{NilpotentSeriesinCube}\label{boolnilpotents}
The coefficients of the following series are nilpotent in $J^{\infty}(R(\mC_l))$ for any $0\le k \leq l-2$:
\begin{equation}\label{Wdefinition}
W_{l,k}=\sum_{I \subset [2,l]}(-1)^{|I|}Y_{I \cup \{1\}}(s)
\frac{\partial^k Y_{[2,l]\backslash I}(s)}{\partial s^k}.
\end{equation}
\end{lem}
\begin{proof}
We prove this lemma by induction on $l$ by showing that the series $W_{l,k}$ is nilpotent in $J^{\infty}(R(\mC_l))[[s]]$, the lemma then follows by Lemma \ref{nilpotentCoefficient}. In the case $l=2$ we only have $k=0$ and $W_{l,k}$ vanishes by~\eqref{archibiboolean}.

For the induction step we first consider the case $k<l-2$. We break up $W_{l,k}$ into two sums and show that each summand can be interpreted as a similar relation for $l-1$. First consider the sum
\begin{equation}
\sum_{J \subset [2,l-1]}(-1)^{|J|+1}Y_{J \cup \{1,l\}}(s)
\frac{\partial^k Y_{[2,l-1]\backslash J}(s)}{\partial s^k}.\label{hypothesis}
\end{equation}
Every $I\subset[l]$ corresponds to a vertex $v^I$ of $\mC_l$ such that $v^I_i$ is 1 if $i\in I$ and 0 otherwise. The subscripts appearing in~\eqref{hypothesis} are those subsets which either contain both 1 and $l$ or neither of them. The corresponding vertices are the vertices of an $(l-1)$-dimensional parallelepiped $D$. Note that $R(D)$ is embedded into $R(\mC_l)$ as the subalgebra generated by the $Y_I$ with $v^I\in D$. Consequently, we can view $J^\infty(R(D))$ as the subalgebra in $J^\infty(R(\mC_l))$ generated by the corresponding $Y_I^{(j)}$.

Now, $D$ can be identified with $\mathcal C_{l-1}$ by projecting along the $l$th coordinate. This projection is unimodular (it identifies the lattice $\spanned(D)\cap\bZ^l$ with $\bZ^{l-1}$), therefore, it induces isomorphisms $R(D)\simeq R(\mathcal C_{l-1})$ and $J^\infty(R(D))\simeq J^\infty(R(\mathcal C_{l-1}))$. The image of~\eqref{hypothesis} under the latter isomorphism is $-W_{l-1,k}$ which is nilpotent by the induction hypothesis.

The summands in $W_{l,k}$ not found in~\eqref{hypothesis} are those featuring subscripts which contain exactly one of 1 and $l$. The corresponding vertices again form the vertex set of an $(l-1)$-dimensional parallelepiped and the same argument shows that the remaining sum is also nilpotent. 


Now we consider $k=l-2$. First, for integers $n\in[2,l]$, $m\ge 0$ and $p\in[0,m]$ consider the expression \[\widetilde W_{n,m,p}=\sum_{I \subset [l-n+2,l]}(-1)^{|I|}
\frac{\partial^{p}Y_{I\cup \{l-n+1\}}(s)}{\partial s^p}
\frac{\partial^{m-p} Y_{[l-n+2,l]\backslash I}(s)}{\partial s^{m-p}}\in J^{\infty}(R(\mC_l))[[s]].\] In particular, $\widetilde W_{l,m,0}=W_{l,m}$ for $m\le l-2$. 

Let $n$ be either $l-1$ or $l$. By, respectively, either the induction hypothesis applied to the Boolean lattice of subsets in $[2,l]$ or the above case $k<l-2$, we have $\widetilde W_{n,m,0}\equiv 0$ for $m\le l-3$ (recall that $a\equiv b$ means that $a-b$ is nilpotent). Then, by Lemma~\ref{NilpotentDerivations} \[\frac\partial{\partial s}\widetilde W_{n,m,0}=\widetilde W_{n,m+1,0}+\widetilde W_{n,m+1,1}\equiv 0\] for $m\le l-3$. We deduce that $\widetilde W_{n,m,1}\equiv 0$ for $m\le l-3$. Iterating this process we obtain that all $\widetilde W_{n,m,p}\equiv 0$ when $m\le l-3$. Hence, by writing the derivative of $\widetilde W_{n,l-3,p}\equiv 0$ for any $p\in [0,l-3]$ we obtain $\widetilde W_{n,l-2,p}=-\widetilde W_{n,l-2,p+1}$. Together these $l-2$ equivalences provide 
\begin{equation}\label{flipequiv}
\widetilde W_{n,l-2,0}\equiv (-1)^{l-2}\widetilde W_{n,l-2,l-2}\text{ for }n=l-1,l.
\end{equation}

Now choose $J \subset [2,l]$. Multiplying $W_{l,l-2}$ by $Y_{J}(s)$ and using the relations \eqref{archibiboolean} we have:
\[Y_{J}(s)\sum_{I \subset [2,l]}(-1)^{|I|}Y_{I \cup \{1\}}(s)
\frac{\partial^{l-2} Y_{[2,l]\backslash I}(s)}{\partial s^{l-2}}=
Y_{J \cup \{1\}}(s)\sum_{I \subset [2,l]}(-1)^{|I|}Y_{I}(s)
\frac{\partial^{l-2} Y_{[2,l]\backslash I}(s)}{\partial s^{l-2}}.\]
In view of~\eqref{flipequiv} for $n=l-1$, we have
\begin{equation}\label{aftermult}
\sum_{I \subset [2,l]}(-1)^{|I|}Y_{I}(s)
\frac{\partial^{l-2} Y_{[2,l]\backslash I}(s)}{\partial s^{l-2}}=-\widetilde W_{l-1,l-2,0}+(-1)^{l-2}\widetilde W_{l-1,l-2,l-2}\equiv 0.
\end{equation}

We see that $Y_J(s)W_{l,l-2}\equiv0$ for any $Y_J(s)$ with $J\subset [2,l]$. By applying~\eqref{flipequiv} for $n=l$ we also have $Y_J(s)\widetilde W_{l,l-2,l-2}\equiv0$. However, it is evident from the relations~\eqref{archibiboolean} that $J^{\infty}(R(\mC_l))[[s]]$ admits an involution exchanging $Y_I^{(j)}$ and $Y_{[l]\backslash I}^{(j)}$. Applying this involution to $Y_J(s)\widetilde W_{l,l-2,l-2}$ we obtain $Y_{J'}(s)W_{l,l-2}$ where $1\in J'$ and we see that this product is also nilpotent. As a result, all $Y_J(s)W_{l,l-2}\equiv0$ and, since $W_{l,l-2}^2$ is a sum of multiples of various $Y_J(s)W_{l,l-2}$, we have $W_{l,l-2}\equiv 0$.\hfill\qedhere


\end{proof}

\begin{rem}\label{cube_rel}
Consider $r\le l$ and an affine embedding $\mC_r\hookrightarrow\mC_l$ mapping vertices to vertices. This induces a map $f:J^\infty(R(\mC_r))\to J^\infty(R(\mC_l))$ and any $f(W_{r,k})$ will lie in the nilradical of $J^\infty(R(\mC_l))$. It can be shown that all elements of this form generate the nilradical, more specifically, a minimal set of such relations is obtained in Subsection~\ref{parallel}.
\end{rem}

\section{Symmetric polynomials and dual spaces of graded components}\label{SymPol}

\subsection{Arcs over the polynomial ring}
In this subsection we use the notations from Subsection~\ref{initial} and discuss the arc ring $\mathbb J_n$ of a polynomial ring. Note that, similarly to $\mathbb J_n$, any arc ring $J^\infty(R)$ has a grading $\grad$ given by $\grad r^{(d)}=d$ for any $r\in R$. If $B$ is any  $\grad$-homogeneous subspace, quotient or subquotient of an arc ring, we will denote by $B^{(d)}$ its component spanned by $b$ with $\grad b=d$. In particular, if all $B^{(d)}$ are finite-dimensional (for example, when $B=\mathbb J_n[\bar a]$, $\bar a\in\mathbb N^n$) we will understand the dual space $B^*$ to be the graded dual with respect to $\grad$.


Choose $\bar a\in\bN^n$ and consider $a_1+ \dots +a_n$ formal variables 
\[s_{1}^{(1)}, \dots, s_{1}^{(a_1)},\dots, s_{n}^{(1)}, \dots, s_{n}^{(a_n)}.\]
Let $\epsilon_i \in \mathbb{N}^n$ be the vector with its $i$th coordinate equal to 1 and all other coordinates 0.
We have
\[X_i(s_{i}^{(j)})\in \mathbb J_n[\epsilon_i]\widehat\otimes \Bbbk[s_{i}^{(j)}]\]
where we view any $\mathbb J_n[\bar a]$ as a topological vector space with a neighborhood base at 0 consisting of the spaces $\bigoplus_{d\ge N}\mathbb J_n[\bar a]^{(d)}$ for $N\ge0$.
Consider the expression
\begin{multline}\label{UDefinition}
U_{\bar a}(X_1,\dots, X_n):=X_1(s_{1}^{(1)})\dots X_1(s_{a_1}^{(1)})\cdot \dots \cdot X_n(s_{1}^{(n)})\dots X_n(s_{a_n}^{(n)})\in\\
\mathbb J_n[\bar a] \widehat\otimes 
\Bbbk[s_{1}^{(1)}, \dots, s_{1}^{(a_1)},\dots, s_{n}^{(1)}, \dots, s_{n}^{(a_n)}].
\end{multline}
We have:
\[U_{\bar a}(X_1,\dots, X_n)=\sum_{(k_i^j|1 \leq i \leq n,1 \leq j \leq a_i)\in\bN^{a_1+\dots+a_n}}\prod_{1 \leq i \leq n,1 \leq j \leq a_i} (s_i^{(j)})^{k_{i}^j}\prod_{1 \leq i \leq n,1 \leq j \leq a_i} X_i^{(k_i^j)}.\]
Collecting all terms with the same $\prod_{1 \leq i \leq n,1 \leq j \leq a_i}
X_i^{(k_i^j)}$ we obtain:
\begin{multline*}
U_{\bar a}(X_1,\dots, X_n)=\\
\sum_{\substack{(k_i^j|1 \leq i \leq n,1 \leq j \leq a_i)\in\bN^{a_1+\dots+a_n}\text{ with }\\k_i^1\leq k_i^2 \leq \dots \leq k_i^{a_i}\text{ for every }i}}\left(\prod_{1 \leq i \leq n,1 \leq j \leq a_i}
X_i^{(k_i^j)}\left(\sum_{(\sigma_1, \dots, \sigma_n) \in \mathfrak{S}_{a_1}\times \dots \times \mathfrak{S}_{a_n}} \prod_{1 \leq i \leq n,1 \leq j \leq a_i}\left(s_i^{(\sigma_i(j))}\right)^{k_{i}^j}\right)\right).
\end{multline*}
Therefore the coefficients at $\prod_{1 \leq i \leq n,1 \leq j \leq a_i}
X_i^{(k_i^j)}$ form a linear basis of the ring
\begin{equation}\label{LambdaDefinition}\Lambda_{\bar a}(\{s_i^{(j)}\}_{1 \leq i \leq n,1 \leq j \leq a_i}):=\Bbbk[s_i^{(j)}]_{1 \leq i \leq n,1 \leq j \leq a_i}^{\mathfrak{S}_{a_1}\times \dots \times \mathfrak{S}_{a_n}}\end{equation}
of partially symmetric functions. Here the factor $\mathfrak{S}_{a_1}$ permutes the variables $s_1^{(1)},\dots,s_1^{(a_1)}$ and so on. We will shorten the above notation to $\Lambda_{\bar a}(\mathbf s)$ where we denote $\mathbf s=\{s_i^{(j)}\}_{1 \leq i \leq n,1 \leq j \leq a_i}$. 
In particular, we have:
\begin{equation}
U_{\bar a}(X_1,\dots, X_n)\in \mathbb J_n[\bar a] \widehat\otimes \Lambda_{\bar a}(\mathbf s).
\end{equation}

Consider a linear functional $\xi \in \mathbb J_n[\bar a]^*$ and let $1$ be the identity isomorphism on $\Lambda_{\bar a}(\mathbf s)$. We have the map $\xi\widehat\otimes 1:\mathbb J_n[\bar a]\widehat\otimes\Lambda_{\bar a}(\mathbf s)\to\Lambda_{\bar a}(\mathbf s)$ and denote
\begin{equation}\label{SumPolynomialforFunctional}
f_{\xi}:=(\xi\widehat\otimes 1)\left(U_{\bar a}(X_1, \dots, X_n)\right)\in \Lambda_{\bar a}(\mathbf s).
\end{equation}

Set $\mathbb J_n^*=\bigoplus_{\bar a \in \mathbb{N}^n}\mathbb J_n[\bar a]^*$. We immediately obtain
\begin{lem}\label{dualPolynomial}
The map $\gamma_{\bar a}:\xi \mapsto f_\xi$ is a linear isomorphism
from $\mathbb J_n[\bar a]^*$ to $\Lambda_{\bar a}(\mathbf s)$. The map $\gamma=\bigoplus_{\bar a \in \mathbb{N}^n}\gamma_{\bar a}$ is an isomorphism  from $\mathbb J_n^*$ to $\bigoplus_{\bar a \in \mathbb{N}^n} \Lambda_{\bar a}(\mathbf s)$.
\end{lem}

For vectors $\bar a, \bar b \in \mathbb{N}^n$ we have a natural map
\begin{equation}\label{comult}
\Delta_{a,b}:\Lambda_{\bar a+\bar b}(\mathbf s)\hookrightarrow
\Lambda_{\bar a}(\mathbf s)\otimes\Lambda_{\bar b}(\{s_i^{(j)}\}_{1 \leq i \leq n,1+a_i \leq j \leq b_i+a_i})
\simeq 
\Lambda_{\bar a}(\mathbf s)\otimes\Lambda_{\bar b}(\mathbf s)
\end{equation}
because any polynomial symmetric in $s_i^{(1)}, \dots, s_i^{(a_i)},s_i^{(a_i+1)},\dots,s_i^{(a_i+b_i)}$ is symmetric separately in
$s_i^{(1)}, \dots, s_i^{(a_i)}$ and $s_i^{(a_i+1)},\dots,s_i^{(a_i+b_i)}$.
We choose a concrete isomorphism on the right in~\eqref{comult} by identifying the variable $s_i^{(a_i+j)}$ in $\Lambda_{\bar b}(\{s_i^{(j)}\}_{1 \leq i \leq n,1+a_i \leq j \leq b_i+a_i})$ with the variable $s_i^{(j)}$ in $\Lambda_{\bar b}(\mathbf s)$. Then~\eqref{comult} provides a comultiplication map $\Delta$ on the space  
$\bigoplus_{\bar a \in \mathbb{N}^n}  \Lambda_{\bar a}(\mathbf s)$. This comultiplication is obviously coassociative.

\begin{prop}\label{MultiplicationDualtoComultiplication}
The comultiplication map $\Delta$ is dual to the multiplication in the ring $\mathbb J_n$ with respect to $\gamma$.
\end{prop}
\begin{proof}
Consider the comultiplication $\widetilde\Delta=(\gamma^{-1}\otimes\gamma^{-1})\Delta\gamma$ on $\mathbb J_n^*$. We need to show that for any $X,Y\in\mathbb J_n$ and $\xi\in \mathbb J_n^*$ we have $\xi(XY)=\widetilde\Delta(\xi)(X\otimes Y)$. By linearity it suffices to consider the case when $X$ and $Y$ are monomials and $\xi(Z)=1$ for some monomial $Z$ while $\xi$ vanishes on all other monomials. We are to check that $\widetilde\Delta(\xi)(X\otimes Y)$ equals 1 if $XY=Z$ and 0 otherwise. This follows directly from the definitions of $\gamma$ and $\Delta$.
\end{proof}

In the following corollary and below we often identify $\mathbb J_n[a]^*$ with $\Lambda_{\bar a}(\mathbf s)$ via $\gamma_{\bar a}$.
\begin{cor}
\hfill
\begin{enumerate}[label=(\alph*)]
\item Consider an ideal $I \subset \mathbb J_n$.
Then the subspace in $\bigoplus_{\bar a \in \mathbb{N}^n}  \Lambda_{\bar a}(\mathbf s)$ which annihilates $I$ is 
a subcoalgebra.
\item Consider a subalgebra $S \subset \mathbb J_n$.
Then the subspace in $\bigoplus_{\bar a \in \mathbb{N}^n}  \Lambda_{\bar a}(\mathbf s)$ which annihilates $S$ is 
a coideal.
\end{enumerate} 
\end{cor}

\subsection{Action of the current Lie algebra on dual components}

Since the polynomial algebra $R=\Bbbk[X_1, \dots, X_n]$ is $\mathbb{N}^n$-graded, it is acted upon by an $n$-dimensional abelian Lie algebra $\fh=\langle h_1, \dots, h_n \rangle$. For $r\in R[\bar a]$ we have $h_i(p)=a_i p$. The Lie algebra $\fh$ acts by derivations and, applying Defintion~\ref{currentDerivationDefinition}, we obtain derivations $h_i^{(k)}$ on $\mathbb J_n$ given by \[h_{i_1}^{(k)}(X_{i_2}^{(j)})=\delta_{i_1,i_2} X_{i_2}^{(j-k)}.\] By Lemma~\ref{currentaction} we have an action of $\fh[s]$ on $\mathbb J_n$ with $h_is^k$ acting as $h_i^{(k)}$. This $\fh[s]$-action is faithful, hence $h_i s^k\mapsto h_i^{(k)}$ provides an embedding $\fh[s]\hookrightarrow \Der\mathbb J_n$ and allows us to identify $h_is^k$ with $h_i^{(k)}$.
This action also preserves the $\bN^n$-grading, i.e.\
\[\mathfrak{h}[s](\mathbb J_n[\bar a])\subset
\mathbb J_n[\bar a].\]

We may extend the action of $\fh[s]$ to $\mathbb J_n[[s]]$ (coefficientwise) and $\mathbb J_n \widehat\otimes \Lambda_{\bar a}(\mathbf s)$ (by acting on the first factor). In particular, we have
\begin{equation}\label{hCurrentGenerating}
h_{i_1}^{(k)}X_{i_2}(s)=\delta_{i_1,i_2}s^k X_{i_2}(s).
\end{equation}
\begin{lem}
\[h_i^{(k)}(U_{\bar a}(X_1, \dots, X_n))=U_{\bar a}(X_1, \dots, X_n)\sum_{j=1}^{a_i}
\left(s_{i}^{(j)}\right)^k.\]
\end{lem}
\begin{proof}
This follows from \eqref{hCurrentGenerating} by the Leibniz rule.
\end{proof}
\begin{prop}
The action of $\fh[s]$ on $\Lambda_{\bar a}(\mathbf s)$ dual (with respect to $\gamma_{\bar a}$) to the $\fh[s]$-action on $\mathbb J_n[\bar a]$ is given by
\[h_i^{(k)}(f)=-\sum_{j=1}^{a_i}
\left(s_{i}^{(j)}\right)^k f\]
for any $f \in \Lambda_{\bar a}(\mathbf s)$.
\end{prop}
\begin{proof}
Choose $f\in \Lambda_{\bar a}(\mathbf s)$ and let $\xi=\gamma_{\bar a}^{-1}(f)$. Then for the dual $\fh[s]$-action we have:
\begin{multline}\label{ProofDualAction}
h_i^{(k)}(f)=(h_i^{(k)}(\xi)\widehat \otimes 1)(U_{\bar a}(X_1, \dots, X_n))=
-(\xi \widehat\otimes 1)\big(h_i^{(k)}(U_{\bar a}(X_1, \dots, X_n))\big)=\\
-(\xi \widehat\otimes 1)(U_{\bar a}(X_1, \dots, X_n))\sum_{j=1}^{a_i} \left(s_{i}^{(j)}\right)^k=-\sum_{j=1}^{a_i}
\left(s_{i}^{(j)}\right)^k f.\qedhere
\end{multline}
\end{proof}
The action \eqref{LieDerivationsAction} of the Lie algebra $\Der^c \Bbbk[[s]]$ also preserves graded components of the arc ring $\mathbb J_n$.
\begin{prop}
The action of $\Der^c \Bbbk[[s]]$ on $\Lambda_{\bar a}(\mathbf s)$ dual to the action on $\mathbb J_n[\bar a]$ is given by
\[d_k(f)=-\left(\sum_{1 \leq i \leq n, 1 \leq j \leq a_i}
\left(s_i^{(j)}\right)^{k+1} 
\frac{\partial}{\partial s_i^{(j)}}\right)f.\]
\end{prop}
\begin{proof}
We have
\[d_k\left(X_i\left(s_i^{(j)}\right)\right)=\left(s_i^{(j)}\right)^{k+1} \frac{\partial}{\partial s_i^{(j)}} X_i\left(s_i^{(j)}\right).\]
Therefore, by the Leibniz rule:
\[d_k(U_{\bar a}(X_1, \dots, X_n))=\left(\sum_{1 \leq i \leq n, 1 \leq j \leq a_i}
\left(s_i^{(j)}\right)^{k+1} 
\frac{\partial}{\partial s_i^{(j)}}\right) U_{\bar a}(X_1, \dots, X_n).\]
The proof is now completed similarly to the previous proposition.
\end{proof}

\subsection{Arcs and quadratic monomial ideals} \label{subsection arcs of monomial ideals}
Consider a set of pairs \[\pi=\{(\alpha_1,\beta_1), \dots, (\alpha_p,\beta_p)\},~
\alpha_i, \beta_i \in [n],~ \alpha_i<\beta_i\] and the monomial ideal  $\mathcal{M}_\pi=\langle X_{\alpha_i}X_{\beta_i} \rangle_{i=1, \dots,p}$ in $\Bbbk[X_1, \dots, X_n]$. The goal of this subsection is to
describe the space dual to the ring $J^\infty(\Bbbk[X_1, \dots, X_n]/\mathcal{M}_\pi)$ in terms of symmetric functions.

We have the surjection $\rho:\mathbb J_n\twoheadrightarrow J^\infty(\Bbbk[X_1, \dots, X_n]/\mathcal{M}_\pi)$ with kernel generated by coefficients of the series
\begin{equation}\label{JetforQuadraticMonomialIdeal}
X_{\alpha_i}(s)X_{\beta_i}(s),~(\alpha_i,\beta_i) \in \pi.
\end{equation} 

This kernel is $\mathbb{N}^n$-homogeneous, hence, the ring $J^\infty(\Bbbk[X_1, \dots, X_n]/\mathcal{M}_\pi)$ is $\mathbb{N}^n$-graded and $\rho$ restricts to $\mathbb{N}^n$-homogeneous surjections $\rho[\bar a]$ from $\mathbb J_n[\bar a]$ to $J^\infty(\Bbbk[X_1, \dots, X_n]/\mathcal{M}_\pi)[\bar a]$.
Therefore, the (graded) dual of the homogeneous component 
$J^\infty(\Bbbk[X_1, \dots, X_n]/\mathcal{M}_\pi)[\bar a]$
can be naturally viewed as a subspace in $\mathbb J_n[\bar a]^*$.

The kernel of $\rho$ is $\mathfrak{h}[s]$ invariant. Therefore, the map $\rho$ is a homomorphism of  $\mathfrak{h}[s]$-modules and thus 
the dual of the homogeneous component 
$J^\infty(\Bbbk[X_1, \dots, X_n]/\mathcal{M}_\pi)[\bar a]$ has a structure of a $\mathfrak{h}[s]$-submodule.

\begin{prop}\label{DualSpaceQuagraticMon}
The subspace \[\gamma_{\bar a}\left(J^\infty(\Bbbk[X_1, \dots, X_n]/\mathcal{M}_\pi)[\bar a]^*\right)\subset \Lambda_{\bar a}(\mathbf s)\] is the principal ideal generated by \begin{equation}\label{principal}
\prod_{\substack{(\alpha, \beta) \in \pi\\ 1 \leq i \leq a_{\alpha},1 \leq j \leq a_{\beta}}}
(s_{\alpha}^{(i)}-s_{\beta}^{(j)}).
\end{equation}
This principal ideal is a cyclic  $\mathfrak{h}[s]$-module generated by \eqref{principal}
that is isomorphic to $\Lambda_{\bar a}(\mathbf s)$. 
\end{prop}
\begin{proof}
In view of \eqref{JetforQuadraticMonomialIdeal}, the kernel of $\rho[\bar a]$ is the linear span of coefficients of the series
\begin{equation}\label{QuadraticMonomialKernel}
U_{\bar a -e_\alpha -e_\beta}(X_1, \dots, X_n)X_\alpha(t)X_\beta(t)=U_{\bar a}(X_1, \dots, X_n)|_{s_{\alpha}^{(a_\alpha)}\mapsto t,s_{\beta}^{(a_\beta)}\mapsto t}
\end{equation}
for $(\alpha,\beta) \in \pi$,
where $\epsilon_\alpha$ is the $\alpha$th basis vector in $\mathbb{N}^n$.
For $\xi \in \mathbb J_n[\bar a]^*$
we have
\begin{multline*}
(\xi\widehat\otimes 1)\left(U_{\bar a}(X_1, \dots, X_n)|_{s_{\alpha}^{(a_\alpha)}\mapsto t,s_{\beta}^{(a_\beta)}\mapsto t}\right)=\\
(\xi\widehat\otimes 1) \left(U(X_1, \dots, X_n)_{\bar a}\right)|_{s_{\alpha}^{(a_\alpha)}\mapsto t,s_{\beta}^{(a_\beta)}\mapsto t}=\gamma_{\bar a}(\xi)|_{s_{\alpha}^{(a_\alpha)}\mapsto t,s_{\beta}^{(a_\beta)}\mapsto t}.
\end{multline*}
Thus, $\xi(\ker\rho[\bar a])=0$ if and only if $\gamma_{\bar a}(\xi)|_{s_{\alpha}^{(a_\alpha)}\mapsto t,s_{\beta}^{(a_\beta)}\mapsto t}=0$. By symmetricity, $\gamma_{\bar a}(\xi)$ must then be divisible  by 
\eqref{principal}.
The last claim follows from \eqref{ProofDualAction}.
This completes the proof.
\end{proof}

We also need the following generalization of this proposition.
For each pair of indices $(\alpha,\beta), 1 \leq \alpha <\beta \leq n$ choose a number $k(\alpha,\beta)\in \mathbb{N}$. Denote the collection of these numbers by $\bar k$.
Consider the ideal $\mathcal{I}_{\bar k}$ generated by the coefficients of
\begin{equation}\label{QuadraticJetMonomialEquations}
\Big\langle  X_{\alpha_i}(s)\frac{\partial^{k'}X_{ \beta_i}(s)}{\partial s^{k'}} \Big\rangle ,\,1 \leq \alpha <\beta \leq n,\, 0 \leq k' <k(\alpha,\beta).
\end{equation}
We denote by $\rho$ the natural surjection $\mathbb J_n \twoheadrightarrow \mathbb J_n/\mathcal{I}_{\bar k}$ and let $\rho[\bar a]:\mathbb J_n[\bar a]\twoheadrightarrow(\mathbb J_n/\mathcal{I}_{\bar k})[\bar a]$ be its $\mathbb N^n$-homogeneous components. As before $\mathcal{I}_{\bar k}$ is $\mathfrak{h}[s]$-invariant. Therefore, $\rho[\bar a]$ is a map of $\mathfrak{h}[s]$-modules.
As in the previous case we identify $(\mathbb J_n/\mathcal{I}_{\bar k})[\bar a]^*$ with a principal ideal in $\Lambda_{\bar a}(\mathbf s)$.

\begin{prop}\label{DualSpaceDerivedQuagraticMon}
The subspace $\gamma_{\bar a}\left((\mathbb J_n/\mathcal{I}_{\bar k})[\bar a]^*\right)\subset \Lambda_{\bar a}(\mathbf s)$ is the principal ideal generated by
\begin{equation}\label{genprincipal}\prod_{\substack{1\le \alpha< \beta\le n\\1 \leq i \leq a_{\alpha},1 \leq j \leq a_{\beta}}}
(s_{\alpha}^{(i)}-s_{\beta}^{(j)})^{k(\alpha,\beta)}.
\end{equation}
This principal ideal is a cyclic  $\mathfrak{h}[s]$-module generated by \eqref{genprincipal}
that is isomorphic to $\Lambda_{\bar a}(\mathbf s)$. 
\end{prop}
\begin{proof}
By definition, $\ker\rho[\bar a]$ is the linear span of coefficients of the series
\begin{equation}\label{QuadraticDerivedMonomialKernel}
U_{\bar a -e_\alpha -e_\beta}(X_1, \dots, X_n)X_\alpha(t)\frac{\partial ^{k'}X_\beta(t)}{\partial t^{k'}}=\left. 
\frac{\partial ^{k'} U_{\bar a}(X_1, \dots, X_n)}{(\partial s_{\beta}^{(a_\beta)})^{k'}}
\right|_{s_{\alpha}^{(a_\alpha)}\mapsto t,s_{\beta}^{(a_\beta)}\mapsto t}
\end{equation}
where $0 \leq k' <k$.
We have for any $\xi \in \mathbb J_n[\bar a]^*$:
\begin{multline*}
(\xi\widehat\otimes 1)\left(\left.
\frac{\partial ^{k'} U_{\bar a}(X_1, \dots, X_n)}{(\partial s_{\beta}^{(a_\beta)})^{k'}}
\right|_{s_{\alpha}^{(a_\alpha)}\mapsto t,s_{\beta}^{(a_\beta)}\mapsto t}\right )=\\\left.
\frac{\partial ^{k'}\xi(U_{\bar a}(X_1, \dots, X_n))}{(\partial s_{\beta}^{(a_\beta)})^{k'}}
\right|_{s_{\alpha}^{(a_\alpha)}\mapsto t,s_{\beta}^{(a_\beta)}\mapsto t}=\left.\frac{\partial ^{k'}\gamma_{\bar a}(\xi)}{(\partial s_{\beta}^{(a_\beta)})^{k'}}\right|_{s_{\alpha}^{(a_\alpha)}\mapsto t,s_{\beta}^{(a_\beta)}\mapsto t}.
\end{multline*}
Thus, $\xi(\ker\rho[\bar a])=0$ if and only if for any $\alpha<\beta$ and $0\le k'<k(\alpha, \beta)$ we have
\[ \left.\frac{\partial ^{k'}\gamma_{\bar a}(\xi)} {(\partial s_{\beta}^{(a_\beta)})^{k'}}\right|_{s_{\alpha}^{(a_\alpha)}\mapsto t,s_{\beta}^{(a_\beta)}\mapsto t}=0.\]
By symmetricity, such $\gamma_{\bar a}(\xi)$ are precisely those divisible by $(s_{\alpha}^{(i)}-s_{\beta}^{(j)})^{k(\alpha,\beta)}$ for any $\alpha<\beta$, $1 \leq i \leq a_{\alpha}$ and  $1 \leq j \leq a_{\beta}$. The last claim is straightforward. 
\end{proof}



\begin{rem}
Consider a finite set $S\subset\bN^n$ and the ideal $I_S \subset \Bbbk[X_1, \dots, X_n]$ generated by monomials $X^{\bar b}$, $\bar b \in S$. 
Then the subspace
$\gamma_{\bar a}\left(J^\infty(\Bbbk[X_1, \dots, X_n]/I_S)[\bar a]^*\right)\subset \Lambda_{\bar a}(\mathbf s)$ is the ideal
\[
\left\{f\text{ such that }
f\left|_{s_1^{(1)},\dots,s_1^{(b_1)}, \dots, s_n^{(1)},\dots,s_n^{(b_n)}\mapsto t}=0\right. \text{ for all } \bar b \in S  \right\}.
\]
This generalizes Propositions~\ref{DualSpaceQuagraticMon} and~\ref{DualSpaceDerivedQuagraticMon} the proof being similar. We note that in this general case the dual space is still an ideal in $\Lambda_{\bar a}(\mathbf s)$ which, however, need not be principal.
\end{rem}

\section{Arc rings of toric rings} 
\label{ToricArcs}




\subsection{Graded components of toric rings}

Consider a normal lattice polytope $P \subset \mathbb{R}_{\ge0}^n$ and the corresponding toric ring $R(P)$. Note that $R(P)$ is invariant under shifts of $P$ so we do not lose generality by assuming that $P$ lies in the positive orthant.

We have an embedding
\[ \eta: R(P) \hookrightarrow\Bbbk[z_1, \dots, z_n,w]\] with $\eta(Y_{\bar \alpha})=z^{\bar \alpha}w$ for $\bar\alpha \in P \cap \mathbb{Z}^n$.
We study the corresponding map of arc rings
\[ J^{\infty}(\eta): J^{\infty}R(P) \to J^{\infty}(\Bbbk[z_1, \dots, z_n,w]).\] By definition, applying $J^{\infty}(\eta)$ coefficientwise we have $J^{\infty}(\eta)(Y_{\bar \alpha}(s))=z(s)^{\bar \alpha}w(s)$.
Let $\{\bar\alpha^1, \dots, \bar\alpha^m\}$ be the set of integer points in $P$. 
We have a surjection \[\tau:\Bbbk[X_{\bar \alpha^1},\dots,X_{\bar \alpha^m}]\twoheadrightarrow R(P)\] where $\tau(X_{\bar \alpha^i})=Y_{\bar \alpha^i}$. We also have the map
\[\eta \circ \tau: \Bbbk[X_{\bar \alpha^1}, \dots, X_{\bar \alpha^m}]\rightarrow \Bbbk[z_1, \dots, z_n,w]\] and for a vector $\bar r =(r_{i})\in \mathbb{N}^m$:
\[\eta \circ \tau(X^{\bar r})=
z^{\sum_{i=1}^m (r_i \bar \alpha^i)}w^{\sum_{i=1}^m r_i}.\]
Thus, $\eta \circ \tau$ maps monomials in the $X_{\bar \alpha^i}$ to monomials in the $z_i$ and $w$. For $\bar a \in \mathbb{N}^n,~ L \in \mathbb{N}$ we define: 
\[\mathcal{R}(\bar a, L)=\left\{\bar r\in\mathbb N^m\big|\sum_{i=1}^m r_i \bar\alpha^i =\bar a\text{ and } \sum_{i=1}^m r_i =L\right\}.\]
Then 
$\eta \circ \tau (X^{\bar r})=z^{\bar a} w^L$ if and only if
$\bar r \in  \mathcal{R}(\bar a, L).$

The map $J^{\infty}(\tau):\mathbb J_m\to J^\infty R(P)$ is surjective, since $J^{\infty}(\tau)(X_{\bar \alpha^i}^{(j)})=Y_{\bar \alpha^i}^{(j)}$ (more generally, $J^\infty$ is easily seen to preserve epimorphisms). Therefore, we have:
\[J^{\infty}(\eta)(J^\infty R(P))= J^{\infty}(\eta \circ \tau)(\mathbb J_m)\subset J^\infty(\Bbbk[z_1, \dots, z_n,w])\]
and we have a surjective map
\begin{equation}\label{nu}\nu:J^{\infty}_\mathrm{red}(R(P))\twoheadrightarrow J^{\infty}(\eta)(J^\infty R(P))= J^{\infty}(\eta \circ \tau)(\mathbb J_m).\end{equation}
By Corollary \ref{jredinj} $\nu$ is an isomorphism.

For the rest of this section we fix a monomial order $\prec$ on $\Bbbk[X_{\bar \alpha^1},\dots,X_{\bar \alpha^m}]$ as well as the corresponding order on $\bN^m$.
As discussed in Subsection~\ref{initial}, we have two $(\bN^m,\prec)$-filtrations with components $\mathbb J_m[\prec\bar r]$ and $\mathbb J_m[\preceq\bar r]$. 
Since both $J^{\infty}_\mathrm{red}(R(P))$ and $J^{\infty}(\eta \circ \tau)(\mathbb J_m)$ are surjective images of $\mathbb J_m$ they admit induced $(\mathbb N^m,\prec)$-filtrations which are seen to be identified by the isomorphism $\nu$.
Since $\nu$ respects the grading $\grad$, we obtain
\begin{lem}\label{redlowerbound}
For any $d\in\mathbb N$ and $\bar r \in \mathbb N^m$ the map $\nu$ induces an isomorphism of vector spaces 
\[\gr_{\prec}J^{\infty}_\mathrm{red}R(P)[\bar r]^{(d)}
\simeq
\gr_{\prec}J^{\infty}(\eta \circ \tau)(\mathbb J_m)[\bar r]^{(d)}.\]
\end{lem}

\subsection{Dual of the inclusion map}\label{DualToInclusion}




Retaining notations from the previous subsection consider $\bar r\in \mathcal{R}(\bar a, L)\subset \mathbb N^m$ for some $\bar a\in\bN^n$ and $L\in\bN$. Consider the subspace $A_{\bar r} \subset J^{\infty}(\Bbbk[z_1, \dots, z_n,w])$ spanned by elements of the form
\[
\eta(Y_{\bar \alpha^1})^{(j_1^1)}\dots \eta(Y_{\bar \alpha^1})^{(j_{r_1}^1)} \dots
\eta(Y_{\bar \alpha^m})^{(j_1^m)}\dots \eta(Y_{\bar \alpha^m})^{(j_{r_m}^m)}=J^{\infty}(\eta)\left(\prod_{i=1}^m\prod_{k=1}^{r_m} Y_{\bar\alpha^i}^{(j^i_k)}\right)
\]
for all $j_u^v \in \bN$.
This subspace coincides with $J^{\infty}(\eta \circ \tau)(\mathbb J_m[\bar r])$.

In this subsection we study the dual map to the inclusion
\begin{equation}\label{componentInjectionDefinition}
\varphi_{\bar r}:A_{\bar r} \hookrightarrow J^{\infty}\left(\Bbbk[z_1, \dots, z_n,w]\right)\left[\bar a, L\right]
\end{equation}
where on the right we have a component of the natural $\mathbb N^{n+1}$-grading.

We consider the following series in variables $t_i^{(j)}$ with $1\le i\le m$ and $1\le j\le r_i$:
\begin{multline}\label{VDefinition}
V_{\bar r}:= \eta(Y_{\bar\alpha^1})(t_1^{(1)})\dots \eta(Y_{\bar\alpha^1})(t_1^{(r_1)}) \dots
\eta(Y_{\bar\alpha^m})(t_m^{(1)})\dots \eta(Y_{\bar\alpha^m})(t_m^{(r_m)})\\=
(z^{\bar \alpha^1}w)(t_1^{(1)})\dots (z^{\bar \alpha^1}w)(t_1^{(r_1)}) 
\dots
(z^{\bar \alpha^m}w)(t_m^{(1)})\dots (z^{\bar \alpha^m}w)(t_m^{(r_m)}).
\end{multline}
From the definitions we have
\begin{prop}\label{ASpan}
The space $A_{\bar r}$ is spanned by the coefficients of $V_{\bar r}$ as of a series in the $t_i^{(j)}$.
\end{prop}

By Lemma \ref{dualPolynomial} we obtain a linear isomorphism
\[\gamma_{\bar a,L}:J^\infty(\Bbbk[z_1,\dots,z_n,w])[\bar a,L ]^*
\to  \Lambda_{\bar a}(\mathbf s)\otimes \Lambda_L(\{s_w^{(j)}\}_{j=1, \dots, L}).\] We will denote the right-hand side above by $\Lambda_{\bar a,L}(\mathbf s)$.

Recall that
$\bar a=\sum_{j=1}^m r_j \bar \alpha^j$ and
$L=\sum_{i=1}^m r_i$. Consider the homomorphism of polynomial algebras
\[
\varphi_{\bar r}^{\vee}:\Bbbk[s_i^{(j)},s_w^{(j')}]_{\substack{1 \leq i \leq n\\ 1 \leq j \leq a_i\\ 1 \leq j' \leq L}}\rightarrow 
\Bbbk[t_i^{(j)}]_{\substack{1 \leq i \leq m\\ 1 \leq j \leq r_i}}
\]
given by
\begin{equation}\label{SymmetricAlgebrasMapDefinition}
\varphi_{\bar r}^{\vee}(s_i^{(j)})=\varphi_{\bar r}^{\vee}(s_w^{(j')})=
\begin{cases}
t_1^{(1)},\text{ if }j\in[1,\alpha^1_i]\text{ and }j'=1,\\
\dots\\
t_1^{(r_1)},\text{ if }j\in[(r_1-1)\alpha^1_i+1,r_1\alpha^1_i]\text{ and }j'=r_1,\\
t_2^{(1)},\text{ if }j\in[r_1\alpha^1_i+1,r_1\alpha^1_i+\alpha^2_i]\text{ and }j'=r_1+1,\\
\dots\\
t_2^{(r_2)},\text{ if }j\in[r_1\alpha^1_i+(r_2-1)\alpha^2_i+1,r_1\alpha^1_i+r_2\alpha^2_i]\text{ and }j'=r_1+r_2,\\
\dots
\end{cases}
\end{equation}
In other words, for a given $1\le i\le n$ we divide the tuple $s_i^{(1)},\dots,s_i^{(a_i)}$ into $L$ groups of consecutive elements so that the first $r_1$ groups contain $\alpha^1_i$ elements, the next $r_2$ groups contain $\alpha^2_i$ elements, etc. Then $\varphi_{\bar r}^{\vee}$ maps the elements in the first $r_1$ groups to, respectively, $t_1^{(1)},\dots,t_1^{(r_1)}$, the elements in the next $r_2$ groups to, respectively, $t_2^{(1)},\dots,t_2^{(r_2)}$ and so on. As for the variables $s_w^{(1)},\dots,s_w^{(L)}$, the first $r_1$ of them  are mapped to, respectively, $t_1^{(1)},\dots,t_1^{(r_1)}$, the next $r_2$ are mapped to, respectively, $t_2^{(1)},\dots,t_2^{(r_2)}$, etc.



\begin{rem}\label{PowerSumsMaps}
The map $\varphi_{\bar r}^\vee$ can be described in the following way.
The ring $\Lambda_{\bar a,L}(\mathbf s)$ is generated by power sums
\[p_{k}\left(s_i^{(1)},\dots, s_i^{(a_i)} \right)=\sum_{j=1}^{a_i} \left(s_i^{(j)}\right)^k\text{ and }
p_{k}\left(s_w^{(1)},\dots, s_w^{(L)} \right)=\sum_{j=1}^{L} \left(s_w^{(j)}\right)^k.\]
We have
\begin{align*}\varphi_{\bar r}^\vee\left(p_{k}\left(s_i^{(1)},\dots, s_i^{(a_i)}\right)\right)&=\sum_{j=1}^m \alpha^{j}_i p_{k}\left(t_j^{(1)},\dots, t_j^{(r_j)}\right), \\
\varphi_{\bar r}^\vee\left(p_{k}\left(s_w^{(1)},\dots, s_w^{(L)}\right)\right)&=\sum_{j=1}^m  p_{k}\left(t_j^{(1)},\dots, t_j^{(r_j)}\right).
\end{align*}
In particular, $\varphi_{\bar r}^\vee(\Lambda_{\bar a,L}(\mathbf s))$ is generated as a ring by the elements 
$\sum_{j=1}^m \alpha^{j}_i p_{k}\left(t_j^{(1)},\dots, t_j^{(r_j)}\right)$, $i=1,\dots,n$ and
$\sum_{j=1}^m  p_{k}\left(t_j^{(1)},\dots, t_j^{(r_j)}\right)$.
\end{rem}

\begin{example}
Let $P$ be the convex hull of $\{(0,0),(1,2),(2,1)\}$ and let $\alpha_1=(0,0),\alpha_2=(1,1),\alpha_3=(2,1),\alpha_4=(1,2)$. Consider $\bar r=(1,1,1,1)\in \mathcal{R}((4,4),4)$.
Then $\varphi^{\vee}_{\bar r}$ is given by
\begin{align*}
s_w^{(1)}&\mapsto t_1^{(1)},\\
s_w^{(2)},s_1^{(1)},s_2^{(1)}&\mapsto t_2^{(1)},\\
s_w^{(3)},s_1^{(2)},s_1^{(3)},s_2^{(2)}&\mapsto t_3^{(1)},\\
s_w^{(4)},s_1^{(4)},s_2^{(3)},s_2^{(4)}&\mapsto t_4^{(1)}\\
\end{align*}
and we have:
\begin{align*}
\varphi_{\bar r}^\vee\left(p_k(s_w^{(1)},s_w^{(2)},s_w^{(3)},s_w^{(4)})\right)&=p_k(t_1^{(1)})+p_k(t_2^{(1)})+p_k(t_3^{(1)})+p_k(t_4^{(1)}),\\
\varphi_{\bar r}^\vee\left(p_k(s_1^{(1)},s_1^{(2)},s_1^{(3)},s_1^{(4)})\right)&=p_k(t_2^{(1)})+2p_k(t_3^{(1)})+p_k(t_4^{(1)}),\\
\varphi_{\bar r}^\vee\left(p_k(s_2^{(1)},s_2^{(2)},s_2^{(3)},s_2^{(4)})\right)&=p_k(t_2^{(1)})+p_k(t_3^{(1)})+2p_k(t_4^{(1)}).
\end{align*}
\end{example}

\begin{lem}\label{ImageIsSymmetric}
The image $\varphi_{\bar r}^{\vee}(\Lambda_{\bar a,L}(\mathbf s))$ is contained in the subring $\Lambda_{\bar r}(\mathbf t)$.
\end{lem}
\begin{proof}
The ring $\Lambda_{\bar a,L}(\mathbf s)$ decomposes into a product of its subspaces as
\begin{equation*}
\Lambda_{\bar a,L}(\mathbf s) =
\Lambda_{a_1}(\{s_1^{(j)}\}_{j=1, \dots, a_1})\dots \Lambda_{a_n}(\{s_n^{(j)}\}_{j=1, \dots, a_n})
\Lambda_L(\{s_w^{(j)}\}_{j=1, \dots, L}).
\end{equation*}
The image under $\varphi_{\bar r}^{\vee}$ of a polynomial in $\Lambda_{a_i}(\{s_i^{(j)}\}_{j=1, \dots, a_i})$ is evidently symmetric in $t_l^{(1)},\dots,t_l^{(r_l)}$ for each $l$. The same holds for any polynomial in the image of $\Lambda_L(\{s_w^{(j)}\}_{j=1, \dots, L})$ and the lemma follows.
\end{proof}

Recall the inclusion map  $\varphi_{\bar r}$ defined in \eqref{componentInjectionDefinition}.
We denote by \[\varphi_{\bar r}^*:J^{\infty}(\Bbbk[z_1, \dots, z_n,w])[\bar a, L]^*\to A_{\bar r}^*\] the (graded) dual map for this inclusion.

\begin{thm}\label{inclusiondual}
The annihilator of $A_{\bar r}$ in $J^{\infty}(\Bbbk[z_1, \dots, z_n,w])[\bar a, L]^*$ is the subspace $\ker(\varphi_{\bar r}^{\vee}\gamma_{\bar a,L})$. In other words, there exists a linear isomorphism \[\varepsilon_{\bar r}:A_{\bar r}^*\xrightarrow{\sim} \varphi_{\bar r}^{\vee}(\Lambda_{\bar a,L}(\mathbf s))\] such that $\varepsilon_{\bar r}\varphi_{\bar r}^*=\varphi_{\bar r}^{\vee}\gamma_{\bar a,L}$.
\end{thm}
\begin{proof}
The two claims are equivalent because the annihilator of $A_{\bar r}$ is $\ker\varphi_r^*$, i.e.\ both state that $\varphi_r^*$ and $\varphi_{\bar r}^{\vee}\gamma_{\bar a,L}$ have the same kernel.

Recall the expression $U_{\bar a}(z_1, \dots, z_n)$ defined in  \eqref{UDefinition} and consider \[U=U_{\bar a}(z_1, \dots, z_n)w(s_w^{(1)})\dots w(s_w^{(L)})\in J^{\infty}(\Bbbk[z_1, \dots, z_n,w])[\bar a, L]\widehat\otimes\Lambda_{\bar a, L}(\mathbf s).\] By definition we have $(1 \widehat\otimes \varphi_{\bar r}^{\vee})(U)=V_{\bar r}$.

Consider a functional $\xi \in  J^{\infty}(\Bbbk[z_1, \dots, z_n,w])[\bar a, L]^*$. We are to show that $\xi(A_{\bar r})=0$ if and only if $\varphi_{\bar r}^{\vee}\gamma_{\bar a,L}(\xi)=0$. Note that we may view $V_{\bar r}$ as an element of 
\[J^{\infty}(\Bbbk[z_1, \dots, z_n,w])[\bar a, L]\widehat\otimes\Lambda_{\bar r}(\mathbf t)\] 
and that, in view of Proposition~\ref{ASpan}, $\xi(A_{\bar r})=0$ if and only if $(\xi\widehat\otimes 1)(V_{\bar r})=0$. However,
\[(\xi\widehat\otimes 1)(V_{\bar r})=(\xi \widehat\otimes 1)(1 \widehat\otimes \varphi_{\bar r}^{\vee})(U)=(1 \widehat\otimes \varphi_{\bar r}^{\vee})(\xi \widehat\otimes 1)(U)=\varphi_{\bar r}^{\vee}(\gamma_{\bar a,L}(\xi)).\qedhere\]




\end{proof}

Recall that the monomial order $\prec$ induces a filtration on and an associated graded for $J^\infty(\eta\circ\tau)(\mathbb J_m)$. We can write the filtration components as \[J^\infty(\eta\circ\tau)(\mathbb J_m)[\preceq\bar r]=\sum_{\bar r'\preceq \bar r}A_{\bar r}\] and similarly for $J^\infty(\eta\circ\tau)(\mathbb J_m)[\prec\bar r]$. This lets us write the components of the associated graded as \[\gr_\prec J^\infty(\eta\circ\tau)(\mathbb J_m)[\bar r]=A_{\bar r}\left/\left(A_{\bar r}\cap\sum_{\bar r'\prec \bar r,\bar r'\in\mathcal R(\bar a,L)}A_{\bar r'}\right)\right..\] We can now apply Theorem~\ref{inclusiondual} to identify the (graded) dual space $\gr_\prec J^\infty(\eta\circ\tau)(\mathbb J_m)[\bar r]^*$ with an image of $\varphi^\vee_{\bar r}$.

\begin{cor}\label{DualToSubauotient}
We have an isomorphism of $\grad$-graded spaces
\[\gr_{\prec}J^{\infty}(\eta \circ \tau)(\mathbb J_m)[\bar r]^*\simeq
\varphi_{\bar r}^\vee\left( \bigcap_{\bar r'\prec \bar r, \bar r' \in \mathcal{R}(\bar a,L)}
\ker(\varphi_{\bar r'}^{\vee})\right)
\]
where $\ker$ denotes the kernel in $\Lambda_{\bar a,L}(\mathbf s)$.
\end{cor}

\subsection{Current algebra action and $\mathcal{A}$-freeness}
In this subsection we denote $R=R(P)$. Recall that each ring 
$S\in\{R,J^\infty(R),J^{\infty}_{\red}(R)\}$
possesses an $\bN^{n+1}$-grading with components $S[\bar a,L]$.
The ring $R$ admits a natural action of a 1-dimensional Lie algebra $\fh=\Bbbk h$ given by $hr=Lr$ for $r\in R[\bar a,L]$. In view of Lemma~\ref{currentaction} we obtain an $\fh[s]$-action on $J^{\infty}(R)$ given by $hs^k(r^{(j)})=(hr)^{(j-k)}$. By Lemma~\ref{NilpotentDerivations} we have an induced $\fh[s]$-action on $J^\infty_\red(R)$.
Clearly, the components $J_{\red}^{\infty}(R)[\bar a,L]$ are preserved by the action of  $\fh[s]$. 

Next, consider the algebra $\Bbbk[z_1,\dots,z_n,w]$, it also admits an $\fh$-action given by $h(M)=LM$ for $M=z_{i_1}\dots z_{i_N}w^L$. By Lemma~\ref{currentaction} we obtain an action of $\fh[s]$ on $J^\infty\Bbbk[z_1,\dots,z_n,w]$. Now recall the isomorphism \eqref{nu} \[\nu:J_{\red}^{\infty}(R)\to J^{\infty}(\eta \circ \tau)(\mathbb J_m).\] Since the right-hand side above is a subspace in $J^\infty\Bbbk[z_1,\dots,z_n,w]$, we will view $\nu$ as an embedding into the latter. It is then seen from the definitions that $\nu$ is a homomorphism of $\fh[s]$-algebras with respect to the $\fh[s]$-action on $J_{\red}^{\infty}(R)$ defined in the previous paragraph.

\begin{lem}
Every subspace $A_{\bar r}\subset \nu(J_{\red}^{\infty}(R))$ is preserved by the $\fh[s]$-action.
\end{lem}
\begin{proof}
For any $\bar\alpha^i$ we have $hs^k((z^{\bar\alpha^i}w)^{(j)})=(z^{\bar\alpha^i}w)^{(j-k)}$. Hence, applying $hs^k$ coefficientwise we have $hs^k((z^{\bar\alpha^i}w)(t))=t^k(z^{\bar\alpha^i}w)(t)$. Therefore,
\begin{equation}\label{derActionOnV}
hs^k(V_{\bar r})=\left(\sum_{i=1}^m 
\sum_{l=1}^{r_i} (t_i^{(l)})^k \right) V_{\bar r}.
\end{equation}
Proposition \ref{ASpan} now provides the claim.
\end{proof}

Let us fix $L\in \bN$ and $\bar a\in \bN^n$.

\begin{cor}
For any $\bar r\in\mathcal{R}(\bar a,L)$ we have an $\mathfrak{h}[s]$-action on the subquotient 
\[\sum_{\bar r' \preceq \bar r, \bar r' \in \mathcal{R}(\bar a,L)}A_{\bar r'}\left/\sum_{\bar r' \prec \bar r, \bar r' \in \mathcal{R}(\bar a,L)}A_{\bar r'}.\right.\] 
\end{cor}

Recall that for every $L$ the space $J_{\red}^{\infty}(R)[\bar a,L]$ is graded by  $\grad$ and the homogeneous components are finite-dimensional. We denote by $J_{\red}^{\infty}(R)[\bar a,L]^*$ the graded dual space. We denote by $\mathcal{A}_L\subset U(\fh[s])$ the subalgebra generated by $hs, \dots, hs^L$. We will be interested in situations when the action of $\mathcal{A}_L$ on $J_{\red}^{\infty}(R)[\bar a,L]^*$ is free for all $\bar a$ and $L$.

\begin{rem}
An analogous freeness property holds for the homogeneous coordinate ring of the Pl\"ucker embedding and it has important representation theoretic meaning (see \cite{FL, FM2, Kato}). This property also holds for the Veronese curve (see \cite{DF}). It is interesting to investigate the graded rank of these free modules. For example, in the case of the Pl\"ucker embedding they are equal to Macdonald polynomials for $t = 0$, see \cite{Ion, N}.   



\end{rem}





For $\bar r\in\mathcal{R}(\bar a,L)$ we define the following $\fh[s]$-action on the algebra $\Lambda_{\bar r}(\mathbf t)$:
\begin{equation}\label{derActionLambda}
hs^l(f)=\left(\sum_{i=1}^m 
\sum_{k=1}^{r_i} (t_i^{(k)})^l \right)f.
\end{equation}
\begin{prop}
The space $\varphi_{\bar r}^\vee(\Lambda_{\bar a,L}(\mathbf s))$ is invariant under the action \eqref{derActionLambda}. The action on this subspace is dual to the action of $\fh[s]$ on $A_{\bar r}$ with respect to the isomorphism $\varepsilon_{\bar r}$ (Theorem~\ref{inclusiondual}).
\end{prop}
\begin{proof}
To prove the first claim note that the first factor in the right-hand side of~\eqref{derActionLambda} lies in $\varphi_{\bar r}^\vee(\Lambda_{\bar a,L}(\mathbf s))$ by Remark~\ref{PowerSumsMaps}. The second claim follows from Equation~\eqref{derActionOnV}.
\end{proof}


\begin{lem}\label{freenessOfAction}
Assume that the natural action of $\mathcal{A}_{L}$ is  free on $\varphi_{\bar r}^\vee\left( \bigcap_{\bar r'\prec \bar r, \bar r' \in \mathcal{R}(\bar a,L)}
\ker(\varphi^\vee_{\bar r'})\right)$ for every $\bar r\in\mathcal{R}(\bar a,L)$. Then it is free on the space
$J^{\infty}(\eta \circ \tau)(\mathbb J_m)[\bar a,L]^*$.
\end{lem}
\begin{proof}
One sees that the isomorphism in Corollary \ref{DualToSubauotient} is one of $\fh[s]$-modules, thus, $J^{\infty}(\eta \circ \tau)(\mathbb J_m)[\bar a,L]^*$ has a filtration with free subquotients. However, free modules are projective, therefore, 
$J^{\infty}(\eta \circ \tau)(\mathbb J_m)[\bar a,L]^*$ is isomorphic to the direct sum of its subqotients.
\end{proof}

We will need the following general fact about symmetric polynomials.

\begin{lem}\label{FreeSymmetric}
Let $n = n_1 + \hdots + n_a$.
The ring $\bC[x_1, \hdots, x_{n}]^{\mathfrak{S}_{n_1} \times \hdots \times \mathfrak{S}_{n_a}}$ is free of $q$-rank $\binom{n}{n_1 \hdots n_a}_q$ over its subring $\bC[x_1, \hdots, x_n]^{\mathfrak{S}_{n}}$.
\end{lem}

\begin{proof}
We are going to prove that it is a direct summand of a free module.
Consider the tower of ring extensions:
\[
\bC[x_1, \hdots, x_{n}]^{\mathfrak{S}_{n}} \subset \bC[x_1, \hdots, x_{n}]^{\mathfrak{S}_{n_1} \times \hdots \times \mathfrak{S}_{n_a}} \subset \bC[x_1, \hdots, x_{n}].
\]
Note that the rightmost ring in this tower is isomorphic to 
\[\bigotimes_{i = 1}^a \bC[x_{n_1 + \hdots + n_{i - 1} + 1}, \hdots, x_{n_1 + \hdots + n_i}]\] 
and hence is free of $q$-rank $[n_1]_q! \hdots [n_a]_q!$ over the second ring in this tower, which is isomorphic to 
\[
\bigotimes_{i = 1}^a \bC[x_{n_1 + \hdots + n_{i - 1} + 1}, \hdots, x_{n_1 + \hdots + n_i}]^{\mathfrak{S}_{n_i}}
\]
(here we use the classical fact that a polynomial ring is free over its subring of symmetric polynomials).
Therefore, $\bC[x_1, \hdots, x_{n}]^{\mathfrak{S}_{n_1} \times \hdots \times \mathfrak{S}_{n_a}}$ is a direct summand in $\bC[x_1, \hdots, x_{n}]$ as a $\bC[x_1, \hdots, x_{n}]^{\mathfrak{S}_{n_1} \times \hdots \times \mathfrak{S}_{n_a}}$-module, and, in particular, as a $\bC[x_1, \hdots, x_{n}]^{\mathfrak{S}_{n}}$-module.
Hence, $\bC[x_1, \hdots, x_{n}]^{\mathfrak{S}_{n_1} \times \hdots \times \mathfrak{S}_{n_a}}$ is a projective module over $\bC[x_1, \hdots, x_{n}]^{\mathfrak{S}_{n}}$ which by the Quillen-Suslin theorem implies freeness.

Since the ranks of free modules in a tower of extensions multiply, the $q$-rank of our extension can be expressed as a quotient of ranks
\[
\frac{[n]_q!}{[n_1]_q! \hdots [n_a]_q!} = \binom{n}{n_1 \hdots n_a}_q.\qedhere
\]
\end{proof}

Denote by $I$ the kernel of the projection $\mathbb J_m\to J_{\red}^{\infty}(R(P))$. By Lemma~\ref{UpperBoundInitial} the quotient $\mathbb J_m/\In_\prec I$ is isomorphic to $\gr_\prec J_{\red}^{\infty}(R(P))$.

\begin{thm}
\label{equalityANDcofree}
Assume that $\In_\prec I$ is generated by relations of the form \eqref{QuadraticJetMonomialEquations}. Then 
\[\gr_\prec(J_{\red}^{\infty}(R(P)))[\bar a,L]\simeq J_{\red}^{\infty}(R(P))[\bar a,L]\]
as $\mathcal{A}_L$-modules and the dual of this module is free.
\end{thm}
\begin{proof}
Consider $\bar r\in\mathcal{R}(\bar a,L)$. By Proposition \ref{DualSpaceDerivedQuagraticMon} the space $\gr_\prec(J_{\red}^{\infty}(R(P)))[\bar r]^*$ 
is a free module over $\Lambda_{\bar r}(\mathbf t)$. By~\eqref{derActionLambda} the algebra $\mathcal A_L$ can be embedded into $\Lambda_{\bar r}(\mathbf t)$ as the subalgebra of polynomials symmetric in all $t_i^{(k)}$.  Therefore, applying Lemma \ref{FreeSymmetric} we deduce that $\gr_\prec(J_{\red}^{\infty}(R(P)))[\bar r]^*$ is free over $\mathcal{A}_L$. Lemma \ref{freenessOfAction} now shows that $\gr_\prec(J_{\red}^{\infty}(R(P))[\bar a,L]^*$ is free.
\end{proof}

\subsection{Technical lemma}\label{TechnicalLemma}
Consider integers $m\ge 1$ and $\zeta_1, \dots, \zeta_m \ge 0$ together with a vector $\bar a =(a_1,\dots, a_m,a) \in \mathbb{N}^{m+1}$. We also consider variables $u_i^{(j)}$ with $1 \leq i \leq m$, $1 \leq j \leq a_i$ and $s^{(j)}$ with $1 \leq j \leq a$. We work with the ring
\[\Lambda_{\bar a}(\mathbf u,\mathbf s):=\Bbbk[u_1^{(1)},\dots,u_1^{(a_1)},\dots, u_m^{(1)},\dots,u_m^{(a_m)},s^{(1)},\dots,s^{(a)}]^{\mathfrak{S}_{a_1}\times\dots \times\mathfrak{S}_{a_m}\times\mathfrak{S}_{a}}.\]

Within the set of collections of positive integers $\bar r=(r_{i,j})_{1 \leq i \leq m, 0 \leq j \leq \zeta_i}$ we distinguish the subset 
\[\mathcal{R}(\bar a):=\left\{ \bar r \left| \sum_{j=0}^{\zeta_i}r_{i,j}=a_i, \sum_{1 \leq i \leq m, 0 \leq j \leq \zeta_i} r_{i,j}j=a \right. \right\}.\]
For each $\bar r \in \mathcal{R}(\bar a)$ we have the ring
\[\Lambda_{\bar r} (\mathbf t):=\Bbbk[t_{i,j}^{(1)}, \dots, t_{i,j}^{(r_{i,j})}]_{1 \leq i \leq m, 0 \leq j \leq \zeta_i}^{\bigtimes_{i,j}\mathfrak{S}_{r_{i,j}}}.\] We now define a map $\psi_{\bar r}:\Lambda_{\bar a}(\mathbf u,\mathbf s)\to \Lambda_{\bar r} (\mathbf t)$. As for $\varphi^\vee_{\bar r}$ above we define $\psi_{\bar r}$ as a map between rings of all polynomials and then restrict to subrings of symmetric polynomials, the fact that the image lies in $\Lambda_{\bar r} (\mathbf t)$ is checked similarly to Lemma~\ref{ImageIsSymmetric}. Consider a variable $u_i^{(k)}$ and note that we have a unique $j_0\in[0,\zeta_i]$ for which $\sum_{j=1}^{j_0-1} r_{i,j}< k\le \sum_{j=1}^{j_0} r_{i,j}$. Denote $l=k-\sum_{j=1}^{j_0} r_{i,j}$, then we set $\psi_{\bar r}(u_i^{(k)})=t_{i,j_0}^{(l)}$. Now consider a variable $s^{(k)}$. We have a unique $i_0\in[1,m]$ for which \[\sum_{1\le i\le i_0-1,0\le j\le\zeta_i} j r_{i,j}< k\le \sum_{1\le i\le i_0,0\le j\le\zeta_i} j r_{i,j}.\] Denote the first sum above by $S$. We also have a unique $j_0\in[0,\zeta_{i_0}]$ for which \[S+\sum_{j=1}^{j_0-1} jr_{i_0,j}< k\le S+\sum_{j=1}^{j_0} jr_{i_0,j}.\] Finally, denote \[l=\left\lceil\frac{k-S-\sum_{j=1}^{j_0-1} jr_{i_0,j}}{j_0}\right\rceil,\] note that $l\in[1,r_{i_0,j_0}]$. We set $\psi_{\bar r}(s^{(k)})=t_{i_0,j_0}^{(l)}$. 

Note that $\psi_{\bar r}$ establishes a bijection between all variables of the form $u_i^{(k)}$ and all variables of the form $t_{i,j}^{(l)}$ for every $i$. Meanwhile, the number of variables of the form $s^{(k)}$ which are mapped into a given $t_{i,j}^{(l)}$ is equal to $j$.

\begin{rem}\label{dualtopsi}
The map $\psi_{\bar r}$ appears naturally in the following context. Consider lattice polytopes $P\subset\mathbb R^n$ and $Q\subset\mathbb R^{n+1}$ such that the following holds. Let $\bar \alpha^1,\dots,\bar \alpha^m$ be the integer points of $P$, then the integer points of $Q$ are precisely the points $\bar\beta_{i,j}=\bar \alpha^i\times j$ with $1 \leq i \leq m$, $0\le j\le\zeta_i$. In particular, the projection of $Q$ along the $(n+1)$st coordinate is $P$. The resulting embedding $Q\hookrightarrow P\times\mathbb R$ gives rise to the map 
\begin{equation}\label{psimeaning}
\Bbbk[X_{i,j}]_{1 \leq i \leq m, 0 \leq j \leq \zeta_i}\to\Bbbk[X_i]_{1\leq i\leq m}\otimes\Bbbk[z]
\end{equation}
taking $X_{i,j}$ to $X_i z^j$. We may now successively apply $J^\infty$ to this map, restrict to the component of grading $\bar r$ on the left and the component of grading $\bar a$ on the right and then dualize both sides. We obtain a map $\psi_{\bar r}:\Lambda_{\bar a}(\mathbf u,\mathbf s)\to \Lambda_{\bar r} (\mathbf t)$, this is precisely the map defined above.
\end{rem}

Consider a total order $\prec$ on pairs $(i,j)$, $1 \leq i \leq m, 0 \leq j \leq \zeta_i$ with the property $(i,j)\prec (i,j')$ if $j<j'$. One has an induced lexicographic order on vectors $\bar r_{i,j}$ such that $\bar r \prec \bar r'$ if and only if the $\prec$-minimal $(i,j)$ for which  $r_{i,j}\neq r'_{i,j}$ satisfies $r_{i,j}<r'_{i,j}$.


\begin{dfn}\label{kappa}
For $(i,j)\prec (i',j')$ let $\varkappa\big( (i,j), (i',j') \big)$ be the number of $l$ such that $(i,j)\prec (i',l)$ and $l<j'$. 
\end{dfn}

In other words, let $e_{min}$ be the minimal number $e$ such that $(i,j)\prec (i',e)$. Then $\varkappa\big( (i,j), (i',j') \big)=j'-e_{min}$. The key claim of this subsection is as follows.

\begin{lem}\label{Split}
\[\psi_{\bar r}\left(\bigcap_{\bar r \succ \bar r'}\ker (\psi_{\bar r'})\right)\supset\prod_{\substack{(i,j)\prec (i',j'), \\ 1 \leq l \leq r_{i,j}, \\1 \leq l' \leq r_{i',j'}}}\left( t_{i,j}^{(l)}-t_{i',j'}^{(l')} \right)^{\varkappa\big((i,j),(i',j')\big)} \Lambda_{\bar r}(\mathbf t).\]
\end{lem}

This lemma will be proved by induction reducing the case of vectors $\bar a$ and $\bar r$ to smaller vectors defined as follows. Choose $1\le\tilde i\le m$ with $\zeta_{\tilde i}>0$. For $\max(a_{\tilde i}-a,0)\le l\le a_{\tilde i}$ we define the truncation operator
\[\tr_l^{\tilde i}(\bar a)=(a_1,\dots, a_{\tilde i-1}, a_{\tilde i}-l,a_{\tilde i+1},\dots,a_m,a- a_{\tilde i}+l)\in\bN^{m+1}.\]
The vector $\bar r$ is transformed by the map $\tr^{\tilde i}$ with
\[\tr^{\tilde i}(\bar r)_{i,j}=
\begin{cases}
r_{i,j}~\text{if}~ i \neq \tilde i,\\
r_{i,j+1}~\text{if}~ i = \tilde i
\end{cases}\]
where $1\le i\le m$ and $1\le j\le\zeta_i$ for $i\neq\tilde i$ while $1\le j\le\zeta_{\tilde i}-1$ otherwise.
Note that we have
\[\tr^{\tilde i}(\bar r) \in \mathcal{R}\left( \tr_{r_{\tilde i,0}}^{\tilde i}(\bar a) \right).\] We now show that $\psi_{\bar r}$ factors through $\psi_{\tr^{\tilde i}(\bar r)}$.

We define the map \[\chi_l^{\tilde i}:\Lambda_{\bar a}(\mathbf u,\mathbf s) \rightarrow \Lambda_{\tr_l^{\tilde i}(\bar a)}(\mathbf u,\mathbf s)\otimes\Bbbk[t_{\tilde i, 0}^{(1)},\dots,t_{\tilde i, 0}^{(l)}]^{\mathfrak{S}_{l}}\]
by
\begin{align*}
u_i^{(j)} &\mapsto u_i^{(j)}~ \text{if}~ i \neq \tilde i,~u_{\tilde i}^{(j)} \mapsto u_{\tilde i}^{(j-l)}~ \text{if}~ j>l,~u_{\tilde i}^{(j)} \mapsto t_{\tilde i,0}^{(j)}~ \text{if}~  j\leq l,\\
s^{(j)}&\mapsto u_{\tilde i}^{(j)}~ \text{if}~j\le a_{\tilde i}-l, ~s^{(j)}\mapsto s^{(j-a_{\tilde i}+l)}~ \text{if}~ j>a_{\tilde i}-l.
\end{align*}
Let $\sh_{\tilde i}$ be the index shift operator:
\[\sh_{\tilde i}(t_{i,j})=t_{i,j} ~\text{if}~ i \neq \tilde i,~
\sh_{\tilde i}(t_{\tilde i,j})=t_{\tilde i,j+1}.\]
Then we have the following identity:
\begin{equation}\label{psidecomp}
    \psi_{\bar r}=\big((\sh_{\tilde i}\circ \psi_{\tr^{\tilde i}(\bar r)})\otimes \id\big)\circ \chi_{r_{\tilde i,0}}^{\tilde i}.
\end{equation}
For clarity let us point out the following. As mentioned, the maps $\psi_\bullet$, $\chi_\bullet^\bullet$ and $\sh_\bullet$ are defined as maps between polynomial rings and then restricted to symmetric polynomials. It should be noted that identity~\eqref{psidecomp} does not hold for these maps between polynomial rings (only for the restrictions) because the left-hand side and the right-hand side may not coincide on some $s^{(i)}$ (but always coincide up to the action of the corresponding symmetric group).

\begin{rem}
Identity~\eqref{psidecomp} can be interpreted in terms of Remark~\ref{dualtopsi} as follows. The map~\eqref{psimeaning} decomposes as \[\Bbbk[X_{i,j}]_{1 \leq i \leq m, 0 \leq j \leq \zeta_i}\to\Bbbk[X_i]_{1\leq i\leq m}\otimes\Bbbk[z]\otimes\Bbbk[X_{\tilde i,0}]\to\Bbbk[X_i]_{1\leq i\leq m}\otimes\Bbbk[z].\] Here the left arrow maps $X_{i,j}$ with $i\neq\tilde i$ to $X_iz^j$, maps $X_{\tilde i,j}$ with $j\ge 1$ to $X_{\tilde i}z^{j-1}$ and preserves $X_{\tilde i,0}$. The right arrow preserves $z$ and $X_i$ with $i\neq\tilde i$ while mapping $X_{\tilde i}$ to $X_{\tilde i}z$ and $X_{\tilde i,0}$ to $X_{\tilde i}$. We again apply $J^\infty$, restrict to components of gradings $\bar r$, $\tr^{\tilde i}_{r_{\tilde i,0}}(\bar a)\oplus r_{\tilde i,0}$ and $\bar a$ and then dualize. One may check that on the right we will obtain the map $\chi_{r_{\tilde i,0}}^{\tilde i}$ and on the left we will obtain $(\sh_{\tilde i}\circ \psi_{\tr^{\tilde i}(\bar r)})\otimes \id$ (the latter can be seen by further decomposing the left arrow).
\end{rem}

For $l < a_{\tilde i}$ we have:
\[\ker \chi_l^{\tilde i} \supset \ker \chi_{l+1}^{\tilde i}.\]
Moreover, we have $\chi_{l}^{\tilde i}=\tau\circ\chi_{l+1}^{\tilde i}$, where
\[\tau:\Lambda_{\tr_{l+1}^{\tilde i}(\bar a)}(\mathbf u,\mathbf s)\otimes\Bbbk[t_{\tilde i, 0}^{(1)},\dots,t_{\tilde i, 0}^{(l+1)}]^{\mathfrak{S}_{l+1}}\rightarrow \Lambda_{\tr_{l+1}^{\tilde i}(\bar a)}(\mathbf u,\mathbf s)\otimes \Bbbk[u_{\tilde i}^{(a_{\tilde i}-l)}]\otimes\Bbbk[t_{\tilde i, 0}^{(1)},\dots,t_{\tilde i, 0}^{(l)}]^{\mathfrak{S}_{l}},\]
\[\tau(t_{\tilde i,0}^{(l+1)})=\tau(s^{(a-a_{\tilde i}+l+1)})=u_{\tilde i}^{(a_{\tilde i}-l)}\]
and $\tau$ sends any variable (on the left) different from $t_{\tilde i,0}^{(l+1)}$ and  $s^{(a-a_{\tilde i}+l+1)}$ to the same variable on the right.

\begin{lem}\label{SplitInductionStep}
\[\chi_{l+1}^{\tilde i} \left(\ker \chi_{l}^{\tilde i}\right) \supset\prod_{j=1,\dots,a-a_{\tilde i}+l+1,j'=1,\dots,l+1}\left( s^{(j)}-t_{\tilde i, 0}^{(j')} \right)\left( \Lambda_{\tr_{l+1}^{\tilde i}(\bar a)}({\bf u,s})\otimes\Bbbk[t_{\tilde i, 0}^{(1)},\dots,t_{\tilde i, 0}^{(l+1)}]^{\mathfrak{S}_{l+1}}\right).\]
\end{lem}
\begin{proof}
One has
\begin{align*}\chi_{l+1}^{\tilde i}(p_k(u_{\tilde i,0}^{(1)},\dots,u_{\tilde i,0}^{(a_{\tilde i})}))&=p_k(u_{\tilde i,0}^{(1)},\dots,u_{\tilde i,0}^{(a_{\tilde i}-l-1)})+p_k(t_{\tilde i,0}^{(1)},\dots, t_{\tilde i,0}^{(l+1)}),\\
\chi_{l+1}^{\tilde i}(p_k(s^{(1)},\dots,s^{(a)}))&=p_k(u_{\tilde i,0}^{(1)},\dots,u_{\tilde i,0}^{(a_{\tilde i}-l-1)})+p_k(s^{(1)},\dots, s^{(a-a_{\tilde i}+l+1)}).
\end{align*}
Therefore, the image of $\chi^{\tilde i}_{l+1}$ contains the differences
\[p_k(t_{\tilde i,0}^{(1)},\dots, t_{\tilde i,0}^{(l+1)})-p_k(s^{(1)},\dots, s^{(a-a_{\tilde i}+l+1)})\]
which generate the ring $\Omega_{l+1,a-a_{\tilde i}+l+1}$ of supersymmetric polynomials in variables $t_{\tilde i, 0}^{(1)}, \dots, t_{\tilde i, 0}^{(l+1)}$ and $s^{(1)},\dots,s^{(a-a_{\tilde i}+l+1)}$ (Proposition \ref{SuperSymGenerators}). 
We deduce that the image of $\chi^{\tilde i}_{l+1}$ contains the ring 
\[
S=\Omega_{l+1,a-a_{\tilde i}+l+1}\otimes \Lambda_{a_{\tilde i}-l-1}(\mathbf{u}_{\tilde i})\otimes \bigotimes_{i \neq \tilde i}\Lambda_{a_i}(\mathbf{s}_i).
\]
The right-hand side in the statement of our lemma is equal to $\ker \tau \cap S$ due to Proposition \ref{resultant}. However, $\chi_{l+1}^{\tilde i}(\ker \chi_{l}^{\tilde i})=\im \chi_{l+1}^{\tilde i}\cap \ker \tau\supset S\cap \ker \tau$.
\end{proof}

\begin{proof}[Proof of Lemma~\ref{Split}]
Let $(\tilde i, 0)$ be the smallest pair $(i,j)$ in the order $\prec$. By definition we have:
\[\bar r \succ \bar r' \Longleftrightarrow r_{\tilde i, 0}>r'_{\tilde i, 0}~\text{or}~(r_{\tilde i, 0}=r'_{\tilde i, 0}~\text{and}~\tr^{\tilde i}(\bar r)\succ \tr^{\tilde i}(\bar r'))\]
where $\tr^{\tilde i}(\bar r)$ and $\tr^{\tilde i}(\bar r')$ are also compared lexicographically.

Using Lemma \ref{SplitInductionStep} and identity \eqref{psidecomp} we have:
\begin{equation}\label{inclusionsmallerchi}
\chi^{\tilde i}_{r_{\tilde i,0}}\left( \bigcap_{\bar r': r_{\tilde i,0} >  r'_{\tilde i,0}}\ker (\psi_{\bar r'}) \right)\supset \prod_{\substack{ l=1,\dots,a-a_{\tilde i}+r_{\tilde i, 0},\\ l'=1,\dots,r_{\tilde i, 0}}}\left( s^{(l)}-t_{\tilde i, 0}^{(l')} \right)\Lambda_{\tr_{r_{\tilde i, 0}}^{\tilde i}(\bar a)}({\bf u,s})\otimes\Bbbk[t_{\tilde i, 0}^{(1)},\dots,t_{\tilde i, 0}^{(r_{\tilde i, 0})}]^{\mathfrak{S}_{r_{\tilde i, 0}}}.
\end{equation}
We proceed by induction on $\sum_{i=1}^m \zeta_i$.
By the induction hypothesis we have
\[\psi_{\tr^{\tilde i}(\bar r)}\left(\bigcap_{\tr^{\tilde i}(\bar r) \succ \tr^{\tilde i}(\bar r')}\ker (\psi_{\tr^{\tilde i}(\bar r')})\right)\supset\prod_{\substack{(i,j)\prec (i',j'), \\ 1 \leq l \leq \tr^{\tilde i}(r)_{i,j}, \\1 \leq l' \leq \tr^{\tilde i}(r)_{i',j'}}}\left( t_{i,j}^{(l)}-t_{i',j'}^{(l')} \right)^{\varkappa\big((i,j),(i',j')\big)} \Lambda_{\tr^{\tilde i}(\bar r)}(\mathbf t).\]
Let $\bar r'\prec \bar r$ be such that $r'_{\tilde i,0}=r_{\tilde i,0}$.
Consider the set $F_{\bar r'}$ consisting 
of all 
$f\in \Lambda_{\tr_{r_{\tilde i,0}}^{\tilde i}(\bar a)}({\bf u,s})\otimes\Bbbk[t_{\tilde i, 0}^{(1)},\dots,t_{\tilde i, 0}^{(r_{\tilde i, 0})}]^{\mathfrak{S}_{r_{\tilde i, 0}}}$
such that 
$\left(\psi_{\tr^{\tilde i}(\bar r')}\otimes \id\right)f=0$.
Then by equation \eqref{psidecomp} 
\[
\psi_{\bar r}\left(\ker (\psi_{\bar r'})\right) = 
\big((\sh_{\tilde i}\circ \psi_{\tr^{\tilde i}(\bar r)})\otimes \id\big)\circ \chi_{r_{\tilde i,0}}^{\tilde i}\left(\ker (\psi_{\bar r'})\right),
\]
and hence  we get (note that $\sh_{\tilde i}$ is bijective)
\begin{align*}
\psi_{\bar r}\left(\ker (\psi_{\bar r'})\right) = 
\big((\sh_{\tilde i}\circ \psi_{\tr^{\tilde i}(\bar r)})\otimes \id\big)\left(\ker \big(\psi_{\tr^{\tilde i}(\bar r')}\otimes \id\big) \cap(\im \chi_{r_{\tilde i,0}}^{\tilde i})\right)\\
\supset \left((\sh_{\tilde i}\circ \psi_{\tr^{\tilde i}(\bar r)})\otimes \id\right)
\prod_{\substack{ l=1,\dots,a-a_{\tilde i}+r_{\tilde i, 0},\\ l'=1,\dots,r_{\tilde i, 0}}}\left( s^{(l)}-t_{\tilde i, 0}^{(l')} \right)F_{\bar r'} 
\end{align*}
due to Lemma \ref{SplitInductionStep}.
Taking the intersection one gets 
\[
\psi_{\bar r}\left(\bigcap_{\bar r' \prec \bar r,r'_{\tilde i,0}=r_{\tilde i,0}}\ker (\psi_{\bar r'})\right)\supset \left((\sh_{\tilde i}\circ \psi_{\tr^{\tilde i}(\bar r)})\otimes \id\right) \prod_{\substack{ l=1,\dots,a-a_{\tilde i}+r_{\tilde i 0},\\ l'=1,\dots,r_{\tilde i, 0}}}\left( s^{(l)}-t_{\tilde i, 0}^{(l')} \right) \bigcap_{\bar r' \prec \bar r,r'_{\tilde i,0}=r'_{\tilde i,0}}F_{\bar r'}.
\]
Therefore, using \eqref{inclusionsmallerchi} we obtain:
\[\psi_{\bar r}\left(\bigcap_{\bar r' \prec \bar r}\ker (\psi_{\bar r'})\right)\supset \left((\sh_{\tilde i}\circ \psi_{\tr^{\tilde i}(\bar r)})\otimes \id\right) \prod_{\substack{ l=1,\dots,a-a_{\tilde i}+r_{\tilde i, 0},\\ l'=1,\dots,r_{\tilde i, 0}}}\left( s^{(l)}-t_{\tilde i, 0}^{(l')} \right)  \bigcap_{\bar r' \prec \bar r,r'_{\tilde i,0}=r'_{\tilde i,0}}F_{\bar r'} .\]
Definition \ref{kappa} implies
\begin{multline*}
(\sh_{\tilde i}\circ \psi_{\tr^{\tilde i}(\bar r)})\otimes \id\left(\prod_{
\substack{l=1,\dots,a_s-a_{\tilde i}+r_{\tilde i, 0},\\ l'=1,\dots,r_{\tilde i, 0}}}\left( s^{(l)}-t_{\tilde i, 0}^{(l')} \right) \right)\\=
\prod_{(i,j)\succ (\tilde i, 0)} \prod_{\substack{l=1, \dots, r_{i,j},\\ l'=1,\dots,r_{\tilde i, 0}}}\left(t_{i,j}^{(l)}-t_{ \tilde i, 0}^{(l')}  \right)^{\varkappa\big((i,j),(\tilde i,0)\big)}.
\end{multline*}
By induction hypothesis we have:
\begin{multline*}
\left((\sh_{\tilde i}\circ \psi_{\tr^{\tilde i}(\bar r)})\otimes \id\right) \left( \bigcap_{\bar r' \prec \bar r,r'_{\tilde i,0}=r'_{\tilde i,0}}F_{\bar r'}\right)\supset\\ \prod_{(i,j)\succ (i',j')\succ(\tilde i,0)} \prod_{\substack{l=1, \dots, r_{i,j},\\ l'=1,\dots,r_{i',j'}}}\left(t_{i,j}^{(l)}-t_{i',j'}^{(l')}  \right)^{\varkappa\big((i,j),(i',j')\big)}\Lambda_{\tr_{r_{\tilde i,0}}^{\tilde i}(\bar a)}({\bf u,s})\otimes\Bbbk[t_{\tilde i, 0}^{(1)},\dots,t_{\tilde i, 0}^{(r_{\tilde i, 0})}]^{\mathfrak{S}_{r_{\tilde i, 0}}}.\qedhere
\end{multline*}
\end{proof}

\section{Inductive construction}\label{Inductive}

\subsection{Veronese curve}\label{Veronese}

In this subsection we study the arc schemes corresponding to the segments  $[0,\zeta]$ (recall that a one-dimensional polytope is affine equivalent to a segment).
By definition, the Veronese ring
$R([0,\zeta])$ is generated by $Y_\alpha=z^\alpha w$, $\alpha=0, \dots, \zeta$. 
The ring $J_{\red}^\infty(R([0,\zeta]))$
of reduced arcs over $R([0,\zeta])$
was studied in the paper \cite{DF}. Here we study this ring by methods developed above. Namely, we use the results from section \ref{hibireduced} on the structure of the reduced arc schemes of the cubes and Lemma \ref{Split}. 

Recall the notation $Y_I$ for the generators of the arc rings of cubes. 

For any $0 \leq \alpha <\beta \leq \zeta$
we define the following linear map of polytopes 
\[[0,1]^{\beta-\alpha}\rightarrow [0,\zeta],\ 
(d_1, \dots, d_{\beta-\alpha})\mapsto \alpha + \sum_{i=1}^{\beta-\alpha}d_i.\]
This map defines a homomorphism of corresponding toric rings $R([0,1])^{\beta-\alpha}\to R([0,\zeta])$:
\[\iota_{\alpha,\beta}: Y_I \mapsto Y_{\alpha+|I|}.\]


\begin{lem}\label{VeroneseDerivativeRelation}
For integers $0\le \alpha<\beta\le \zeta$ and $0\le k'\le\beta-\alpha-2$ the coefficients of the following series are zero in the ring $J_{\red}^{\infty}(R([0,\zeta]))$:
\[W'_{\alpha,\beta,k'}=\sum_{i=0}^{\beta-\alpha-1} (-1)^i \binom{\beta-\alpha-1}{i}\frac{\partial^{k'} Y_{\alpha + i}(s)}{\partial s^{k'}}Y_{\beta -i}(s).\]
\end{lem}
\begin{proof}
Consider the map $J^{\infty}(\iota_{\alpha,\beta}): 
J^{\infty}(R([0,1]^{\beta-\alpha})) \rightarrow J^{\infty}(R([0,\zeta]))$.
Recall the elements $W_{\beta-\alpha,k'} \in J^{\infty}(R([0,1]^{\beta-\alpha}))$. By Lemma \ref{NilpotentSeriesinCube}
the coefficients of $W_{\beta-\alpha,k'}$ are nilpotent. However, 
\[J^{\infty}(\iota_{\alpha,\beta})(W_{\beta-\alpha,k'})=W'_{\alpha,\beta,k'}.\]
Thus, the coefficients of $W'_{\beta-\alpha,k'}$ are nilpotent. Hence, they are equal to $0$ in $J_{\red}^{\infty}(R([0,\zeta]))$.
\end{proof}
\begin{dfn}
The  ring $\overline{J_{\red}^\infty(R([0,\zeta]))}$ is the quotient of $\Bbbk[X_i^{(j)}]_{0\le i\le \zeta, j\ge 0}$ modulo the ideal generated by coefficients of the series
\begin{equation}\label{VeroneseEquations}
\sum_{i=0}^{\beta-\alpha-1} (-1)^i \binom{\beta-\alpha-1}{i}\frac{\partial^{k'} X_{\alpha + i}(s)}{\partial s^{k'}}X_{\beta -i}(s),
\end{equation}
for all $0\le \alpha<\beta\le \zeta,~0\le k'\le\beta-\alpha-2$.
\end{dfn}

Clearly, there is a surjective homomorphism 
$\overline{J_{\red}^\infty(R([0,\zeta]))}\rightarrow J_{\red}^\infty(R([0,\zeta]))$.

Consider the following  (weight-restricted lexicographic) order on the semigroup $\mathbb{N}^{\zeta+1}$: we say that $(u_0,\dots,u_{\zeta})$ is less than $(v_0,\dots,v_{\zeta})$ if
$u_0+\dots+u_{\zeta}<v_0+\dots+v_{\zeta}$ or 
$u_0+\dots+u_{\zeta}=v_0+\dots+v_{\zeta}$ and
there exists an $i$ such that $u_i<v_i$ and $u_j=v_j$ for $j<i$. We have an induced order $\prec$ on the graded components of $J^{\infty}(\Bbbk[X_0, \dots, X_\zeta])$. The next lemma is obvious.
 \begin{lem}
 The initial part of \eqref{VeroneseEquations} is equal to 
 \begin{equation}\label{VeroneseDegeneratedRelation}
\frac{\partial^{k'}X_{\alpha}(s)}{\partial s^{k'}}X_{\beta}(s).
 \end{equation}
 \end{lem}
 
We define the ring $\widehat{J_{\red}^\infty(R([0,\zeta]))}$
  to be a quotient of $\Bbbk[X_i^{(j)}]_{i=0, \dots, \zeta, j\ge 0}$ by the ideal generated by coefficients of \eqref{VeroneseDegeneratedRelation}. These relations are homogeneous with respect to the grading $\prec$. Thus the ring $\widehat{J_{\red}^\infty(R([0,\zeta]))}$ has a natural $\mathbb{N}^{\zeta+1}$-grading and we denote the $\bar r$th graded component by 
  $\widehat{J_{\red}^\infty(R([0,\zeta]))}[\bar r]$.
Proposition \ref{DualSpaceDerivedQuagraticMon} applied to this graded component gives us the following.
\begin{prop}
We have the following isomorphism of $\grad$-graded spaces:
\[\left(\widehat{J_{\red}^\infty(R([0,\zeta]))}[\bar r]\right)^* \simeq \prod_{0 \leq \alpha<\beta\leq \zeta, 1 \leq i \leq r_{\alpha}, 1 \leq j \leq r_{\beta} }\left(t_{\alpha}^{(i)}-t_{\beta}^{(j)} \right)^{\beta-\alpha-1}
\Bbbk[t_\alpha^{(j)}]_{0 \leq \alpha \leq \zeta, 1 \leq j \leq r_\alpha}^{\mathfrak{S}_{r_0}\times \dots \times \mathfrak{S}_{r_\zeta}}.\] 
\end{prop}


Taking into account the grading  $\grad$ we extend the $\mathbb{N}^{\zeta+1}$-grading above to the $\mathbb{N}^{\zeta+2}$-grading.  Corollary \ref{UpperBoundSubquotient} implies
\begin{lem}\label{subquotientUpperBound}
There is the following surjective map of $\mathbb{N}^{\zeta+2}$-graded rings:
\[\widehat{J_{\red}^\infty(R([0,\zeta]))}\twoheadrightarrow
\gr_\prec J_{\red}^\infty(R([0,\zeta])).\]
\end{lem}


We take a vector $\bar r =(r_0, \dots, r_\zeta)$ and denote
\[\sum_{i=0}^\zeta ir_i=a, \sum_{i=0}^\zeta r_i=L.\]
Recall the notation $\bar r \in \mathcal{R}(a,L)$. We consider the rings of symmetric functions $\Lambda_{a,L}(\mathbf s)$
and $\Lambda_{\bar r}(\mathbf t)$.
We have the map $\varphi_{\bar r}: A_{\bar r} \rightarrow J^{\infty}(\Bbbk[z_1,w])[a,L]$ for $\bar r \in \mathcal{R}(a,L)$
and the map $\varphi_{\bar r}^\vee:\Lambda_{a,L}(\mathbf s) \rightarrow \Lambda_{\bar r}(\mathbf t)$ defined by the the following formulas (see \eqref{SymmetricAlgebrasMapDefinition}):
\begin{multline}
s_w^{(1)} \mapsto t_0^{(1)};\dots;s_w^{(r_0)} \mapsto t_0^{(r_0)};
s_w^{(r_0+1)} , s_1^{(1)} \mapsto t_1^{(1)};\dots;s_w^{(r_0+r_1)} ,
s_1^{(r_1)} \mapsto t_1^{(r_1)};\\
s_w^{(r_0+r_1+1)} , s_1^{(r_1+1)} , s_1^{(r_1+2)} \mapsto t_2^{(1)};\dots,s_w^{(r_0+r_1+r_2)},
s_1^{(r_1+2r_2-1)},
s_1^{(r_1+2r_2)} \mapsto  t_2^{(r_2)};\\
s_w^{(r_0+r_1+r_2+1)}, s_1^{(r_1+2r_2+1)},  s_1^{(r_1+2r_2+2)}, s_1^{(r_1+2r_2+3)} \mapsto t_3^{(1)}\dots
\end{multline}
One has
\begin{gather*}
 \varphi_{\bar r}^\vee\left(p_k\left(s_w^{(1)},\dots,s_w^{(L)} \right)\right)=\sum_{i=0}^\zeta p_k\left(t_i^{(1)},\dots,t_i^{(r_i)} \right),\\
 \varphi_{\bar r}^\vee\left(p_k\left(s_1^{(1)},\dots,s_1^{(a)} \right)\right)=\sum_{i=0}^\zeta i p_k\left(t_i^{(1)},\dots,t_i^{(r_i)} \right).
\end{gather*}

Lemma \ref{Split} with $m=1$ implies

\begin{prop}\label{veroneseImageofKernel}
For any $\bar r \in \mathcal{R}(a,L)$
\begin{equation}
\varphi_{\bar r}^\vee\left(\bigcap_{\bar r'\in \mathcal{R}(a,L),\bar r' \prec \bar r} \ker (\varphi_{\bar r'}^\vee)\right)\supset
\prod_{0 \leq \alpha<\beta\leq \zeta, 1 \leq i \leq r_{\alpha}, 1 \leq j \leq r_{\beta} }\left(t_{\alpha}^{(i)}-t_{\beta}^{(j)} \right)^{\beta-\alpha-1}
\Lambda_{\bar r}(\mathbf t).
\end{equation}
\end{prop}




Thus we reobtain the following theorem from \cite{DF}.
\begin{thm}
For any $\bar r$ $\grad$-graded vector spaces $\gr_\prec J_{\red}^{\infty}(R([0,\zeta]))[\bar r]
$ and $\widehat{\big(J_{\red}^\infty(R([0,\zeta]))\big)}[\bar r]$ are naturally  isomorphic. One has the isomorphism of graded rings
\[J_{\red}^{\infty}(R([0,\zeta]))\simeq \overline{J_{\red}^{\infty}(R([0,\zeta]))}.
\]
\end{thm}
\begin{cor}
The action of $\mathcal{A}_L$ is free on $J_{\red}^{\infty}(R([0,\zeta])[L]^*$.
\end{cor}
\begin{proof}
It follows from Theorem \ref{equalityANDcofree}.
\end{proof}

\subsection{Two dimensional case}
Let $f:[0,\eta] \rightarrow \mathbb{R}_{\geq 0}$ be a convex function such that for any integer $i \in [0,\eta]$: $\zeta_i:=f(i)\in \mathbb{N}$.
In this subsection we consider the case of the two-dimensional polygon obtained as the convex hull of the segment  $[0,\eta]\times \{0\}$ and all points $(i,f(i)), i \in [0,\eta]$. In other words, we consider a positive integer $\eta$ and let $\zeta_0, \dots, \zeta_\eta$ be a tuple of nonnegative integers such that $\zeta_0-\zeta_1 \leq \zeta_1-\zeta_2 \leq \dots \leq \zeta_{\eta-1}-\zeta_\eta$. Then the integer points  of the corresponding convex lattice polygon are $(i,j)|0 \leq i \leq \eta, 0 \leq j \leq \zeta_i$. 
\begin{example}
The case $\zeta=\begin{cases}
\zeta(\alpha)=\alpha/b,~ \alpha \in [0,b];\\
\zeta(\alpha)=1,~\alpha \in [b,\eta]
\end{cases},$
$0 \leq b \leq \eta$ corresponds to Hirzebruch surface.
\end{example}

\begin{example}
The case $\zeta(\alpha)=\alpha$, $\alpha \in [0,\eta]$ corresponds to
Veronese embedding of $\mathbb{P}^2$.
\end{example}

It will be convenient for us to reflect this polygon across the horizontal line $j=\zeta_{\max}/2$, where  $\zeta_{\max}=\max\{\zeta_i|i=0, \dots,\eta\}$. The integer points of the reflected polytope $P$ are $(i,j)$ with $0 \leq i \leq \eta$ and $\zeta_{\max} -\zeta_i \leq j \leq \zeta_{\max}$. An example of such a reflected polygon $P$ is seen in the figures below.

Consider two distinct integer points $\bar\alpha=(\alpha_1,\alpha_2),\bar\beta=(\beta_1,\beta_2) \in P$. We  assume that $\alpha_1 \leq \beta_1$ and if $\alpha_1 = \beta_1$, then $\alpha_2 \leq \beta_2$. We will now define a certain polygonal curve with integer vertices connecting these two points. Denote the vertices of this curve by $(i_0,j_0),\dots,(i_{k+1},j_{k+1})$, we require the following to hold.
\begin{itemize}
\item $(i_0,j_0)=\bar \alpha$ and $(i_{k+1},j_{k+1})=\bar\beta$.
\item The sequence $i_0,\dots,i_{k+1}$ is weakly increasing.
\item The sequence $j_0,\dots,j_{k+1}$ is weakly monotonic, it increases (weakly) if $\alpha_2\le\beta_2$ and decreases (weakly) if $\alpha_2\ge\beta_2$.
\item For every $1\le l\le {k+1}$ we either have $i_l-i_{l-1}=1$ or $i_l=i_{l-1}$ and $|j_{l}-j_{l-1}|=1$.
\item The curve is concave in the sense that a segment connecting two points of the curve cannot contain points lying below the curve.
\item All $(i_l,j_l)\in P$.
\end{itemize}
This provides a concave polygonal curve with integer vertices in $P$ which is directed either towards the top-right or the bottom-right. However, these conditions still leave plenty of degrees of freedom, the curve we choose can be described as follows (modulo one exceptional case which will be discussed below). The curve contains as many vertices as possible on the horizontal line $i=\min(\alpha_2,\beta_2)$ (i.e.\ that which contains the lower endpoint) and as many vertices as possible on the vertical line which contains the higher endpoint (i.e.\ either $j=\alpha_1$ or $j=\beta_1$ depending on $\max(\alpha_2,\beta_2)$). All intermediate vertices (if any) are chosen as vertices of $P$.
In other words, for $1\le l\le {k+1}$ let $e_l=(i_l,j_l)-(i_{l-1},j_{l-1})$ be the vector corresponding to the $l$th segment, explicit formulas for these vectors can be given as follows.




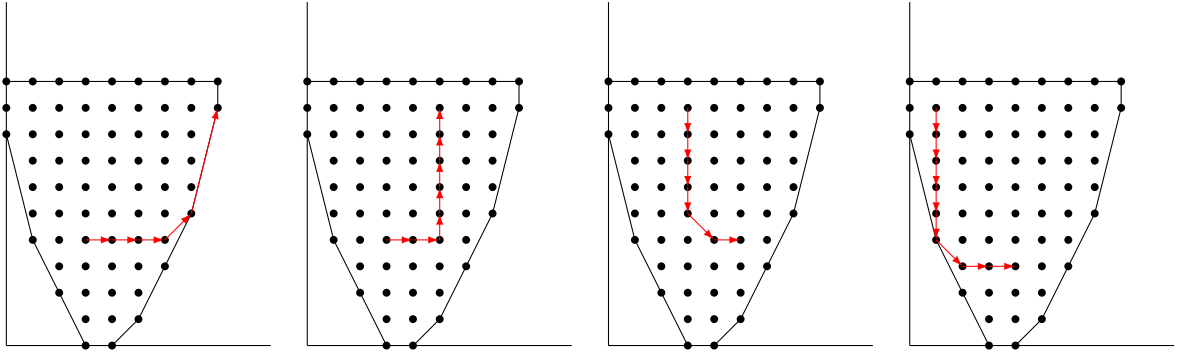
\begin{picture}(110.00,140.00)(-5,00)
\put(5.00,5.00){\line(1,0){100}}
\put(5.00,5.00){\line(0,1){130}}
\put(85.00,105.00){\line(0,-1){10}}
\put(85.00,95.00){\line(-1,-4){10}}
\put(75.00,55.00){\line(-1,-2){20}}
\put(55.00,15.00){\line(-1,-1){10}}
\put(5.00,105.00){\line(1,0){80}}
\put(5.00,85.00){\line(1,-4){10}}
\put(15.00,45.00){\line(1,-2){20}}
\put(5.00,105.00){\circle*{3}}
\put(5.00,95.00){\circle*{3}}
\put(5.00,85.00){\circle*{3}}
\put(15.00,105.00){\circle*{3}}
\put(15.00,95.00){\circle*{3}}
\put(15.00,85.00){\circle*{3}}
\put(15.00,75.00){\circle*{3}}
\put(15.00,65.00){\circle*{3}}
\put(15.00,55.00){\circle*{3}}
\put(15.00,45.00){\circle*{3}}
\put(25.00,105.00){\circle*{3}}
\put(25.00,95.00){\circle*{3}}
\put(25.00,85.00){\circle*{3}}
\put(25.00,75.00){\circle*{3}}
\put(25.00,65.00){\circle*{3}}
\put(25.00,55.00){\circle*{3}}
\put(25.00,45.00){\circle*{3}}
\put(25.00,35.00){\circle*{3}}
\put(25.00,25.00){\circle*{3}}
\put(35.00,105.00){\circle*{3}}
\put(35.00,95.00){\circle*{3}}
\put(35.00,85.00){\circle*{3}}
\put(35.00,75.00){\circle*{3}}
\put(35.00,65.00){\circle*{3}}
\put(35.00,55.00){\circle*{3}}
\put(35.00,45.00){\circle*{3}}
\put(35.00,35.00){\circle*{3}}
\put(35.00,25.00){\circle*{3}}
\put(35.00,15.00){\circle*{3}}
\put(35.00,5.00){\circle*{3}}
\put(45.00,105.00){\circle*{3}}
\put(45.00,95.00){\circle*{3}}
\put(45.00,85.00){\circle*{3}}
\put(45.00,75.00){\circle*{3}}
\put(45.00,65.00){\circle*{3}}
\put(45.00,55.00){\circle*{3}}
\put(45.00,45.00){\circle*{3}}
\put(45.00,35.00){\circle*{3}}
\put(45.00,25.00){\circle*{3}}
\put(45.00,15.00){\circle*{3}}
\put(45.00,5.00){\circle*{3}}
\put(55.00,105.00){\circle*{3}}
\put(55.00,95.00){\circle*{3}}
\put(55.00,85.00){\circle*{3}}
\put(55.00,75.00){\circle*{3}}
\put(55.00,65.00){\circle*{3}}
\put(55.00,55.00){\circle*{3}}
\put(55.00,45.00){\circle*{3}}
\put(55.00,35.00){\circle*{3}}
\put(55.00,25.00){\circle*{3}}
\put(55.00,15.00){\circle*{3}}
\put(65.00,105.00){\circle*{3}}
\put(65.00,95.00){\circle*{3}}
\put(65.00,85.00){\circle*{3}}
\put(65.00,75.00){\circle*{3}}
\put(65.00,65.00){\circle*{3}}
\put(65.00,55.00){\circle*{3}}
\put(65.00,45.00){\circle*{3}}
\put(65.00,35.00){\circle*{3}}
\put(75.00,105.00){\circle*{3}}
\put(75.00,95.00){\circle*{3}}
\put(75.00,85.00){\circle*{3}}
\put(75.00,75.00){\circle*{3}}
\put(75.00,65.00){\circle*{3}}
\put(75.00,55.00){\circle*{3}}
\put(85.00,105.00){\circle*{3}}
\put(85.00,95.00){\circle*{3}}
{\color{red}
\put(35.00,45.00){\vector(1,0){10}}
\put(45.00,45.00){\vector(1,0){10}}
\put(55.00,45.00){\vector(1,0){10}}
\put(65.00,45.00){\vector(1,1){10}}
\put(75.00,55.00){\vector(1,4){10}}
}
\end{picture}
\begin{picture}(110.00,140.00)(-5,00)
\put(5.00,5.00){\line(1,0){100}}
\put(5.00,5.00){\line(0,1){130}}
\put(85.00,105.00){\line(0,-1){10}}
\put(85.00,95.00){\line(-1,-4){10}}
\put(75.00,55.00){\line(-1,-2){20}}
\put(55.00,15.00){\line(-1,-1){10}}
\put(5.00,105.00){\line(1,0){80}}
\put(5.00,85.00){\line(1,-4){10}}
\put(15.00,45.00){\line(1,-2){20}}
\put(5.00,105.00){\circle*{3}}
\put(5.00,95.00){\circle*{3}}
\put(5.00,85.00){\circle*{3}}
\put(15.00,105.00){\circle*{3}}
\put(15.00,95.00){\circle*{3}}
\put(15.00,85.00){\circle*{3}}
\put(15.00,75.00){\circle*{3}}
\put(15.00,65.00){\circle*{3}}
\put(15.00,55.00){\circle*{3}}
\put(15.00,45.00){\circle*{3}}
\put(25.00,105.00){\circle*{3}}
\put(25.00,95.00){\circle*{3}}
\put(25.00,85.00){\circle*{3}}
\put(25.00,75.00){\circle*{3}}
\put(25.00,65.00){\circle*{3}}
\put(25.00,55.00){\circle*{3}}
\put(25.00,45.00){\circle*{3}}
\put(25.00,35.00){\circle*{3}}
\put(25.00,25.00){\circle*{3}}
\put(35.00,105.00){\circle*{3}}
\put(35.00,95.00){\circle*{3}}
\put(35.00,85.00){\circle*{3}}
\put(35.00,75.00){\circle*{3}}
\put(35.00,65.00){\circle*{3}}
\put(35.00,55.00){\circle*{3}}
\put(35.00,45.00){\circle*{3}}
\put(35.00,35.00){\circle*{3}}
\put(35.00,25.00){\circle*{3}}
\put(35.00,15.00){\circle*{3}}
\put(35.00,5.00){\circle*{3}}
\put(45.00,105.00){\circle*{3}}
\put(45.00,95.00){\circle*{3}}
\put(45.00,85.00){\circle*{3}}
\put(45.00,75.00){\circle*{3}}
\put(45.00,65.00){\circle*{3}}
\put(45.00,55.00){\circle*{3}}
\put(45.00,45.00){\circle*{3}}
\put(45.00,35.00){\circle*{3}}
\put(45.00,25.00){\circle*{3}}
\put(45.00,15.00){\circle*{3}}
\put(45.00,5.00){\circle*{3}}
\put(55.00,105.00){\circle*{3}}
\put(55.00,95.00){\circle*{3}}
\put(55.00,85.00){\circle*{3}}
\put(55.00,75.00){\circle*{3}}
\put(55.00,65.00){\circle*{3}}
\put(55.00,55.00){\circle*{3}}
\put(55.00,45.00){\circle*{3}}
\put(55.00,35.00){\circle*{3}}
\put(55.00,25.00){\circle*{3}}
\put(55.00,15.00){\circle*{3}}
\put(65.00,105.00){\circle*{3}}
\put(65.00,95.00){\circle*{3}}
\put(65.00,85.00){\circle*{3}}
\put(65.00,75.00){\circle*{3}}
\put(65.00,65.00){\circle*{3}}
\put(65.00,55.00){\circle*{3}}
\put(65.00,45.00){\circle*{3}}
\put(65.00,35.00){\circle*{3}}
\put(75.00,105.00){\circle*{3}}
\put(75.00,95.00){\circle*{3}}
\put(75.00,85.00){\circle*{3}}
\put(75.00,75.00){\circle*{3}}
\put(75.00,65.00){\circle*{3}}
\put(75.00,55.00){\circle*{3}}
\put(85.00,105.00){\circle*{3}}
\put(85.00,95.00){\circle*{3}}
{\color{red}
\put(55.00,85.00){\vector(0,1){10}}
\put(55.00,75.00){\vector(0,1){10}}
\put(55.00,65.00){\vector(0,1){10}}
\put(55.00,55.00){\vector(0,1){10}}
\put(55.00,45.00){\vector(0,1){10}}
\put(35.00,45.00){\vector(1,0){10}}
\put(45.00,45.00){\vector(1,0){10}}
}
\end{picture}
\begin{picture}(110.00,140.00)(-5,00)
\put(5.00,5.00){\line(1,0){100}}
\put(5.00,5.00){\line(0,1){130}}
\put(85.00,105.00){\line(0,-1){10}}
\put(85.00,95.00){\line(-1,-4){10}}
\put(75.00,55.00){\line(-1,-2){20}}
\put(55.00,15.00){\line(-1,-1){10}}
\put(5.00,105.00){\line(1,0){80}}
\put(5.00,85.00){\line(1,-4){10}}
\put(15.00,45.00){\line(1,-2){20}}
\put(5.00,105.00){\circle*{3}}
\put(5.00,95.00){\circle*{3}}
\put(5.00,85.00){\circle*{3}}
\put(15.00,105.00){\circle*{3}}
\put(15.00,95.00){\circle*{3}}
\put(15.00,85.00){\circle*{3}}
\put(15.00,75.00){\circle*{3}}
\put(15.00,65.00){\circle*{3}}
\put(15.00,55.00){\circle*{3}}
\put(15.00,45.00){\circle*{3}}
\put(25.00,105.00){\circle*{3}}
\put(25.00,95.00){\circle*{3}}
\put(25.00,85.00){\circle*{3}}
\put(25.00,75.00){\circle*{3}}
\put(25.00,65.00){\circle*{3}}
\put(25.00,55.00){\circle*{3}}
\put(25.00,45.00){\circle*{3}}
\put(25.00,35.00){\circle*{3}}
\put(25.00,25.00){\circle*{3}}
\put(35.00,105.00){\circle*{3}}
\put(35.00,95.00){\circle*{3}}
\put(35.00,85.00){\circle*{3}}
\put(35.00,75.00){\circle*{3}}
\put(35.00,65.00){\circle*{3}}
\put(35.00,55.00){\circle*{3}}
\put(35.00,45.00){\circle*{3}}
\put(35.00,35.00){\circle*{3}}
\put(35.00,25.00){\circle*{3}}
\put(35.00,15.00){\circle*{3}}
\put(35.00,5.00){\circle*{3}}
\put(45.00,105.00){\circle*{3}}
\put(45.00,95.00){\circle*{3}}
\put(45.00,85.00){\circle*{3}}
\put(45.00,75.00){\circle*{3}}
\put(45.00,65.00){\circle*{3}}
\put(45.00,55.00){\circle*{3}}
\put(45.00,45.00){\circle*{3}}
\put(45.00,35.00){\circle*{3}}
\put(45.00,25.00){\circle*{3}}
\put(45.00,15.00){\circle*{3}}
\put(45.00,5.00){\circle*{3}}
\put(55.00,105.00){\circle*{3}}
\put(55.00,95.00){\circle*{3}}
\put(55.00,85.00){\circle*{3}}
\put(55.00,75.00){\circle*{3}}
\put(55.00,65.00){\circle*{3}}
\put(55.00,55.00){\circle*{3}}
\put(55.00,45.00){\circle*{3}}
\put(55.00,35.00){\circle*{3}}
\put(55.00,25.00){\circle*{3}}
\put(55.00,15.00){\circle*{3}}
\put(65.00,105.00){\circle*{3}}
\put(65.00,95.00){\circle*{3}}
\put(65.00,85.00){\circle*{3}}
\put(65.00,75.00){\circle*{3}}
\put(65.00,65.00){\circle*{3}}
\put(65.00,55.00){\circle*{3}}
\put(65.00,45.00){\circle*{3}}
\put(65.00,35.00){\circle*{3}}
\put(75.00,105.00){\circle*{3}}
\put(75.00,95.00){\circle*{3}}
\put(75.00,85.00){\circle*{3}}
\put(75.00,75.00){\circle*{3}}
\put(75.00,65.00){\circle*{3}}
\put(75.00,55.00){\circle*{3}}
\put(85.00,105.00){\circle*{3}}
\put(85.00,95.00){\circle*{3}}
{\color{red}
\put(35.00,95.00){\vector(0,-1){10}}
\put(35.00,85.00){\vector(0,-1){10}}
\put(35.00,75.00){\vector(0,-1){10}}
\put(35.00,65.00){\vector(0,-1){10}}
\put(35.00,55.00){\vector(1,-1){10}}
\put(45.00,45.00){\vector(1,0){10}}
}
\end{picture}
\begin{picture}(110.00,140.00)(-5,00)
\put(5.00,5.00){\line(1,0){100}}
\put(5.00,5.00){\line(0,1){130}}
\put(85.00,105.00){\line(0,-1){10}}
\put(85.00,95.00){\line(-1,-4){10}}
\put(75.00,55.00){\line(-1,-2){20}}
\put(55.00,15.00){\line(-1,-1){10}}
\put(5.00,105.00){\line(1,0){80}}
\put(5.00,85.00){\line(1,-4){10}}
\put(15.00,45.00){\line(1,-2){20}}
\put(5.00,105.00){\circle*{3}}
\put(5.00,95.00){\circle*{3}}
\put(5.00,85.00){\circle*{3}}
\put(15.00,105.00){\circle*{3}}
\put(15.00,95.00){\circle*{3}}
\put(15.00,85.00){\circle*{3}}
\put(15.00,75.00){\circle*{3}}
\put(15.00,65.00){\circle*{3}}
\put(15.00,55.00){\circle*{3}}
\put(15.00,45.00){\circle*{3}}
\put(25.00,105.00){\circle*{3}}
\put(25.00,95.00){\circle*{3}}
\put(25.00,85.00){\circle*{3}}
\put(25.00,75.00){\circle*{3}}
\put(25.00,65.00){\circle*{3}}
\put(25.00,55.00){\circle*{3}}
\put(25.00,45.00){\circle*{3}}
\put(25.00,35.00){\circle*{3}}
\put(25.00,25.00){\circle*{3}}
\put(35.00,105.00){\circle*{3}}
\put(35.00,95.00){\circle*{3}}
\put(35.00,85.00){\circle*{3}}
\put(35.00,75.00){\circle*{3}}
\put(35.00,65.00){\circle*{3}}
\put(35.00,55.00){\circle*{3}}
\put(35.00,45.00){\circle*{3}}
\put(35.00,35.00){\circle*{3}}
\put(35.00,25.00){\circle*{3}}
\put(35.00,15.00){\circle*{3}}
\put(35.00,5.00){\circle*{3}}
\put(45.00,105.00){\circle*{3}}
\put(45.00,95.00){\circle*{3}}
\put(45.00,85.00){\circle*{3}}
\put(45.00,75.00){\circle*{3}}
\put(45.00,65.00){\circle*{3}}
\put(45.00,55.00){\circle*{3}}
\put(45.00,45.00){\circle*{3}}
\put(45.00,35.00){\circle*{3}}
\put(45.00,25.00){\circle*{3}}
\put(45.00,15.00){\circle*{3}}
\put(45.00,5.00){\circle*{3}}
\put(55.00,105.00){\circle*{3}}
\put(55.00,95.00){\circle*{3}}
\put(55.00,85.00){\circle*{3}}
\put(55.00,75.00){\circle*{3}}
\put(55.00,65.00){\circle*{3}}
\put(55.00,55.00){\circle*{3}}
\put(55.00,45.00){\circle*{3}}
\put(55.00,35.00){\circle*{3}}
\put(55.00,25.00){\circle*{3}}
\put(55.00,15.00){\circle*{3}}
\put(65.00,105.00){\circle*{3}}
\put(65.00,95.00){\circle*{3}}
\put(65.00,85.00){\circle*{3}}
\put(65.00,75.00){\circle*{3}}
\put(65.00,65.00){\circle*{3}}
\put(65.00,55.00){\circle*{3}}
\put(65.00,45.00){\circle*{3}}
\put(65.00,35.00){\circle*{3}}
\put(75.00,105.00){\circle*{3}}
\put(75.00,95.00){\circle*{3}}
\put(75.00,85.00){\circle*{3}}
\put(75.00,75.00){\circle*{3}}
\put(75.00,65.00){\circle*{3}}
\put(75.00,55.00){\circle*{3}}
\put(85.00,105.00){\circle*{3}}
\put(85.00,95.00){\circle*{3}}
{\color{red}
\put(15.00,95.00){\vector(0,-1){10}}
\put(15.00,85.00){\vector(0,-1){10}}
\put(15.00,75.00){\vector(0,-1){10}}
\put(15.00,65.00){\vector(0,-1){10}}
\put(15.00,55.00){\vector(0,-1){10}}
\put(15.00,45.00){\vector(1,-1){10}}
\put(25.00,35.00){\vector(1,0){10}}
\put(35.00,35.00){\vector(1,0){10}}
}
\end{picture}

Suppose $\alpha_2\le\beta_2$, i.e.\ the curve is directed towards the top-right (first and second figure above). Choose the maximal $u\le\beta_1$ for which $(\alpha_1,u)\in P$.
\begin{align*}
&e_1=\dots=e_{u-\alpha_1 }=(1,0), \\
&e_{u-\alpha_1 +1}=(1,\zeta_{\max}-\zeta_{u+1}-\alpha_2),\\
&e_{u-\alpha_1 +2}=(1,\zeta_{u+2}-\zeta_{u+1}),\dots,e_{\beta_1-\alpha_1}=(1,\zeta_{\beta_1-1}-\zeta_{\beta_1}) \\
&e_{\beta_1-\alpha_1+1}=\dots = e_{k+1}=(0,1).
\end{align*}

Now suppose $\alpha_2 >\beta_2$ and choose the minimal $u\ge\alpha_1$ for which $(u,\beta_2)\in P$. First suppose that $u=\alpha_1$ (so that $(\alpha_1,\beta_2)\in P$, see third figure), this is the exceptional case mentioned above, note vector $e_{\alpha_2-\beta_2}$. Here we set
\[e_1=\dots=e_{\alpha_2-\beta_2-1}=(0,-1),e_{\alpha_2-\beta_2}=(1,-1), e_{\alpha_2-\beta_2+1}=\dots=e_{k+1}=(1,0).\]
Finally, suppose $(\alpha_1,\beta_2)\notin P$ (fourth figure), then we proceed similarly to the case of $\alpha_2 <\beta_2$.
\begin{align*}
&e_1=\dots=e_{\alpha_2 -(\zeta_{\max}-\zeta_{\alpha_1})}=(0,-1),\\
&e_{\alpha_2 -(\zeta_{\max}-\zeta_{\alpha_1})+1}=(1,\zeta_{\alpha_1}-\zeta_{\alpha_1+1}),\dots,e_{\alpha_2 -(\zeta_{\max}-\zeta_{\alpha_1})+u-1-\alpha_1}=(1,\zeta_{u-2}-\zeta_{u-1}), \\
&e_{\alpha_2 -(\zeta_{\max}-\zeta_{\alpha_1})+u-\alpha_1}=(1,\zeta_{\max}-\zeta_{u-1}-\beta_2), \\
&e_{\alpha_2 -(\zeta_{\max}-\zeta_{\alpha_1})+u-\alpha_1+1}=\dots=e_{k+1}=(1,0).
\end{align*}

We will use the notations $k(\bar\alpha,\bar\beta)=k(\bar\beta,\bar\alpha)=k$ and $e_l(\bar\alpha,\bar \beta)=e_l$. Note the following two properties of these vectors, both immediate from concavity and monotonicity:
\begin{itemize}
    \item for any subset $I \subset \{1, \dots, k(\bar \alpha, \bar\beta)+1\}$ one has $\bar \alpha + \sum_{j\in I}e_j(\bar \alpha, \bar \beta)\in P$;
    \item $\bar \alpha + \sum_{j=1}^{k(\bar\alpha,\bar\beta)+1}e_j(\bar \alpha, \bar \beta)=\bar \beta$ and there is no proper nonempty $I \subset \{1, \dots, k( \bar \alpha, \bar \beta)+1\}$ such that $\bar \alpha + \sum_{j\in I}e_j(\bar \alpha, \bar \beta) \in \{\bar \alpha, \bar \beta \}$.
\end{itemize}
We will also write $\varkappa(\alpha,\beta)$ for the number of vertical segments, i.e.\ we set
\begin{equation}\label{dim2kappa}
\varkappa(\alpha,\beta)=k(\alpha,\beta)+1-(\beta_1-\alpha_1)=
\begin{cases}
\beta_2-\alpha_2,\text{ if }\alpha_2\le\beta_2\text{ and }(\beta_1,\alpha_2)\in P,\\
\beta_2-(\zeta_{\max}-\zeta_{\beta_1}),\text{ if }\alpha_2\le\beta_2\text{ and }(\beta_1,\alpha_2)\notin P,\\
\alpha_2-\beta_2-1,\text{ if }\alpha_2>\beta_2\text{ and }(\alpha_1,\beta_2)\in P,\\
\alpha_2-(\zeta_{\max}-\zeta_{\alpha_1}),\text{ if }\alpha_2>\beta_2\text{ and }(\alpha_1,\beta_2)\notin P.
\end{cases}
\end{equation}

\begin{prop}
\label{CurvetriangleDerivativeRelation}
The coefficients of the following series are zero in the ring $J_{\red}^{\infty}(R(P))$:
\begin{equation}\label{W(Y)}
W'_{\bar\alpha,\bar\beta,k'}=\sum_{I \subset \{1, \dots, k(\bar\alpha, \bar\beta)\}} (-1)^{|I|}\frac{\partial^{k'} Y_{\bar\alpha +\sum_{j \in I}e_j(\bar \alpha, \bar \beta)}(s)}{\partial s^{k'}}Y_{\bar \beta -\sum_{j \in I}e_j(\bar \alpha, \bar \beta)}(s),~ k'=0,\dots, k(\bar\alpha,\bar\beta)-1.
\end{equation}
\end{prop}
\begin{proof}
We define the map $[0,1]^{k(\bar \alpha, \bar \beta)+1} \rightarrow P$:
\[(x_1, \dots, x_{k(\bar \alpha, \bar \beta)+1})\mapsto \bar \alpha+
x_1 e_1(\bar \alpha, \bar \beta)+  \dots + x_{k(\bar \alpha, \bar \beta)+1}e_{k(\bar \alpha, \bar \beta)+1}(\bar \alpha, \bar \beta).\]

This map of polytopes defines the map of rings $\eta:R(\mathcal{B}_{k(\bar \alpha,\bar \beta)+1})\rightarrow R(P)$. Applying the corresponding arc map to $W_{k(\bar \alpha, \bar \beta)+1,k'} $ \eqref{Wdefinition} we have:
\[J^\infty(\eta)(W_{k(\bar \alpha, \bar \beta)+1,k'})=W'_{\bar\alpha,\bar\beta,k'}.\]
Thus the element $W'_{\bar\alpha,\bar\beta,k'}$ is nilpotent. This completes the proof.
\end{proof}

We define a lexicographic order on points $\bar\alpha \in P$:
\[\bar\alpha<\bar \beta~ \text{ if }~ \alpha_2<\beta_2~\text{or}~\alpha_2=\beta_2 ~\text{and}~ \alpha_1<\beta_1.\]
We define the following monomial order on the monomials in the variables $X_{\bar\alpha}, \bar \alpha \in P$. Let $X^{\bar r}:=\prod_{\bar \alpha \in P}X_{\bar \alpha}^{r_{\bar\alpha}}$. We define $r_{i}:=\sum_{j} r_{(i,j)}$, $L(r):=\sum_{\bar\alpha \in P} r_{\bar\alpha}$. If 
\[(L(\bar r), r_{0}, r_{1},\dots, r_{\eta})>(L(\bar r'), {r'}_{0}, {r'}_{1},\dots, {r'}_{\eta})\] 
in the standard lexicographic order, then $X^{\bar r} \succ X^{\bar r'}$. For $\bar r=(r_{\bar \alpha})_{\bar \alpha \in P}, \bar r'=(r'_{\bar \alpha})_{\bar \alpha \in P}$ such that
$(L(\bar r), r_{0}, r_{1},\dots, r_{\eta})=(L(\bar r'), {r'}_{0}, {r'}_{1},\dots, {r'}_{\eta})$ the order $\prec$ is lexicographic: we consider the minimal $\alpha$ with respect $<$ for which $r_{\bar\alpha}\neq r'_{\bar\alpha}$ and set $X^{\bar r}\prec X^{\bar r'}$ if $r_{\bar\alpha}< r'_{\bar\alpha}$.

As before, we have an induced order on the graded parts of the ring $J^{\infty}(\Bbbk[X_{\alpha}])$ and a filtration on the ring  $J^{\infty}(R(P))$.
Let $W'_{\bar\alpha,\bar\beta,k'}(X)$ be defined as the right-hand side of \eqref{W(Y)} with $Y_\bullet$ replaced by $X_\bullet$.
\begin{lem}
The initial part $\In_\prec W'_{\bar\alpha,\bar\beta,k'}(X)$ is equal to 
\[\frac{\partial^{k'}X_{\bar \alpha}(s)}{\partial s^{k'}}X_{\bar \beta}(s).\]
\end{lem}
\begin{proof}
We need to prove that $X_{\bar \alpha +\sum_{j \in I}e_j(\bar \alpha, \bar \beta)}X_{\bar \beta-\sum_{j \in I}e_j(\bar \alpha, \bar \beta)}\prec X_{\bar \alpha}X_{\bar \beta}$. It is clear from the definition of the order.
\end{proof}
\begin{cor}\label{twodimupperboud}
The dual of the graded component $\gr_{\prec}J^{\infty}_{\red}R(P)[\bar r]$ has a natural homogeneous $\mathfrak{h}[s]$-equivariant embedding into the space \[\prod_{\bar \alpha < \bar \beta \in P}\prod_{1 \leq i \leq r_{\bar \alpha}, 1 \leq j \leq r_{\bar \beta}}\left(t_{\bar\alpha}^{(i)}-t_{\bar\beta}^{(j)} \right)^{k(\bar \alpha, \bar \beta)}\Lambda_{\bar r}(\mathbf t).\]
\end{cor}

The rest of this section is dedicated to the proof of the following theorem.
\begin{thm}\label{dim2}
The coefficients of the series $W'_{\bar\alpha,\bar\beta,k'}(X)$ with distinct $\bar\alpha,\bar\beta\in P$ and $k'\in[0,k(\bar\alpha,\bar\beta)-1]$ generate the ideal of relations in $J^{\infty}_{\red}(R(P))$. Initial parts $\In_\prec W'_{\bar\alpha,\bar\beta,k'}(X)$ generate the corresponding initial ideal. 
\end{thm}
\begin{proof}
As before we study the dual map for the inclusion $\varphi_{\bar r}:A_{\bar r} \hookrightarrow J^\infty(\Bbbk[z_1,z_2,w])[\bar a,L]$, where $\bar a =\sum_{\bar \alpha \in P}r_{\bar\alpha}\bar \alpha$, $L=L(\bar r)$. Recall the notation $\bar r \in \mathcal{R}(\bar a,L)$.
The image of the map $\varphi^\vee_{\bar r}: \Lambda_{a_1,a_2,L}({\bf s})\rightarrow \Lambda_{\bar r}({\bf t})$,
\[s_1^{(\alpha_1^1 (j-1) +1)},\dots, s_1^{(\alpha_1^1 (j-1) + \alpha_1^1)},
s_2^{(\alpha_2^1 (j-1) +1)},\dots, s_2^{(\alpha_2^1 (j-1) + \alpha_2^1)},s_w^{(j)}\mapsto t_{\bar \alpha^1}^{(j)},~j=1, \dots, r_{\bar \alpha^1};\]
\begin{multline*}
s_1^{(r_{\bar \alpha^1}\alpha_1^1+  \alpha_1^2(j-1) +1)},\dots, s_1^{(r_{\bar \alpha^1}\alpha_1^1+\alpha_1^2 (j-1) + \alpha_1^2)},\\
s_2^{(r_{\bar \alpha^1}\alpha_2^1+\alpha_2^2 (j-1) +1)},\dots, s_2^{(r_{\bar \alpha^1}\alpha_2^1+\alpha_2^2 (j-1) + \alpha_2^2)},s_w^{(r_{\bar \alpha^1}+j)}\mapsto t_{\bar \alpha^2}^{(j)},~j=1, \dots, r_{\bar \alpha^2};
\end{multline*}
\[\dots\]
is isomorphic to the dual space of $A_{\bar r}$.
The dual space to $\sum_{\bar r' \preceq \bar r}A_{\bar r}/\sum_{\bar r' \prec \bar r}A_{\bar r}$ is isomorphic to the space $\varphi_{\bar r}^{\vee}(\bigcap_{\bar r' \in \mathcal{R}(\bar a,L),\bar r' \prec \bar r} \ker(\varphi_{\bar r'}^{\vee}))$. Our goal is to compute this space.

We consider the set of variables $u_i^{(j)}, i=0, \dots, \eta, j=1, \dots, r_i$. Consider the following map $\varphi^{\vee}_{(r_0, \dots, r_\eta)}\otimes \id:\Lambda_{a_1,L}({\bf s})\otimes\Lambda_{a_2}({\bf s_2})=\Lambda_{a_1,a_2,L}({\bf s})\rightarrow \Lambda_{(r_0, \dots, r_\eta,a_2)}({\bf u,s_2})$ (see Subsection \ref{Veronese}):
\begin{align*}
s_w^{(j)} &\mapsto u_0^{(j)}, j=1, \dots, r_0; \\
s_1^{(j)},s_w^{(r_0+j)} &\mapsto u_1^{(r_0+j)}, j=1, \dots, r_1; \\
s_1^{(r_1+2j-1)},s_1^{(r_1+2j)},s_w^{(r_0+r_1+j)} &\mapsto u_1^{(r_0+r_1+j)}, j=1, \dots, r_2;\dots;  &\hspace{-15mm}s_2^{(i)} \mapsto s_2^{(i)}.
\end{align*}
We have an order on vectors $\pr_2(\bar r):=(r_0,\dots, r_\eta)$ as in the previous section. By definition we have: 
\begin{itemize}
    \item if $\bar r' \prec \bar r$, then 
$\pr_2(\bar r')\leq \pr_2(\bar r)$;
    \item the map $\varphi^{\vee}_{\bar r}$ is right divisible by $\varphi^{\vee}_{\pr_2(\bar r)}$. 
\end{itemize}

More precisely we define the following map $\tilde\psi_{\bar r}: \Lambda_{(r_0, \dots, r_\eta,a_2)}({\bf u,s_2}) \rightarrow \Lambda_{\bar r}({\bf t})$,
\begin{equation*}
u_{i}^{(j)},s_2^{((\zeta_{\max}-\zeta_i)(j-1)+1)}, \dots, s_2^{((\zeta_{\max}-\zeta_i)j)}\mapsto t_{i,(\zeta_{\max}-\zeta_i)}^{(j)}, j=1, \dots, r_{i,\zeta_{\max}-\zeta_i};
\end{equation*}
\begin{multline*}
u_{i}^{(r_{i,\zeta_{\max}-\zeta_i}+j)},s_2^{(r_{i,\zeta_0-\zeta_i}(\zeta_{\max}-\zeta_i)+(\zeta_{\max}-\zeta_i+1)(j-1)+1)}, \dots, s_2^{(r_{i,\zeta_{\max}-\zeta_i}(\zeta_{\max}-\zeta_i)+(\zeta_{\max}-\zeta_i+1)j)}\mapsto\\ t_{i,(\zeta_{\max}-\zeta_i+1)}^{(j)}, j=1, \dots, r_{i,\zeta_{\max}-\zeta_i+1};\dots
\end{multline*}
Then:
\[\varphi^\vee_{\bar r}=\tilde\psi_{\bar r }\circ (\varphi^{\vee}_{\pr_2(\bar r)}\otimes \id).\]
Thus we have:
\begin{multline*}\varphi^{\vee}_{\pr_2(\bar r)}\otimes \id\left(\bigcap_{\bar r'\prec \bar r}\ker \varphi^{\vee}_{\bar r'}\right) \supset \\\prod_{1 \leq i<j\leq \eta}\prod_{1 \leq l \leq r_i,1 \leq l' \leq r_j}\left(u_{i}^{(l)}-u_{j}^{(l')} \right)^{j-i-1}
\Lambda_{(r_0,\dots,r_\eta,a_2)}({\mathbf u,\mathbf s_{2}})\cap  \bigcap_{\bar r' \prec \bar r, \pr_{2}(\bar r')=(r_0,\dots,r_\eta)}\ker (\tilde\psi_{\bar r'})\supset\\
\prod_{1 \leq i<j\leq \eta}\prod_{1 \leq l \leq r_i,1 \leq l'\leq r_j}\left(u_{i}^{(l)}-u_{j}^{(l')} \right)^{j-i-1}
\left(\bigcap_{\bar r' \prec  \bar r, \pr_{2}(\bar r')=(r_0,\dots,r_\eta)}\ker (\tilde\psi_{\bar r'})\right).
\end{multline*}
Therefore
\[\varphi^{\vee}_{\bar r}\left( \bigcap_{r' < r}\ker \varphi^{\vee}_{\bar\rho'} \right)\supset
\tilde\psi_{\bar r}\left(\prod_{0 \leq a<b\leq \eta, 1 \leq l \leq r_{i}, 1 \leq l' \leq r_{j} }\left(u_{i}^{(l)}-u_{j}^{(l')} \right)^{b-a-1}\right)
\tilde\psi_{\bar r}\left(\Lambda_{(r_0, \dots, r_\eta,a_2)}({\bf u,s_2})\right).\]

Consider now the set of $\bar r'\prec \bar r$ such that $\pr_2(\bar r')=\bar\rho$, where $\bar\rho=\pr_2(\bar r)$.
 One has:
\begin{multline}\label{twofactors}
\varphi^{\vee}_{\bar r}\left( \bigcap_{\bar r' \prec \bar r, \pr_2(\bar r')=\bar\rho}\ker \varphi^{\vee}_{\bar\rho'} \right)\supset\\\tilde\psi_{\bar r}\left(\prod_{0 \leq a<b\leq \eta, 1 \leq i \leq r_{a}, 1 \leq j \leq r_{b} }\left(u_{a}^{(i)}-u_{b}^{(j)} \right)^{b-a-1}\right)
\tilde\psi_{\bar r}\left( \bigcap_{\bar r' \prec \bar r, \pr_2(\bar r')=\bar\rho}\ker (\tilde\psi_{\bar r'}) \right).
\end{multline}

Let us study the two factors in \eqref{twofactors} separately. For the first factor one has
\begin{multline*}
\tilde\psi_{\bar r}\left(\prod_{0 \leq a<b\leq \eta, 1 \leq i \leq r_{a}, 1 \leq j \leq r_{b} }\left(u_{a}^{(i)}-u_{b}^{(j)} \right)^{b-a-1}\right)\\=
\prod_{\bar \alpha,\bar \beta \in P, \alpha_1<\beta_1}
\prod_{1 \leq l\leq r_{\bar\alpha},1 \leq l'\leq r_{\bar \beta}}
\left( t_{\bar\alpha}^{(l)}-t_{\bar\beta}^{(l')} \right)^{\beta_1-\alpha_1-1}.
\end{multline*}

Now let us compute the second factor. We define the following map:
\[g_{\bar\rho}:\Lambda_{(\bar\rho,a_2)}({\bf u,s_2})\rightarrow \Lambda_{(\bar\rho,a_2-\sum_{i=0}^\eta r_i(\zeta_{\max}-\zeta_i))}({\bf u,s_2}),\]
\[s_2^{(\sum_{i'=0}^{i-1} r_{i'}(\zeta_{\max}-\zeta_{i'})+(l-1)(\zeta_{\max}-\zeta_{i})+1)}, \dots, s_2^{(\sum_{i'=0}^{i-1} r_{i'}(\zeta_{\max}-\zeta_{i'})+l(\zeta_{\max}-\zeta_{i}))},u_i^{(l)}\mapsto u_i^{(l)},\]
\[s_2^{(l)}\mapsto s_2^{(l-\sum_{i=0}^\eta r_i(\zeta_{\max}-\zeta_i))},~l \geq \sum_{i=0}^\eta r_i(\zeta_{\max}-\zeta_i).\]

Then we have:
\[\tilde \psi_{\bar r}=\psi_{d(\bar r)} \circ g_{\bar\rho},\]
where $d(\bar r)=d(\bar r)_{i,j}$, $i=0,\dots,\eta$, $j=0,\dots,\zeta_i$,
$d(\bar r)_{i,j}=r_{i,j+\zeta_{\max}-\zeta_i}$ and $\psi_{d(\bar r)} $ is defined in the following way (it differs from the $\psi_{d(\bar r)}$ defined in Subsection \ref{TechnicalLemma} only by a change of the subscripts):
\begin{equation*}
u_{i}^{(j)}\mapsto t_{i,(\zeta_{\max}-\zeta_i)}^{(j)}, j=1, \dots, r_{i,\zeta_{\max}-\zeta_i};
\end{equation*}
\begin{equation*}
u_{i}^{(r_{i,\zeta_{\max}-\zeta_i}+j)},s_2^{(j)}\mapsto\\ t_{i,(\zeta_{\max}-\zeta_i+1)}^{(j)}, j=1, \dots, r_{i,\zeta_{\max}-\zeta_i+1};
\end{equation*}
\begin{multline*}
u_{i}^{(r_{i,\zeta_{\max}-\zeta_i}+r_{i,\zeta_{\max}-\zeta_i+1}+j)},s_2^{(r_{i,\zeta_{\max}-\zeta_i+1}+2j-1)},s_2^{(r_{i,\zeta_{\max}-\zeta_i+1}+2j)}\\\mapsto t_{i,(\zeta_{\max}-\zeta_i+2)}^{(j)}, j=1, \dots, r_{i,\zeta_{\max}-\zeta_i+2};\dots
\end{multline*}

The map $g_{\bar\rho}$ is surjective. Thus we have
\[\tilde\psi_{\bar r}\left( \bigcap_{\bar r' \prec \bar r, \pr_2(\bar r')=\rho}\ker (\tilde\psi_{\bar r'}) \right)=
\psi_{d(\bar r)}\left( \bigcap_{\bar r' \succ \bar r, \pr_2(\bar r')=\rho}\ker (\psi_{d(\bar r')}) \right).\]

Note that the order on vectors $\bar r$ gives the lexicographic order on vectors $d(\bar r)$ with respect to the following order on coordinates: $(i,j) \prec(i',j')$ iff $j+ (\zeta_{\max}-\zeta_i)<j'+ (\zeta_{\max}-\zeta_{i'})$ or  $j+ (\zeta_{\max}-\zeta_i)=j'+ (\zeta_{\max}-\zeta_{i'})$ and $i<i'$.
So using Lemma \ref{Split} we have:
\[\psi_{d(\bar r)}\left( \bigcap_{\bar r' \prec \bar r, \pr_2(\bar r')=\bar\rho}\ker (\psi_{d(\bar r')}) \right)\supset \prod_{\bar \alpha<\bar \beta}\prod_{1 \leq l \leq r_{\bar \alpha},1 \leq l' \leq r_{\bar \beta}}\left( t_{\bar\alpha}^{(l)}-t_{\bar \beta}^{(l)} \right)^{\varkappa(\bar \alpha,\bar \beta)} \Lambda_{\bar r}({\bf t})\]
(see \eqref{dim2kappa}).
Multiplying the two factors of \eqref{twofactors} we obtain 
\[\varphi_{\bar r}^{\vee}\left(\bigcap_{r' \prec r} \ker(\varphi_{\bar r'}^{\vee})\right)\supset \prod_{\bar\alpha<\bar\beta}\prod_{1 \leq l \leq r_{\bar \alpha},1 \leq l' \leq r_{\bar \beta}}\left( t_{\bar\alpha}^{(l)}-t_{\bar\beta}^{(l')} \right)^{k(\bar\alpha,\bar\beta)} \Lambda_{\bar r}({\bf t}).\]

Therefore, due to Corollary \ref{twodimupperboud} these spaces are equal.  Thus  $\In_\prec W'_{\bar\alpha,\bar\beta,k'}(X)$ generate the defining ideal of $\gr_{\prec}J^{\infty}_{\red}(R(P))$. This completes the proof of the theorem.
\end{proof}
\begin{cor}
The action of $\mathcal{A}_L$ is free on $J_{\red}^{\infty}(R(P)[L]^*$.
\end{cor}

\subsection{Higher-dimensional case}
In this subsection we study a family of polytopes such that the (dual of the) reduced arc rings of the corresponding toric varieties admit a free action of the polynomial algebras. 

Consider a convex lattice polytope $P \subset \mathbb{R}^n$. Let $\bar \alpha^1,\dots,\bar \alpha^m$ be the tuple of integer points of $P$, $\prec$ be a monomial order on $\Bbbk[X_{\bar \alpha^1},\dots,X_{\bar \alpha^m}]$, $\gamma(\bar \alpha^i,\bar \alpha^j)=\gamma(\bar \alpha^i,\bar \alpha^j)$ with $1\le i,j\le m$ be a collection of nonnegative integers. Assume that for any $i\ne j$ there exists a set of vectors $e_1(\bar \alpha^i,\bar \alpha^j), \dots, e_{\gamma(\bar \alpha^i,\bar \alpha^j)+1}(\bar \alpha^i,\bar \alpha^j)$ with the following properties:
\begin{itemize}\label{eConstruction}
    \item $\bar \alpha^i+\sum_{l=1}^{\gamma(\bar \alpha^i,\bar \alpha^j)+1}e_l(\bar \alpha^i,\bar \alpha^j)=\bar \alpha^j$;
    \item for any $J\subset \{1, \dots,\gamma(\bar \alpha^i,\bar \alpha^j)+1\}$: $\bar \alpha^i+\sum_{l\in J}e_l(\bar \alpha^i,\bar \alpha^j)\in P$;
    \item for any nonempty proper subset $J\subset \{1, \dots,\gamma(\bar \alpha^i,\bar \alpha^j)+1\}$:
    $$X_{\bar \alpha^i}X_{\bar \alpha^j}\prec X_{\bar \alpha^i+\sum_{l\in J}e_l(\bar \alpha^i,\bar \alpha^j)}X_{\bar \alpha^j-\sum_{l\in J}e_l(\bar \alpha^i,\bar \alpha^j)};$$
    \item $e_l(\bar \alpha^j,\bar \alpha^i)=-e_l(\bar \alpha^i,\bar \alpha^j)$.
\end{itemize}
We call $(P,\prec,\gamma,\{e_l\})$ a {\it cube generating data}.

We define the linear maps of polytopes
$
\eta_{\bar \alpha^i,\bar \alpha^j}:[0,1]^{\gamma(\bar \alpha^i,\bar \alpha^j)+1}\rightarrow P$: 
\begin{equation} \label{etaDefinition}
\eta_{\bar \alpha^i,\bar \alpha^j}\left(\sum_{l=1}^{\gamma(\bar \alpha^i,\bar \alpha^j)+1}x_le_l\right)=\bar \alpha^i+\sum_{l=1}^{\gamma(\bar \alpha^i,\bar \alpha^j)+1} x_le_l(\bar \alpha^i,\bar \alpha^j),
\end{equation}
where $\{e_l\}$ is the basis vectors of $\mathbb{R}^{\gamma(\bar \alpha^i,\bar \alpha^j)+1}$, $0\leq x_l \leq 1$.
The corresponding map of rings $R([0,1]^{\gamma(\bar \alpha^i,\bar \alpha^j)+1})\rightarrow R(P)$ will also be denoted by $\eta_{\bar \alpha^i,\bar \alpha^j}$.

For any $k \leq \gamma(\bar \alpha^i,\bar \alpha^j)-1$
\[J^{\infty}(\eta_{\bar \alpha^i,\bar \alpha^j})\left( W_{\gamma(\bar \alpha^i,\bar \alpha^j)+1,k} \right)\]
is nilpotent. Therefore the coefficients of the series
\[W_{\bar \alpha^i,\bar \alpha^j,k}(X):=\sum_{I \subset \{2,\dots,\gamma(\bar \alpha^i,\bar \alpha^j)+1\}}(-1)^{|I|}X_{\bar \alpha^i+\sum_{l \in I\cup \{1\}}e_l(\bar \alpha^i,\bar \alpha^j)}(s)\frac{\partial^k X_{\bar \alpha^j-\sum_{l \in I\cup \{1\}}e_l(\bar \alpha^i,\bar \alpha^j)}(s)}{\partial s^k}\]
belong to the defining ideal of $J_{\red}^{\infty}(R(P))$ for $k\le \gamma(\bar \alpha^i,\bar \alpha^j)-1$.
By construction we have
\begin{equation}\label{reducedEquatuinAssumption}
\In_\prec W_{\bar \alpha^i,\bar \alpha^j,k}(X)=X_{\bar \alpha^i}(s)\frac{\partial^k X_{\bar \alpha^j}(s)}{\partial s^k}.
\end{equation}
Recall that in view of Proposition~\ref{UpperBoundInitial} the coefficients of  \eqref{reducedEquatuinAssumption} belong to the defining ideal of the associated graded ring $\gr_{\prec}J_{\red}^{\infty}(R(P))$. Suppose that, moreover, they generate this ideal. 
In this case 
we call the data $(P,\prec~,\gamma,\{e_l\})$ a {\it strict cube generating data}.
\begin{prop}
Consider a strict cube generating data $(P,\prec,\gamma(\bar \alpha^i,\bar \alpha^j),\{e_l(\bar \alpha^i,\bar \alpha^j)\})$. Then for any $\bar\rho\in\bN^m$ we have the following isomorphism of $\grad$-graded spaces:
\[\gr_{\prec}J_{\red}^{\infty}(R(P))[\bar \rho]^*\simeq
\prod_{i<j}\prod_{\substack{1 \leq l \leq \rho_i\\ 1 \leq l' \leq \rho_j}}\left(  t_i^{(l)}-t_j^{(l')}\right)^{\gamma(\bar \alpha^i,\bar \alpha^j)}\Lambda_{\bar \rho}({\bf t}).\]
\end{prop}
\begin{proof}
This follows from Proposition \ref{DualSpaceDerivedQuagraticMon}.
\end{proof}
\begin{cor}
For $L\in\mathbb N$ the action of $\mathcal{A}_L$ is free on $J_{\red}^{\infty}(R(P))[L]^*$.
\end{cor}
\begin{proof}
It follows from Theorem \ref{equalityANDcofree}.
\end{proof}

Consider a convex function $\zeta: P\rightarrow \mathbb R_{\ge0}$ with $\zeta(\bar \alpha^i)=\zeta_i\in\bZ$, set $\zeta_{\max}=\max\{\zeta_i\}$.
Denote by $P^{\zeta}$ the polytope in $\mathbb{R}^{n+1}$ obtained as the convex hull of $P\times \zeta_{\max}$ and all $\bar \alpha^i\times(\zeta_{\max}-\zeta_i)$. In what follows we construct a cube generating data for $P^\zeta$ starting from a cube generating data for $P$ satisfying the following condition. 
For any $\bar \alpha^i,\bar \alpha^j \in P$ such that $\zeta_i\leq \zeta_j$, there exist  nonnegative integers $f_1(\bar \alpha^i,\bar \alpha^j),\dots, f_{\gamma(\bar \alpha^i,\bar \alpha^j)+1}(\bar \alpha^i,\bar \alpha^j)$ such that
\begin{itemize}
    \item $\zeta_i+\sum_{l=1}^{\gamma(\bar \alpha^i,\bar \alpha^j)+1} f_l(\bar \alpha^i,\bar \alpha^j)=\zeta_j;$
    \item for any $J\subset \{1, \dots,\gamma(\bar \alpha^i,\bar \alpha^j)+1\}$:
    $\zeta_i+\sum_{l \in J}f_l(\bar \alpha^i,\bar \alpha^j)\le \zeta(\bar \alpha^i+\sum_{l \in J}e_l(\bar \alpha^i,\bar \alpha^j))$.
\end{itemize}

\begin{example}\label{ConvexFunctionExample}
Consider the cube generating data $(P,\prec,\gamma,\{e_l\})$. Assume that for any $\bar \alpha^i,\bar \alpha^j$ the $l$th coordinate of every $e_p(\bar \alpha^i,\bar \alpha^j)$ is either nonnegative or nonpositive. Then any convex function of the form $\zeta(\bar\alpha)=\zeta(\alpha_l)$, i.e.\ depending only on the $l$th coordinate admits such a set of integers $f_p(
\bar \alpha^i,\bar \alpha^j)$. Indeed, without loss of generality we may assume that $\zeta(\bar\alpha^i) \geq  \zeta(\bar\alpha^j)$. We may consider an enumeration of vectors $e_p(\bar \alpha^i,\bar \alpha^j)$ such that the absolute value of the $l$th coordinate is weakly increasing. 
Let $x$ be the smallest number such that $\zeta(\bar \alpha^i+\sum_{p=1}^x e_p(\bar \alpha^i,\bar \alpha^j))<\zeta(\bar \alpha^i)$. We put 
\begin{align*}&f_q(\bar \alpha^i,\bar \alpha^j)=0, 1 \leq q <x;\\&
f_x(\bar \alpha^i,\bar \alpha^j)=\zeta(\bar \alpha^i+\sum_{p=1}^x e_p(\bar \alpha^i,\bar \alpha^j))-\zeta(\bar \alpha^i);\\
&f_q(\bar \alpha^i,\bar \alpha^j)=\zeta(\bar \alpha^i+\sum_{p=1}^q e_p(\bar \alpha^i,\bar \alpha^j))-\zeta(\bar \alpha^i+\sum_{p=1}^{q-1} e_p(\bar \alpha^i,\bar \alpha^j)),x<q\leq \gamma(\bar \alpha^i,\bar \alpha^j) +1.
\end{align*}
Note that the numbers $f_q(\bar \alpha^i,\bar \alpha^j)$ are all either nonnegative or nonpositive.
\end{example}

For $\bar r =(r_{(\bar\alpha^i,a)})_{(\bar\alpha^i,a)\in P^{\zeta}}$ we define 
$\pr_{n+1}(\bar r)=(\pr_{n+1}(\bar r))_{\bar\alpha^i}$, $i=1,\dots,m$ by the formula
\[(\pr_{n+1}(\bar r))_{\bar\alpha^i}=\sum_{a=\zeta_{\max}-\zeta_i}^{\zeta_{\max}}r_{(\bar\alpha^i,a)}.\]

For $(\bar\alpha^i,a), (\bar\alpha^j,b)\in P^{\zeta}\cap \mathbb{Z}^{n+1}$ we let
\begin{equation}\label{order}
(\bar\alpha^i,a)< (\bar\alpha^j,b) \text{ if } a<b \text{ or } a=b \text{ and } i<j.
\end{equation}
We define an order $\prec^\zeta$ on monomials in $\Bbbk[X_{(\bar\alpha^i,a)}]$, $(\bar\alpha^i,a) \in P^\zeta \cap \mathbb{Z}^{n+1}$ in the following way.
For $r,r'\in\mathbb N^{P^\zeta\cap\bZ^{n+1}}$ with $\pr_{n+1}(\bar r) \prec \pr_{n+1}(\bar r')$ we set $\bar r \prec^\zeta \bar r'$ and $X^{\bar r}\prec^\zeta X^{\bar r'}$. If $\pr_{n+1}(\bar r) = \pr_{n+1}(\bar r')$, then $\prec^\zeta$ compares $r$ to $r'$ and $X^{\bar r}$ to  $X^{\bar r'}$ lexicographically with respect to the order \eqref{order}.

Suppose that $(\bar \alpha^i,a)\leq (\bar \alpha^j,b)$. We set:
\begin{equation*}
    \gamma((\bar\alpha^i,a),(\bar \alpha^j,b))=\gamma(\bar \alpha^i,
    \bar \alpha^j)+\varkappa((\bar\alpha^i,a),(\bar \alpha^j,b)),
\end{equation*}
where
\begin{equation}\label{kappahigherdim}
    \varkappa((\bar\alpha^i,a),(\bar \alpha^j,b))=\begin{cases}
    b-(\zeta_{\max}-\zeta_j),\text{if}~ \zeta_{\max}-\zeta_j>a;\\
    b-a, \text{if}~ \zeta_{\max}-\zeta_j\leq a ~\text{and}~ i \leq j;\\
    b-a-1, \text{if}~  \zeta_{\max}-\zeta_j\leq a ~\text{and}~ i>j.
    \end{cases}
\end{equation}
We construct the vectors $e_l((\bar\alpha^i,a),(\bar \alpha^j,b))$ in the following way.

Suppose that $\zeta_{\max}-\zeta_j>a$. 
Let $y$ be the number such that $\sum_{l=1}^{y-1}f_l(\bar\alpha^i,\bar\alpha^j)<a-(\zeta_{\max}-\zeta_i)\leq \sum_{l=1}^{y}f_l(\bar\alpha^i,\bar\alpha^j)$.
Then 
\begin{align*}
&e_l((\bar\alpha^i,a),(\bar \alpha^j,b))=(e_l(\bar \alpha^i,\bar \alpha^j),0), l=1,\dots,y-1,\\
&e_y=((\bar\alpha^i,a),(\bar \alpha^j,b))=(e_y(\bar \alpha^i,\bar \alpha^j),f_y(\bar \alpha^i,\bar \alpha^j)-(\zeta_{\max}-\zeta_i)),\\
&e_l((\bar\alpha^i,a),(\bar \alpha^j,b))=(e_l(\bar \alpha^i,\bar \alpha^j),f_l(\bar \alpha^i,\bar \alpha^j)), l=y+1,\dots,\gamma(\bar \alpha^i,\bar \alpha^j)+1,\\
&e_l((\bar\alpha^i,a),(\bar \alpha^j,b))=(0,1), l=\gamma(\bar \alpha^i,\bar \alpha^j)+2, \dots, \gamma((\bar\alpha^i,a),(\bar \alpha^j,b))+1. 
\end{align*}
If $a\geq\zeta_{\max}-\zeta_j$ and $i\leq j$,
then 
\begin{align*}
&e_l((\bar\alpha^i,a),(\bar \alpha^j,b))=(e_l(\bar \alpha^i,\bar \alpha^j),0), l=1,\dots,\gamma(\bar \alpha^i,\bar \alpha^j)+1,\\ 
&e_l((\bar\alpha^i,a),(\bar \alpha^j,b))=(0,1), l=\gamma(\bar \alpha^i,\bar \alpha^j)+2, \dots, \gamma(\bar \alpha^i,\bar \alpha^j)+b-a+1.    
\end{align*} 
If $a\geq\zeta_{\max}-\zeta_j$ and $i>j$, then 
\begin{align*}
&e_l((\bar\alpha^i,a),(\bar \alpha^j,b))=(e_l(\bar \alpha^i,\bar \alpha^j),0), l=1,\dots,\gamma(\bar \alpha^i,\bar \alpha^j),\\
&e_{\gamma(\bar \alpha^i,\bar \alpha^j)+1}((\bar\alpha^i,a),(\bar \alpha^j,b))=
(e_{\gamma(\bar \alpha^i,\bar \alpha^j)+1},1),\\
&e_l((\bar\alpha^i,a),(\bar \alpha^j,b))=(0,1), l=\gamma(\bar \alpha^i,\bar \alpha^j)+2, \dots, \gamma(\bar \alpha^i,\bar \alpha^j)+b-a.    
\end{align*}
For $a>b$ we put $e_l((\bar\alpha^i,a),(\bar \alpha^j,b))=-e_l((\bar \alpha^j,b),(\bar\alpha^i,a))$. 
The following lemma follows directly from the definitions.
\begin{lem}
$(P^{\zeta},\prec^\zeta, \gamma,\{e_l\})$ satisfy the conditions for cube generating data.
\end{lem}

We consider the maps 
\[\eta_{(\bar\alpha^i,a),(\bar\alpha^j,b)}:[0,1]^{\gamma((\bar\alpha^i,a),(\bar\alpha^j,b))+1}\rightarrow P^\zeta\]
corresponding to this data. By construction and Proposition \ref{DualSpaceQuagraticMon} we have 

\begin{equation}\label{inclusionToIdeal}
(\gr_{\prec}J_{\red}^{\infty}(R(P^{\zeta}))[\bar r])^*\hookrightarrow \prod_{(\bar\alpha^i,a)<(\bar\alpha^j,b)}\prod_{1 \leq l \leq r_{(\bar\alpha^i,a)},1\leq l' \leq r_{(\bar\alpha^j,b)}}\left(t_{(\bar\alpha^i,a)}^{(l)}-t_{(\bar\alpha^j,b)}^{(l')} \right)^{\gamma((\bar\alpha^i,a),(\bar\alpha^j,b))}\Lambda_{\bar r}({\bf t}).
\end{equation}
We now show that this inclusion is, in fact, an equality.

\begin{thm}Assume that \label{main} 
$(P,\prec, \gamma,\{e_l\})$ is a strict cube generating data. 
Then the data
$(P^{\zeta},\prec^\zeta~,\gamma,\{e_{l}\})$ is a strict cube generating data.
\end{thm}
\begin{proof}
In this proof we compute the space $(\gr_{\prec}J_{\red}^{\infty}(R(P^{\zeta}))[\bar r])^*$ in the same way as in the proof of Theorem \ref{dim2}.
As before we study the dual map for the inclusion $\varphi_{\bar r}:A_{\bar r} \hookrightarrow J^\infty(\Bbbk[z_1,z_2,\dots,z_{n+1},w])[(\bar a,a_{n+1}),L]$, where $(\bar a,a_{n+1}) =\sum_{(\bar \alpha,\alpha_{n+1}) \in P^{\zeta}}r_{(\bar\alpha,\alpha_{n+1})}(\bar \alpha,\alpha_{n+1})$, $L=\sum_{(\bar \alpha,\alpha_{n+1}) \in P^{\zeta}}r_{(\bar\alpha,\alpha_{n+1})}$.
The image of the map $\varphi^\vee_{\bar r}: \Lambda_{(\bar a,a_{n+1}),L}({\bf s})\rightarrow \Lambda_{\bar r}({\bf t})$,
\begin{multline*}s_1^{(\alpha_1^1 (j-1) +1)},\dots, s_1^{(\alpha_1^1 (j-1) + \alpha_1^1)},\dots,
s_{n+1}^{((\zeta_{\max}-\zeta_1) (j-1) +1)},\dots, s_{n+1}^{((\zeta_{\max}-\zeta_1) j )},s_w^{(j)}\\\mapsto t_{(\bar \alpha^1,\zeta_{\max}-\zeta_i)}^{(j)},~j=1, \dots, r_{(\bar \alpha^1,\zeta_{\max}-\zeta_1)};
\end{multline*}
\begin{multline*}s_1^{(r_{(\bar \alpha^1,\zeta_{\max}-\zeta_1)} \alpha_1^1+ \alpha_1^1 (j-1) +1)},\dots, s_{n+1}^{(r_{(\bar \alpha^1,\zeta_{\max}-\zeta_1)}(\zeta_{\max}-\zeta_i)+(\zeta_{\max}-\zeta_1+1) j)},s_w^{(j)}\\\mapsto t_{(\bar \alpha^1,\zeta_{\max}-\zeta_i)}^{(j)},~j=r_{(\bar \alpha^1,\zeta_{\max}-\zeta_1)}+1, \dots, r_{(\bar \alpha^1,\zeta_{\max}-\zeta_1)}+r_{(\bar \alpha^1,\zeta_{\max}-\zeta_1-1)};
\end{multline*}
\[\dots\]
is isomorphic to the dual space of $A_{\bar r}$. Here for every integer point $(\bar\alpha,\alpha_{n+1})\in P^\zeta$ and $1\le i\le n+1$ the map $\varphi^\vee_{\bar r}$ takes $r_{(\bar\alpha,\alpha_{n+1})}\alpha_i$ variables of the form $s_i^\bullet$ to variables of the form $t_{(\bar\alpha,\alpha_{n+1})}^\bullet$, it also takes $r_{(\bar\alpha,\alpha_{n+1})}$ of the $s_w^\bullet$ to the  $t_{(\bar\alpha,\alpha_{n+1})}^\bullet$.
The dual space to $A_{\bar r}\left/\left(A_{\bar r}\cap(\oplus_{r' \prec r}A_{\bar r})\right)\right.$ is isomorphic to the space $\varphi_{\bar r}^{\vee}(\cap_{r' \prec r} \ker(\varphi_{\bar r'}^{\vee}))$. Our goal is to compute this space.

Denote $\bar\rho:=\pr_{n+1}(\bar r)$.
We consider the set of variables $u_i^{(j)}=u_{\alpha^i}^{(j)}, i\in P, j=1, \dots, \rho_i$. Consider the map $\varphi^{\vee}_{\bar \rho}\otimes \id:\Lambda_{(\bar a,a_{n+1}),L}({\bf s})\rightarrow \Lambda_{{(\bar \rho,a_{n+1})}}({\mathbf u,\mathbf s_{n+1}})$, i.e.\ $s_{n+1}^{(i)} \mapsto s_{n+1}^{(i)}$ and it is equal to $\varphi^{\vee}_{\bar \rho}$ on $s_1^{(i)},\dots,s_n^{(i)},s_w^{(i)}$ (but with $t$ replaced by $u$).

By definition we have: 
\begin{itemize}
    \item if $\bar r' \prec \bar r$, then 
$\pr_{n+1}(\bar r')\preceq \pr_{n+1}(\bar r)$;
    \item the map $\varphi^{\vee}_{\bar r}$ is right divisible by $\varphi^{\vee}_{\rho}$.
\end{itemize}
More precisely, we define the following map $\tilde\psi_{\bar r}: \Lambda_{(\bar \rho,a_{n+1})}({\mathbf u,\mathbf s_{n+1}}) \rightarrow \Lambda_{\bar r}({\bf t})$,
\begin{align*}
u_{i}^{(j)},s_{n+1}^{((\zeta_{\max}-\zeta_i)(j-1)+1)}, \dots, s_{n+1}^{((\zeta_{\max}-\zeta_i)j)}\mapsto t_{(\bar\alpha^i,(\zeta_{\max}-\zeta_i))}^{(j)},& j=1, \dots, r_{i,\zeta_{\max}-\zeta_i};\\
u_{i}^{(r_{i,\zeta_{\max}-\zeta_i}+j)},s_{n+1}^{(r_{i,\zeta_{\max}-\zeta_i}(\zeta_{\max}-\zeta_i)+(\zeta_{\max}-\zeta_i+1)(j-1)+1)}, &\dots, s_{n+1}^{(r_{i,\zeta_{\max}-\zeta_i}(\zeta_{\max}-\zeta_i)+(\zeta_{\max}-\zeta_i+1)j)}\mapsto\\ t_{(\bar\alpha^i(\zeta_{\max}-\zeta_i+1))}^{(j)}, j=1, &\dots, r_{i,\zeta_{\max}-\zeta_i+1};\\&\dots
\end{align*}
Here for every $(\bar\alpha_i,\alpha_{n+1})\in P^\zeta$ we have $r_{(\bar\alpha_i,a_{n+1})}$ variables of the form $u_i^\bullet$ and $\alpha_{n+1}$ variables of the form $s_{n+1}^\bullet$ being mapped to variables of the form $t_{(\bar\alpha_i,\alpha_{n+1})}^\bullet$.
Then $\varphi^\vee_{\bar r}=\tilde\psi_{\bar r }\circ (\varphi^{\vee}_{\rho}\otimes \id)$.

By the definition of strict cube generating data and Proposition~\ref{UpperBoundInitial}, the right-hand sides of~\eqref{reducedEquatuinAssumption} generate the ring $\gr_{\prec}J_{\red}^{\infty}(R(P))$. A component of the latter ring is considered in Corollary \ref{DualToSubauotient} and Proposition \ref{DualSpaceDerivedQuagraticMon} provides:
\[\varphi^{\vee}_{\bar\rho}\otimes \id\left( \bigcap_{\bar\rho' \prec \rho}\ker( \varphi^{\vee}_{\bar\rho'}\otimes \id) \right)=
\prod_{1 \leq i<j\leq m}\prod_{1 \leq l \leq \rho_i,1 \leq l' \leq \rho_j}\left(u_{\bar\alpha^i}^{(l)}-u_{\bar\alpha^j}^{(l')} \right)^{\gamma(\bar\alpha^i,\bar\alpha^j)}
\Lambda_{(\rho,\alpha_{n+1})}({\mathbf u,\mathbf s_{n+1}}).\]
Consequently, we have:
\[\varphi^{\vee}_{\rho}\otimes \id\left( \bigcap_{\pr_{n+1}(\bar r') \prec \bar\rho}\ker \varphi^{\vee}_{\bar r'} \right)\supset
\prod_{1 \leq i<j\leq m}\prod_{1 \leq l \leq \rho_i,1 \leq l' \leq \rho_j}\left(u_{\bar\alpha_i}^{(l)}-u_{\bar\alpha_j}^{(l')} \right)^{\gamma(\bar\alpha^i,\bar\alpha^j)}
\Lambda_{(\rho,\alpha_{n+1})}({\mathbf u,\mathbf s_{n+1}}).\]
Therefore, 
\begin{multline*}\varphi^{\vee}_{\rho}\otimes \id\left(\bigcap_{\bar r'\prec \bar r}\ker \varphi^{\vee}_{\bar r'}\right) \supset \\\prod_{1 \leq i<j\leq m}\prod_{1 \leq l \leq \rho_i,1 \leq l' \leq \rho_j}\left(u_{\bar\alpha_i}^{(l)}-u_{\bar\alpha_j}^{(l')} \right)^{\gamma(\bar\alpha^i,\bar\alpha^j)}
\Lambda_{(\rho,\alpha_{n+1})}({\mathbf u,\mathbf s_{n+1}})\cap  \bigcap_{\bar r' \prec^\zeta \bar r, \pr_{n+1}(\bar r')=\bar\rho}\ker (\tilde\psi_{\bar r'})\supset\\
\prod_{1 \leq i<j\leq m}\prod_{1 \leq l \leq \rho_i,1 \leq l' \leq \rho_j}\left(u_{\bar\alpha_i}^{(l)}-u_{\bar\alpha_j}^{(l')} \right)^{\gamma(\bar\alpha^i,\bar\alpha^j)}
\left(\bigcap_{\bar r' \prec^\zeta \bar r, \pr_{n+1}(\bar r')=\bar\rho}\ker (\tilde\psi_{\bar r'})\right).
\end{multline*}
Hence, we have:
\begin{multline}\label{twofactorshigherdim}
\varphi^{\vee}_{\bar r}\left( \bigcap_{r'\prec^{\zeta} r}\ker \varphi^{\vee}_{\bar r'} \right)
\supset\\
\tilde\psi_{\bar r}\left(\prod_{1 \leq i<j\leq m}\prod_{1 \leq l \leq \rho_i,1 \leq l' \leq \rho_j}\left(u_{\bar\alpha_i}^{(l)}-u_{\bar\alpha_j}^{(l')} \right)^{\gamma(\bar\alpha^i,\bar\alpha^j)}\right)
\tilde\psi_{\bar r}\left( \bigcap_{\bar r' \prec^\zeta \bar r, \pr_{n+1}(\bar r')=\bar\rho}\ker (\tilde\psi_{\bar r'}) \right).
\end{multline}

Let us study the two factors above separately.
By definition of $\tilde \psi_{\bar r}$ we have:
\begin{multline*}
\tilde\psi_{\bar r}\left(\prod_{1 \leq i<j\leq m}\prod_{1 \leq l \leq \rho_i,1 \leq l' \leq \rho_j}\left(u_{\bar\alpha_i}^{(l)}-u_{\bar\alpha_j}^{(l')} \right)^{\gamma(\bar\alpha^i,\bar\alpha^j)}\right)\\=
\prod_{\substack{(\bar \alpha^i,\alpha^i_{n+1}),(\bar \alpha^j,\alpha^j_{n+1}) \in P^{\zeta}\\\text{where }i<j}}
\prod_{\substack{1 \leq l\leq r_{(\bar\alpha^i,\alpha^i_{n+1})},\\1 \leq l'\leq r_{(\bar \alpha^j,\alpha^j_{n+1})}}}
\left( t_{(\bar \alpha^i,\alpha^i_{n+1})}^{(l)}-t_{(\bar \alpha^j,\alpha^j_{n+1})}^{(l')} \right)^{\gamma(\bar\alpha^i,\bar\alpha^j)}.
\end{multline*}

Define the following map:
\[g_{\bar \rho}:\Lambda_{(\bar \rho,a_{n+1})}({\mathbf u, \mathbf s_{n+1}})\rightarrow \Lambda_{(\bar \rho,a_{n+1}-\sum_{\alpha_i \in P} \rho_i(\zeta_{\max}-\zeta_i))}({\mathbf u,\mathbf s_{n+1}}),\]
\[s_{n+1}^{(\sum_{i'=1}^{i-1} \rho_{i'}(\zeta_{\max}-\zeta_{i'})+(l-1)(\zeta_{\max}-\zeta_{i})+1)}, \dots, s_{n+1}^{(\sum_{i'=1}^{i-1} \rho_{i'}(\zeta_{\max}-\zeta_{i'})+l(\zeta_{\max}-\zeta_{i}))},u_i^{(l)}\mapsto u_i^{(l)},\]
\[s_{n+1}^{(l)}\mapsto s_{n+1}^{(l-\sum_{i=0}^m \rho_i(\zeta_{\max}-\zeta_i))},~l \geq \sum_{\alpha_i 
\in P} \rho_i(\zeta_{\max}-\zeta_i).\]
Here $\zeta_{\max}-\zeta_i$ variables of the form $s_{n+1}^\bullet$ are mapped in to every variable of the form $u_i^\bullet$ while the remaining $s_{n+1}^\bullet$ are mapped into (distinct) $s_{n+1}^\bullet$.
We have:
\[\tilde \psi_{\bar r}=\psi_{d(\bar r)} \circ g_{\bar \rho},\]
where $d(\bar r)=(d(\bar r)_{(\bar\alpha^i,a)})$, $i=0,\dots,m$, $a=0,\dots,\zeta_i$,
$d(\bar r)_{(\bar\alpha^i,a)}=\bar r_{(\bar \alpha^i,a+\zeta_{\max}-\zeta_i)}$ and $\psi_{d(\bar r)} $ is defined in the following way (it differs from the $\psi_{d(\bar r)}$ defined in Subsection \ref{TechnicalLemma} only by a change of the subscripts):
\begin{align*}
u_{\bar\alpha^i}^{(j)}&\mapsto t_{(\bar\alpha^i,\zeta_{\max}-\zeta_i)}^{(j)}, j=1, \dots, r_{\bar\alpha^i,(\zeta_{\max}-\zeta_i)};\\
u_{\bar \alpha^i}^{(r_{(\bar \alpha^i,\zeta_{\max}-\zeta_i)}+j)},s_{n+1}^{(j)}&\mapsto t_{(\bar \alpha^i,\zeta_{\max}-\zeta_i+1)}^{(j)}, j=1, \dots, r_{(\bar \alpha^i,\zeta_{\max}-\zeta_i+1)};\\
u_{i}^{(r_{(\bar \alpha^i,\zeta_{\max}-\zeta_i)}+r_{(\bar \alpha^i,\zeta_{\max}-\zeta_i+1)}+j)},&s_{n+1}^{(r_{(\bar \alpha^i,\zeta_{\max}-\zeta_i+1)}+2j-1)},s_{n+1}^{(r_{(\bar \alpha^i,\zeta_{\max}-\zeta_i+1)}+2j)}\\&\mapsto t_{(\bar \alpha^i,\zeta_{\max}-\zeta_i+2)}^{(j)}, j=1, \dots, r_{(\bar \alpha^i,\zeta_{\max}-\zeta_i+2)};\\&\dots
\end{align*}
Here distinct variables of the form $u_{\bar \alpha_i}^{\bullet}$ are mapped to distinct variables of the form $t_{(\bar \alpha_i,a)}^{\bullet}$ and $a-(\zeta_{\max}-\zeta_i)$ variables of the form $s_{n+1}^{\bullet}$ are mapped to each variable $t_{(\bar \alpha_i,a)}^{\bullet}$.
The map $g_{\bar\rho}$ is surjective. Thus, we have
\[\tilde\psi_{\bar r}\left( \bigcap_{\bar r' \prec \bar r, \pr_{n+1}(\bar r')=\bar\rho}\ker (\tilde\psi_{\bar r'}) \right)=
\psi_{d(\bar r)}\left( \bigcap_{\bar r' \prec \bar r, \pr_{n+1}(\bar r')=\bar\rho}\ker (\psi_{d(\bar r')}) \right).\]

Note that the order $\prec^\zeta$ on vectors $\bar r$ defines an order on the vectors $d(\bar r)$. This order is lexicographic with respect to the following order on coordinates: $(\bar \alpha^i,a) <_d(\bar \alpha^j,b)$ if and only if $a+ (\zeta_{\max}-\zeta_i)<b+ (\zeta_{\max}-\zeta_{j})$ or  $a+ (\zeta_{\max}-\zeta_i)=b+ (\zeta_{\max}-\zeta_{j})$ and $i<j$.
Thus, using Lemma \ref{Split} we have:
\begin{multline*}\psi_{d(\bar r)}\left( \bigcap_{\bar r' \prec \bar r, \pr_{n+1}(\bar r')=\bar\rho}\ker (\psi_{d(\bar r')}) \right)\supset \\
\prod_{\substack{(\bar \alpha^i,\alpha^i_{n+1})<(\bar \alpha^j,\alpha^j_{n+1})\\\in P^\zeta}}\prod_{\substack{1 \leq l\leq r_{(\bar\alpha^i,\alpha^i_{n+1})},\\1 \leq l'\leq r_{(\bar \alpha^j,\alpha^j_{n+1})}}}\left( t_{(\bar \alpha^i,\alpha^i_{n+1})}^{(l)}-t_{(\bar \alpha^j,\alpha^j_{n+1})}^{(l')} \right)^{\varkappa((\bar \alpha^i,\alpha^i_{n+1}),(\bar \alpha^j,\alpha^j_{n+1}))} \Lambda_{\bar r}({\bf t})
\end{multline*}
(see \eqref{kappahigherdim}).
We obtain:
\begin{multline*}\varphi_{\bar r}^{\vee}\left(\bigcap_{\bar r' \prec \bar r} \ker(\varphi_{\bar r'}^{\vee})\right)\supset\\ \prod_{\substack{(\bar \alpha^i,\alpha^i_{n+1})<(\bar \alpha^j,\alpha^j_{n+1})\\\in P^\zeta}}\prod_{\substack{1 \leq l\leq r_{(\bar\alpha^i,\alpha^i_{n+1})},\\1 \leq l'\leq r_{(\bar \alpha^j,\alpha^j_{n+1})}}}\left(  t_{(\bar \alpha^i,\alpha^i_{n+1})}^{(l)}-t_{(\bar \alpha^j,\alpha^j_{n+1})}^{(l')} \right)^{\gamma((\bar \alpha^i,\alpha^i_{n+1}),(\bar \alpha^j,\alpha^j_{n+1}))} \Lambda_{\bar r}({\bf t}).\end{multline*}
Therefore, the map \eqref{inclusionToIdeal} is bijective.
\end{proof}

\section{Examples}\label{Examples}

\subsection{Parallelepiped}\label{parallel}
Consider the polytope $P_{d_1,\dots,d_n}:=(x_1,\dots,x_n)\in \mathbb{R}^n,0 \leq x_i \leq d_i$, $d_1,\dots,d_n \in \mathbb{N}$.
Choose $d_{n+1}\in\bN$, define $\zeta(\bar \alpha)=d_{n+1}, \bar \alpha \in P_{d_1,\dots,d_n}$. Then we have:
\[P_{d_1,\dots,d_n,d_{n+1}}=P_{d_1,\dots,d_n}^\zeta.\]
By induction we can construct a strict cube generating data for this polytope (in this example $f_i(\bar \alpha,\bar \beta)=0$ for all $\bar \alpha, \bar \beta, i$). In particular, the action of $\mathcal{A}_L$ on the $L$th graded component of $J^{\infty}_{\red}(R(P_{d_1,\dots,d_n}))^*$ is free. Let us compute the graded dimension of this component.

We define an order on $P_{d_1,\dots,d_n}\cap\bZ^n$ by setting $\bar \alpha<_n\bar \beta$ if for the largest $k$ such that $\alpha_k \neq \beta_k$ we have $\alpha_k < \beta_k$. Note that this order corresponds to the order  \eqref{order}, i.e.\ $<_{n+1}=<_n^\zeta$. Let us compute $\gamma(\bar \alpha, \bar \beta)$ for $\bar \alpha, \bar \beta\in P_{d_1,\dots,d_n}$. 
First, by \eqref{kappahigherdim} for $0\le i\le n-1$ and $(\bar\alpha,\alpha_{i+1})\leq_{i+1}(\bar\beta,\beta_{i+1})\in P_{d_1,\dots,d_i}^{\zeta_i}$ where $\zeta_i\equiv d_{i+1}$ we have: 
\begin{equation}
    \varkappa((\bar\alpha,\alpha_{i+1}),(\bar\beta,\beta_{i+1}))=\begin{cases}
    \beta_{i+1}-\alpha_{i+1}, \text{if}~ \bar\alpha \leq_i \bar\beta, \\
    \bar\beta_{i+1}-\bar\alpha_{i+1}-1, \text{if}~  \bar\beta<_i \bar\alpha
    \end{cases}
\end{equation}
Therefore, for $\bar \alpha,\bar \beta \in P_{d_1,\dots,d_n}$ we have
\[\gamma(\bar \alpha,\bar \beta)=\sum_{i=1}^n|\beta_i-\alpha_i|-S,\]
where $S$ is the number of $i$ for which the largest $j<i$ with $\alpha_j\neq\beta_j$ satisfies $(\alpha_i-\beta_i)(\alpha_j-\beta_j)<0$ (or no such $j$ exists). We note that if all $d_i=1$, then we recover Remark \ref{cube_rel}.

Note that the graded dimension of $\Lambda_{\bar r}(\mathbf{s})$ is equal to 
\[\prod_{\bar \alpha\in P_{d_1,\dots,d_n}}\prod_{l=1}^{r_{\bar \alpha}}\frac{1}{1-q^l}.\]
Using Proposition~\ref{DualSpaceDerivedQuagraticMon} we can now write the graded dimension of $J^\infty_{\red}R(P_{d_1,\dots,d_n})[L]$ as
\[\sum_{\bar r|\sum_{\bar \alpha \in P_{d_1,\dots,d_n}}r_{\bar\alpha}=L}v^{\sum_{\bar \alpha \in P_{d_1,\dots,d_n}} r_{\bar\alpha}\bar\alpha}\prod_{\bar \alpha\in P_{d_1,\dots,d_n}}\prod_{l=1}^{r_{\bar \alpha}}\frac{1}{1-q^l}\prod_{\bar \alpha<\bar \beta \in P_{d_1,\dots,d_n}}q^{r_{\bar \alpha}r_{\bar \beta}\gamma(\bar \alpha,\bar \beta)}.\]
Here variables $v_1,\dots,v_n$ correspond to the $n$-dimensional torus action on $R(P_{d_1,\dots,d_n})$ (and, subsequently, on $J^\infty_{\red}R(P_{d_1,\dots,d_n})$ while $q$ corresponds to the grading $\grad$.

\subsection{Simplex}
Consider the polytope 
\[P_{n,d}=\{(x_1,\dots,x_n)\in \mathbb{R}^n|d \geq x_n \geq \dots \geq x_1\geq 0\}\] for $d \in \mathbb{N}$. 
Consider $\zeta((x_1,\dots,x_n))=d-x_n$. Then $P_{n,d}^\zeta$ is  equal to the polytope $P_{n+1,d}$. Thus, by induction we can construct a strict cube generating data for all the polytopes $P_{n,d}$ using the construction from Example \ref{ConvexFunctionExample}. In particular, the action of $\mathcal{A}_L$ on the $L$th graded component of $J^{\infty}_{\red}(R(P_{n,d}))^*$ is free.

We define an order on $P_{n,d}\cap\bZ^n$ by setting $\bar \alpha<_n\bar \beta$ if for the largest $k$ such that $\alpha_k \neq \beta_k$ we have $\alpha_k < \beta_k$. Note that this order corresponds to the order  \eqref{order}, i.e.\ $<_{n+1}=<^\zeta$. Let us compute $\gamma(\bar \alpha, \bar \beta)$ for $\bar \alpha, \bar \beta\in P_{n,d}$. 
First, by \eqref{kappahigherdim} for $1\le i\le n-1$ and $(\bar\alpha,\alpha_{i+1})\le_{i+1}(\bar\beta,\beta_{i+1})\in P_{i,d}$ we have: 
\begin{equation}
    \varkappa((\bar\alpha,\alpha_{i+1}),(\bar\beta,\beta_{i+1}))=\begin{cases}
    \beta_{i+1}-\beta_i,\text{if}~ \beta_i>\alpha_{i+1};\\
    \beta_{i+1}-\alpha_{i+1}, \text{if}~ \beta_i\leq \alpha_{i+1} ~\text{and}~ \bar \alpha \leq \bar \beta;\\
    \beta_{i+1}-\alpha_{i+1}-1, \text{if}~  \beta_i\leq \alpha_{i+1} ~\text{and}~ \bar \alpha>\bar \beta.
    \end{cases}
\end{equation}
Now, for $\bar\alpha,\bar\beta\in P_{n,d}$ and $i=1,\dots,n-1$ we define
\begin{equation}\label{kappasimplexcomputation}
    \varkappa_i(\bar\alpha,\bar\beta)=\begin{cases}
    \beta_{i+1}-\beta_i,\text{if}~ \beta_i>\alpha_{i+1};\\
    \beta_{i+1}-\alpha_{i+1}, \text{if}~ \beta_i\leq \alpha_{i+1} ~\text{and}~ (\alpha_1,\dots,\alpha_i) \leq (\beta_1,\dots, \beta_i);\\
    \beta_{i+1}-\alpha_{i+1}-1, \text{if}~  \beta_i\leq \alpha_{i+1} ~\text{and}~ (\alpha_1,\dots,\alpha_i)>(\beta_1,\dots, \beta_i),
    \end{cases}
\end{equation}
and $\varkappa_0(\bar\alpha,\bar\beta)=|\alpha_1-\beta_1|-1$, if $\alpha_1\neq\beta_1$ and $0$ otherwise. Then by induction we have:
\[\gamma(\bar \alpha,\bar \beta)=\sum_{i=0}^{n-1}\varkappa_i(\bar\alpha,\bar\beta).\]

Note that $R(P_{n,d})$ is naturally embedded into $R(P_{n,1})$ as a sum of homogeneous components and $R(P_{n,1})$ is a polynomial algebra in $n+1$ variables which is acted upon by $\msl_{n+1}(\bC)$. Consider the restriction of this action to $R(P_{n,d})$. By Lemmas \ref{currentaction} and~\ref{NilpotentDerivations} we have an action  of $\msl_{n+1}(\bC)[s]$ on $J^\infty_{\red}R(P_{n,d})$ and for any $L\in \mathbb{N}$ the graded component $J^\infty_{\red}R(Q)[L]$ is preserved by this action. Due to \cite{DF} we have that $J^\infty_{\red}R(Q)[L]$ is isomorphic to the global Demazure module $\mathbb{D}_{d,L\omega_1}$.
The graded dimension of $J^\infty_{\red}R(Q)[L]$ is equal to
\begin{equation} \label{graded dimension for simplex}
\sum_{\bar r|\sum_{\bar \alpha \in P_{n,d}}r_{\bar\alpha}=L}v^{\sum_{\bar \alpha \in P_{d_1,\dots,d_n}}r_{\bar\alpha}}\prod_{\bar \alpha\in P_{n,d}}\prod_{l=1}^{r_{\bar \alpha}}\frac{1}{1-q^l}\prod_{\bar \alpha<\bar \beta \in P_{n,d}}q^{r_{\bar \alpha}r_{\bar \beta}\gamma(\bar \alpha,\bar \beta)}
\end{equation}
where the $v_i$ correspond to coordinates in a maximal torus of $\msl_{n+1}(\bC)$.


\appendix
\section{Veronese--Segre embeddings}
\label{appendix}

In the appendix we deal with special kind of toric projective embeddings. We use an approach, completely different to the one in the main body of the paper, to study their (reduced) coordinate rings and compute their characters. This approach is based on the representation theory of the current algebras. It would be interesting to investigate the combinatorial identities, obtained by comparing the formulae for the graded characters of the form \eqref{graded dimension for simplex} with the ones of the form \eqref{graded dimension of veronese-segre}, \eqref{character of affine demazure}.

For a semisimple Lie algebra $\g$ we denote by $V_\la^\g$ and $D_{d, \la}^\g$ the irreducible $\g$-module of highest weight $\la$ and the $\g[t]$ affine Demazure module of level $d$ and highest weight $d \la$ (see \cite{FL}, \cite{CV} for details on affine Demazure modules).

Let $\g_1, \hdots, \g_m$ be simple Lie algebras. 
To the algebra $\g_1 \oplus \hdots \oplus \g_m$ one can associate the affine Kac-Moody Lie algebra 
$\widehat{\g_1} \oplus \hdots \oplus \widehat{\g_m}$.
Consider the $m$-tuple $(\La_1, \hdots \La_m)$, where $\La_i$ is an affine integrable $\widehat{\g_i}$-weight. Then $L(\La_1) \T \hdots \T L(\La_m)$ naturally is an irreducible 
$\widehat{\g_1} \oplus \hdots \oplus \widehat{\g_m}$-module, where $L(\La_i)$ is an irreducible integrable $\widehat{\g_i}$-module.

Let $D_{(d_1, \hdots, d_m), (\la_1, \hdots, \la_m)}^{\g_1 \oplus \hdots \oplus \g_m}\subset L(\La_1) \T \hdots \T L(\La_m)$
be an affine Demazure module, where 
$d_1, \hdots, d_k$ are positive integers and $\la_i$ is a dominant weight of a simple Lie algebra $\g_i$.  It is clear that:
\[
D_{(d_1, \hdots, d_m), (\la_1, \hdots, \la_m)}^{\g_1 \oplus \hdots \oplus \g_m} \simeq D_{d_1, \la_1}^{\g_1} \T \hdots \T D_{d_m, \la_m}^{\g_m}.
\]



In this subsection we deal with the arc space of the following composition of Veronese and Segre embeddings:
\begin{equation} \label{veronese-segre map}
\begin{aligned} 
\bP(\Bbbk^{n_1}) \times \hdots \times \bP(\Bbbk^{n_m}) \hookrightarrow \bP(S^{d_1} \Bbbk^{n_1}) \times \hdots \times \bP(S^{d_m} \Bbbk^{n_m}) \\
\hookrightarrow \bP((S^{d_1} \Bbbk^{n_1}) \otimes \hdots \otimes (S^{d_m} \Bbbk^{n_m})).
\end{aligned}
\end{equation}
The image of this map is a projective toric variety, corresponding to the product of simplices. 

Consider the semisimple group $SL_{n_1} \times \hdots \times SL_{n_m}$ and its flag variety $\mathcal{B} = (SL_{n_1}/B_{n_1}) \times \hdots \times (SL_{n_m}/B_{n_m})$.
We consider each space $\Bbbk^{n_i}$, appearing in \eqref{veronese-segre map} as the first fundamental $\msl_{n_i}$-module, so $\Bbbk^{n_i} \simeq V_{\om_1}^{\msl_{n_i}}$ and $S^{d_i} \Bbbk^{n_i} \simeq V_{d_i \om_1}^{\msl_{n_i}}$. 
Therefore, the image of the map $\eqref{veronese-segre map}$ is isomorphic to the image of the map
\[
\mathcal{B} \rightarrow \bP(V_{d_1\om_1}^{\msl_{n_1}} \otimes \hdots \otimes V_{d_m\om_1}^{\msl_{n_m}}).
\]

In order to study the arc space of this image, one should replace the flag variety $\mathcal{B}$ by the semi-infinite flag variety $\fQ$ of the semisimple group $SL_{n_1} \times \hdots \times SL_{n_m}$ (see \cite{FiMi}, \cite{BF1}, \cite{Kato}) and study the homogeneous coordinate ring of the map
\[
\fQ \rightarrow \bP((V_{d_1\om_1}^{\msl_{n_1}} \otimes \hdots \otimes V_{d_m\om_1}^{\msl_{n_m}})[[t]]).
\]

Using Corollary \ref{ArcHomogenousSpace} as well as \cite[Proposition 4.1]{DF} applied to the algebra  $\msl_{n_1} \oplus \hdots \oplus \msl_{n_m}$ and its weight $(d_1\om_1, \hdots, d_m\om_1)$ we can describe the coordinate ring as \[
\Gamma = \bigoplus_{\ell \geq 0} \bigl( D_{(d_1, \hdots, d_m), (\om_1, \hdots, \om_1)}^{\msl_{n_1} \oplus \hdots \oplus \msl_{n_m}}[t]^{\odot \ell} \bigr)^*
\]
(we use that $V_{d_i \om_1}^{\msl_{n_i}} \simeq D_{d_i, \om_i}^{\msl_{n_i}}$). Here for a Lie algebra $\mathfrak{a}$ and its cyclic module $W$ with a cyclic vector $w$  we denote by $W^{\odot l}$ the Cartan component ${\rm U}(\mathfrak{a}).w^{\T l}\subset W^{\T l}$.  
We are going to compute the character of the $l$-th homogeneous component $\Gamma_l\subset \Gamma$.

The modules of the form $\Ga_\ell = D_{(d_1, \hdots, d_m), (\om_1, \hdots, \om_1)}^{\msl_{n_1} \oplus \hdots \oplus \msl_{n_m}}[t]^{\odot \ell}$ were studied in \cite{DF, DFF}. In particular, it was shown that $\Ga_\ell$ admits an action of the algebra $\cA_\ell$, which in this case is isomorphic to the algebra of symmetric polynomials in $\ell$ variables. It was also shown that the fiber of $\Ga_\ell$ with respect to this algebra at a generic point $\bc = (c_1, \hdots, c_\ell)$ is isomorphic to
\begin{equation}\label{genericfiber}
\Ga_\ell \T_{\cA_\ell} \bC_\bc \simeq \bigotimes_{i = 1}^\ell D_{(d_1, \hdots, d_m), (\om_1, \hdots, \om_1)}^{\msl_{n_1} \oplus \hdots \oplus \msl_{n_m}} (c_i),
\end{equation}
(where $\bC_\bc$ is the quotient of $\cA_\ell$ by the maximal ideal, corresponding to $\bc$).
It was also proved that to show that $\Ga_\ell$ is free over $\cA_\ell$ it is sufficient to obtain the surjection
\[
\Ga_\ell \T_{\cA_\ell} \bC_0 \twoheadleftarrow D_{(d_1, \hdots, d_m), (\ell \om_1, \hdots, \ell \om_1)}^{\msl_{n_1} \oplus \hdots \oplus \msl_{n_m}}
\]
where the latter has the same dimension as the right-hand side of~\eqref{genericfiber}. But this was shown for restriction to any $\msl_{n_i}$ in \cite[Proposition 3.2]{DF}. The surjectivity of the whole map follows. Therefore, we proved the following:
\begin{lem}
Module $\Ga_\ell$ is free over the algebra of symmetric polynomials in $\ell$ variables $\cA_\ell$. The fiber at 0 with respect to this algebra is isomorphic to
\[
D_{(d_1, \hdots, d_m), (\ell \om_1, \hdots, \ell \om_1)}^{\msl_{n_1} \oplus \hdots \oplus \msl_{n_m}}.
\]
\end{lem}

Thus, we obtain the desired
\begin{cor}
The character of the $\ell$-th homogeneous component of the arc space of the projective embedding \eqref{veronese-segre map} is equal to:
\begin{equation} \label{graded dimension of veronese-segre}
\ch \cA_\ell \ch D_{(d_1, \hdots, d_m), (\ell \om_1, \hdots, \ell \om_1)}^{\msl_{n_1} \oplus \hdots \oplus \msl_{n_m}} = \frac{1}{(q)_\ell} \prod_{i = 1}^{m} \ch D_{d_i, \ell \om_1}^{\msl_{n_i}}.
\end{equation}
\end{cor}
The last factor in this expression was found in \cite[Theorem 2.11]{FJKLM}. Namely, in the notations of \cite{FJKLM}, $D_{d_i, l\om_1}^{\msl_{n_i}} = \underbrace{V^{\msl_{n_i}}_{d_i \om_1} \ast \hdots \ast V^{\msl_{n_i}}_{d_i \om_1}}_{\ell} = V(\mathbf{n}_i, \mathbf d_i)$, where $\mathbf n_i = (\underbrace{n_i, \hdots, n_i}_{\ell}), \mathbf d_i =( \underbrace{d_i, \hdots, d_i}_{\ell})$, and $*$ denotes the fusion product. 
Now \cite[Formula (2.52)]{FJKLM} gives that
\begin{equation} \label{character of affine demazure}
\ch D^{\msl_{n_i}}_{d_i, \ell \om_1} = \sum_{\substack{\la \in \bZ^{n_i} \\ |\la| = d_i \ell}} e^{\la_2 \om_1} \hdots e^{\la_{n_i} \om_{n_i - 1}} \tilde S_{\la, (\ell^{d_i})}(q),
\end{equation}
where $\la = (\la_1, \hdots, \la_{n_i})$, $e^{\om_i}$ stay for the $\msl_{n_i}$-weights and $\tilde S_{\la, \mu}$ is the $q$-supernomial coefficient, defined in \cite[Formula (2.1)]{Sch}.

\end{document}